\documentclass[12pt]{amsart}

\DeclareFontFamily{U}{wasy}{}
\DeclareFontShape{U}{wasy}{m}{n}{
     <-5.5> wasy5
  <5.5-6.5> wasy6
  <6.5-7.5> wasy7
  <7.5-8.5> wasy8
  <8.5-9.5> wasy9
     <9.5-> wasy10
}{}
\DeclareFontShape{U}{wasy}{b}{n}{
 <-10> ssub * wasy/m/n
 <10-> wasyb10
 }{}
\DeclareFontShape{U}{wasy}{bx}{n}{ <-> ssub * wasy/b/n}{}
\DeclareFontShape{U}{wasy}{m}{sl}{ <-> wasysl10 }{}
\DeclareFontShape{U}{wasy}{m}{it}{ <-> ssub * wasy/m/sl }{}

\usepackage{amsmath,amsthm,amscd,amsfonts,wasysym,amssymb,epic,eepic,bm,bbm,tikz-cd,mathtools,mathdots,multicol,enumitem}
\usepackage[pagebackref,colorlinks=true,linkcolor=blue,citecolor=blue]{hyperref}
\usepackage[noabbrev]{cleveref}
\usepackage{MnSymbol}

\allowdisplaybreaks
\newtheorem{thm}{Theorem}[section]
\newtheorem{cor}[thm]{Corollary}

\newtheorem{lem}[thm]{Lemma}

\newtheorem{prop}[thm]{Proposition}
\theoremstyle{remark}
\newtheorem*{rem}{Remark}
\newtheorem*{example}{Example}

\newcounter{remarkscounter}

\numberwithin{equation}{section}
\newcommand{\A}{\mathbb{A}}
\newcommand{\GL}{\mathrm{GL}}

\newcommand{\PGL}{\textrm{PGL}}
\newcommand{\FF}{\mathbb{F}}

\newcommand{\QQ}{\mathbb{Q}}

\newcommand{\lto}{\longrightarrow}

\newcommand{\OO}{\mathcal{O}}
\newcommand{\CC}{\mathbb{C}}
\newcommand{\RR}{\mathbb{R}}
\newcommand{\GG}{\mathbb{G}}
\newcommand{\ZZ}{\mathbb{Z}}

\newcommand{\quash}[1]{}
\renewcommand\mod[1]{\ (\mathop{\rm mod}#1)}
\theoremstyle{definition}

\newcommand{\trd}{\mathrm{trd}}
\newcommand{\nrd}{\mathrm{nrd}}
\newcommand{\cusp}{\mathrm{cusp}}
\newcommand{\res}{\mathrm{res}}
\newcommand{\cont}{\mathrm{cont}}
\newcommand{\sw}{\Phi^{\operatorname{sw}}}
\newcommand{\diff}{\Phi^{\operatorname{diff}}}
\newcommand{\im}{\mathrm{im}}
\newcommand{\rank}{\mathrm{rank}}
\newcommand{\coeff}{\upsilon}
\newcommand{\defeq}{:=}

\newcommand{\grad}{\nabla}
\newcommand{\vol}{\mathrm{vol}} 
\newcommand{\covol}{\mathrm{covol}}
\newenvironment{psmatrix}
  {\left(\begin{smallmatrix}}
  {\end{smallmatrix}\right)}

\renewcommand{\bar}{\overline}
\numberwithin{equation}{subsection}

\newcommand{\one}{\mathbf{1}}

\newcommand{\norm}[1]{\left\lVert#1\right\rVert}
\newcommand{\abs}[1]{\lvert#1\rvert}

\newcommand{\lcm}{\operatorname{lcm}}

\allowdisplaybreaks

\begin{document}

\title{A nonabelian circle method}

\author{Nuno Arala}
\address{Mathematics Institute, Zeeman Building, University of Warwick, Coventry CV4 7AL, United Kingdom}
\email{Nuno.Arala-Santos@warwick.ac.uk}

\author{Jayce R. Getz}
\address{Department of Mathematics\\
Duke University\\
Durham, NC 27708}
\email{jgetz@math.duke.edu}

\author{Jiaqi Hou}
\address{Department of Mathematics\\
University of Wisconsin-Madison\\
Madison, WI 53706}
\email{jhou39@wisc.edu}

\author{Chun-Hsien Hsu}
\address{Department of Mathematics\\
Duke University\\
Durham, NC 27708}
\email{simonhsu@math.duke.edu}

\author{Huajie Li}
\address{Department of Mathematics\\
Johns Hopkins University\\
Baltimore, MD 21218}
\email{hli213@jhu.edu}

\author{Victor Y. Wang}
\address{Department of Mathematics\\
IST Austria}
\email{victor.wang@ist.ac.at}

\begin{abstract}
We count integral quaternion zeros of
$\gamma_1^2 \pm \dots \pm \gamma_n^2$,
giving an asymptotic when $n\ge 9$,
and a likely near-optimal bound when $n=8$.
To do so, we introduce a new, nonabelian delta symbol method, which is of independent interest.
Our asymptotic at height $X$ takes the form $cX^{4n-8} + O(X^{3n+\varepsilon})$ for suitable $c \in \CC$ and any $\varepsilon>0.$ We construct special subvarieties implying that, in general, $3n+\varepsilon$ can be at best improved to $3n-2.$
\end{abstract}

\maketitle

\setcounter{tocdepth}{1}
{\small
\begin{multicols}{2}
\tableofcontents
\end{multicols}
}

\section{Introduction}

Let $R$ be a ring and let
\begin{align*}
    \delta_a:=\begin{cases} 1 &\textrm{ if }a =0, \\ 0
    &\textrm{ if }a \neq 0,\end{cases}
\end{align*}
be the characteristic function of $\{0\} \subset R.$  A $\delta$-symbol method amounts to an expression for $\delta$ that is ``analytically tractable.'' One uses it to estimate how often $a =0$ as $a$ varies over a family of arithmetically interesting sets.   

When $R=\ZZ$
a $\delta$-symbol method was introduced by Duke, Friedlander, and Iwaniec and refined by Heath-Brown
\cite{DFI,HBnewcircle}.
It yields a form of the circle method that has been used to establish
a quantitative form of the Hasse principle for quadratic equations in as few as $4$ variables,
or $3$ variables if the equation is homogeneous.
In \cite{BrowningVisheCircleMethod} the $\delta$-symbol method was generalized to the case $R=\OO_F,$ where $\OO_F$ is the ring of integers in a number field.

An adelic version of the $\delta$-symbol method was introduced in \cite{GetzQuad}.  In this case the ring $R$ is a number field $F.$
A similar method is possible for global function fields $R$, though for $R=\FF_q[t]$ a simpler substitute is available, as illustrated by \cite{browning2015rational}.

In all previous settings, the $\delta$-symbol was used to study arithmetic questions related to commutative rings (typically $\ZZ$).
In this paper we move into the nonabelian setting.
More specifically, we consider the case that $R=D$ where $D$ is a central division algebra over a number field $F.$
Prospects for further generalizations are discussed in \S~\ref{SEC:generalize-delta}.

Thus let $D$ be a central division algebra over a number field $F.$
Let 
\begin{align}
    \RR_{>0} \lto F_\infty^\times
\end{align}
be the diagonal embedding and let $A_{\GG_m}$ be its image.  We often identify $F_\infty^\times$ with the center of $D_{F_\infty}^\times$ and hence $A_{\GG_m}$ with a central subgroup of $D_{F_\infty}^\times.$  

Let $F\langle Y_1,\dots,Y_n \rangle$ be the associative algebra of noncommutative polynomials in $Y_1,\dots,Y_n$ with coefficients in $F$ and let $P \in F\langle Y_1,\dots,Y_n \rangle$.  Let $f \in \mathcal{S}(D_{\A_F}^{n}).$  Our motivation for the expression for the $\delta$-symbol obtained in this work is to study the asymptotic behavior of expressions of the form
\begin{align}
    \Sigma(X):=\sum_{\substack{\gamma \in D^{n} \\P(\gamma)=0}}f\left(\frac{\gamma}{X} \right)
\end{align}
for $X\in A_{\GG_m}.$
To our knowledge, equations in division algebras have so far mainly been studied
qualitatively (see e.g.~\cite{niven1946note,pollackwaring,bankswatson} and references within).
For instance, \cite{niven1946note} constructs a solution to  $\gamma_1^2+\gamma_2^2+\gamma_3^2 = t$
for each $t$ in the ring of \emph{Hurwitz quaternions}
$$
\mathbb{H} = \left\{\frac{a+bi+cj+dk}{2}: a\equiv b\equiv c\equiv d\bmod{2}\right\}
\quad\textnormal{(detailed in \cite[Chapter~11]{voight2021quaternion})},
$$
with the caveat that the constructed solutions $\gamma_1,\gamma_2,\gamma_3\in \mathbb{H}$ generally have coordinates that are larger than the classical Waring problem for squares would allow.

We focus on quantitative questions.
For $\gamma = w+xi+yj+zk\in \mathbb{H},$ we let
$$
\norm{\gamma} = \max(\abs{w},\abs{x},\abs{y},\abs{z}),
\quad \trd(\gamma) = 2w,
\quad \nrd(\gamma) = w^2+x^2+y^2+z^2.
$$
The simplest problem we study below is to estimate the number of solutions to
\begin{equation}
\label{pm}
\gamma_1^2 \pm \dots \pm \gamma_n^2 = 0
\end{equation}
in Hurwitz quaternions $\gamma_1,\dots,\gamma_n$ with $\norm{\gamma_1},\dots,\norm{\gamma_n}\le X$ as $X \to \infty$.
Let us briefly explain some new difficulties of this problem compared to the classical abelian case, i.e.~when $D=F$.

One difficulty is that the set of solutions of \eqref{pm} seems to have a less useful automorphism group than in the abelian case.   Even the diagonal action of $D^\times$ on quadratic forms in $n$ variables over $D$ given by
\begin{align}
    g.P( Y_1,\dots,Y_n):=P(Y_1g,\dots,Y_ng)
\end{align} for $g \in D^\times$
 is nontrivial (i.e. $g.P \neq Pg^2$).  Moreover, we know of no transitive group action on the set of solutions of \eqref{pm} when $D^\times$ is not abelian. 

Another difficulty stems from an excess of points on \eqref{pm} when $n$ is small.
In the abelian case, it is well known that \emph{Hypothesis~K}
\cite[(17.2)]{waringsurvey}
holds for squares in $F$.
This allows for clever use of H\"{o}lder's inequality,
as is illustrated in \cite[Theorem~4.1, (4.5), and Lemma~4.5]{browning2021cubic} in the setting of cubic norm forms.
However, taking $n=4$ in the following example shows that
\emph{Hypothesis~K$^\ast$} (coined by \cite{hooley_greaves_harman_huxley_1997},  \cite[\S 1.3.2]{Browning:Quant}) fails for squares in $\mathbb{H}$,
and therefore the stronger \emph{Hypothesis~K} also fails for squares in $\mathbb{H}$.
\begin{example}
By Cayley--Hamilton, $x^2 = \trd(x)x - \nrd(x)$ for $x\in \mathbb{H}$.
Therefore, if $n\ge 2$, then
any equation of the form $\gamma_1^2 \pm \dots \pm \gamma_n^2 = 0$, with not all signs positive,
has $\asymp X^{3n-2}$ solutions of height $\le X$
with $\trd(\gamma_1) = \dots = \trd(\gamma_n) = 0$.\footnote{Indeed, these linear equations cut out a single quadric in $3n$ variables.}
This justifies the assertion on special subvarieties at the end of the abstract.  If $n\le 5$ (resp.~$n=6$), this exceeds (resp.~matches) 
$$
X^{\#\textrm{variables}-\textrm{degree}.\#\textrm{equations}}=X^{4n-8},
$$
the size of the usual heuristic main term in the the circle method.
\end{example}

\noindent Thus, at least for $n \leq 5,$ we cannot expect
an asymptotic point count for \eqref{pm} of the shape $c X^{4n-8} + O_\epsilon(X^{2n+\epsilon})$ for all $\epsilon>0$, with an error term matching the square-root of the trivial bound $X^{4n}$.  
Briefly,  square-root cancellation fails.
This is in stark contrast to the abelian case, where square-root cancellation holds in a strong form by \cite{HBnewcircle,GetzQuad,tran_thesis}.

In the setting of the Hamiltonians our main theorem is the following:

\begin{thm}
\label{thm:simplest-case}
Fix $n\ge 8$ and $\coeff_1,\dots,\coeff_n \in \{\pm 1\}$.
For $f_\infty \in C_c^\infty(D_\infty)$ there is a constant $c(f_\infty\one_{\widehat{\mathbb{H}}^n})$ such that
\begin{equation}
\label{EQN:main}
\sum_{\substack{\gamma\in \mathbb{H}^n \\
\coeff_1\gamma_1^2 + \dots + \coeff_n\gamma_n^2 = 0}} f_\infty(\gamma/X)
= c(f_\infty\one_{\widehat{\mathbb{H}}^n}) X^{4n-8} + O_{f_\infty,\epsilon}(X^{3n+\epsilon}).
\end{equation}
\end{thm}

\noindent
Here $\widehat{\mathbb{H}} \defeq \prod_p \mathbb{H}_{\ZZ_p}$, where $p$ ranges over all finite primes.
We discuss the constant $c(f_\infty\one_{\widehat{\mathbb{H}}^n})$ in Theorem \ref{thm:constant} below.
The exponent $3n+\epsilon$ in \eqref{EQN:main} (and in \eqref{EQN:maingen} below) might improve to $3n-2+\epsilon$ for $n\ge 6$ if we further optimized our methods.
But given the length of the paper, we have chosen to optimize our error term only asymptotically in $n$.

In the discussion above we restricted to $\mathbb{H}$ merely for concreteness.  We actually work much more generally
with quaternion algebras $D$ over $\QQ$ that are nonsplit at $2$ and $\infty.$  
\begin{thm}
\label{thm:general-case}
Let $D$ be a quaternion $\QQ$-algebra that is nonsplit at a set of places $S \supseteq \{2,\infty\}$.
Fix a maximal order $\OO_D\subseteq D$.
Fix $n\ge 8$ and $\coeff_1,\dots,\coeff_n \in \{\pm 1\}$.
Let $P(\gamma) := \coeff_1\gamma_1^2 + \dots + \coeff_n\gamma_n^2$.
For $f_\infty \in C_c^\infty(D_\infty)$ we have
\begin{equation}
\label{EQN:maingen}
\sum_{\substack{\gamma\in \OO_D^n \\ P(\gamma) = 0}} f_\infty(\gamma/X)
= c(f_\infty\one_{\widehat{\OO}_D^n}) X^{4n-8} + O_{f_\infty,\epsilon}(X^{3n+\epsilon}).
\end{equation}
\end{thm}
\noindent It should be possible, just with some more technical effort, to treat arbitrary test functions $f \in C_c^\infty(D_{\A_\QQ}^n).$  From a classical perspective this amounts to adding congruence conditions.
Moreover one should be able to replace the assumption $S \supseteq \{2,\infty\}$ with the simpler assumption $S \ne \emptyset$ (which, by the Hasse--Minkowski theorem, is equivalent to the statement that $D$ is nonsplit).
With yet more technical effort one should even be able to work over an arbitrary number field.  In fact the only place where we use the fact that $F=\QQ$ is in the geometry of numbers estimates in \S \ref{sec:geo:num}; this is why the majority of the paper is written in the setting of a general number field.

Let $\psi:=\otimes_v \psi_v:\QQ \backslash \A_\QQ \to \CC^\times$ be an additive character that is unramified at all finite places.

\begin{thm} \label{thm:constant}
Under the assumptions of Theorem \ref{thm:general-case}, one has
\begin{align*}
    c(f_\infty\one_{\widehat{\OO}_D^n})=c(f_\infty)\prod_{p}c(\one_{\OO_{D_p}^n}),
    \end{align*}
where
\begin{align} \label{c}
\begin{split}
    c(f_\infty)&=\int_{D_\infty} \int_{D_\infty^n}f_\infty(Y)\psi_\infty\left(\langle P(Y),Z \rangle \right)dY dZ,\\
    c(\one_{\OO_{D_p}^n})&=\int_{D_p} \int_{\OO_{D_p}^n}\psi_p\left(\langle P(Y),Z \rangle \right)dY dZ\in \RR_{>0}.
\end{split}
\end{align}
Here the measure on $D_v$ is normalized so the Fourier inversion holds. Moreover, $f_\infty$ can be chosen so that $c(f_\infty)\neq 0$.
\end{thm}

\noindent  The series \eqref{c} is reminiscent of singular series appearing in the abelian circle method. It also equals
\begin{align*}
    \int_{U_p}^{\mathrm{reg}} \one_{\OO_{D_p}^n}(u) du
\end{align*}
defined in the sense of \cite[\S~7]{Getz:Hsu:Leslie}, where $U_p<D^n_p$ is the closed subset cut out by $P(Y)=0$ and $du$ is a properly normalized measure on $U_p$.

In the rest of this introduction we will refer to the supplementary document \cite{Appendix}.  None of the results in this paper require the results of \cite{Appendix}, but the latter document may be helpful for researchers interested in generalizations.

\subsection*{Further context}

We now embed our result into the literature.
Taking $\coeff_1=\dots=\coeff_n=1$ for simplicity, and writing $\gamma = w+xi+yj+zk$ with suitable integrality conditions, the equation $P(\gamma):=\gamma_1^2 + \dots + \gamma_n^2 = 0$ treated in Theorem \ref{thm:simplest-case} becomes
\begin{equation}
\label{EQN:explicit-system}
\sum_{1\le i\le n} (w_i^2 - x_i^2-y_i^2-z_i^2)= \sum_{1\le i\le n} w_ix_i = \sum_{1\le i\le n} w_iy_i = \sum_{1\le i\le n} w_iz_i = 0.
\end{equation}
A natural approach to this specific system of equations would be to fix $w$ and exploit the linearity of the last three equations in $x,y,z$,
which when combined with uniform results on quadrics,
might lead to a version of Theorem \ref{thm:simplest-case} with a worse error term as $n\to \infty$.
Since our main goal is to introduce a more general method for nonabelian equations,
with Theorem \ref{thm:simplest-case} being the simplest possible example illustrating the method,
we do not discuss this approach further.

There is much work on systems of $2$ or $3$ quadrics;
see e.g.~\cite{Munshi:Pairs:Quadrics,HBP,PSW,vishe,arala2023analytic},
and references within.
The best results on systems of $\ge 4$ quadrics that are applicable to our situation are due to Myerson \cite{myerson}, which for $\ge 4$ quadrics always supersedes Birch's earlier work \cite{birch}.
Let $Q_1,Q_2,Q_3,Q_4\in \QQ[w_i,x_i,y_i,z_i]_{1\le i\le n}$ be the four quadratic forms that are set equal to $0$ in \eqref{EQN:explicit-system}.
Let
\begin{equation*}
\varrho_\RR \defeq \min_{W\in \RR^4 - \{0\}}
\rank(W_1Q_1+W_2Q_2+W_3Q_3+W_4Q_4)
= \min_{W\in D_\infty - \{0\}} \rank(\trd(WP(\gamma))).
\end{equation*}
The form $Q_2 = \sum_{1\le i\le n} w_ix_i$ has rank $2n$, so $\varrho_\RR \le 2n$.
In Lemma~\ref{LEM:general-rank-bounds}(2), we will see that $\varrho_\RR \ge 2n$.
So $\varrho_\RR = 2n$.
The quantity $\varrho_\RR$ is closely related to singular loci discussed in \cite[Appendix~A]{Appendix}.
Myerson's work \cite[Theorem~1.2, (1.5)]{myerson} is applicable when
\begin{equation}
\varrho_{\RR} > 8\cdot 4,
\end{equation}
i.e.~$n\ge 17$.
In comparison, Theorem~\ref{thm:simplest-case} provides an asymptotic for $n\ge 9$.

As a last contextual comment, we discuss generalizations of  Theorems \ref{thm:general-case} and \ref{thm:constant} to arbitrary quadratic polynomials.   
Our present work should immediately generalize to polynomials $P(\gamma) = \gamma^t M \gamma$, for \emph{symmetric} matrices $M\in \GL_n(F)$. 
These polynomials are precisely the polynomials that can be diagonalized by an $F$-linear 
change of variables $\gamma \mapsto A\gamma$ with $A\in \GL_n(F)$.  Not all (noncommutative) quadratic polynomials on $D^n$ are of this form.   

There is some chance that our methods could be adapted to handle general $M\in \GL_n(F)$, or even $M\in \GL_n(D)$.
We have focused on diagonal $P$ for simplicity.

\subsection*{Proof strategy}
As previously indicated, we develop and apply a nonabelian delta method to give a convenient analytic expression isolating the $\gamma' \in \OO_D^n$ satisfying  $P(\gamma')=0.$
The delta method is initiated in \S~\ref{SEC:delta}.
In \S~\ref{SEC:poisson}, we then apply Poisson summation in $\gamma'.$
After Poisson summation,
we obtain a dual sum over $\gamma\in D^n$ in which
there are two flavors of error terms we must treat.  

The first flavor of error term is a new feature of our nonabelian setting.  After Poisson summation, we require a spectral expansion of a sum over $\delta \in D^\times$ to further refine the $\gamma=0$ term.  This leads naturally to a spectral expansion in terms of automorphic representations of $D^\times_{\A_\mathbb{Q}}.$
The relevant local spectral estimates are established in \S\S~\ref{sec:nonarch:sp} and~\ref{sec:arch:bounds}.
\quash{In some special cases, classical abelian Mellin inversion does suffice, as we briefly explain in Appendix~\ref{SEC:abelian}.  However, usually there are infinite-dimensional automorphic representations of $D_{\mathbb{A}_{\QQ}}^\times$ that are everywhere unramified, and whenever they exist, abelian Mellin inversion is not sufficient.}

The second flavor of error term comes, as in the abelian case, from $\gamma \neq 0.$  However the behavior of these terms is markedly different.  The reduced norm scales quadratically, i.e.~with degree $d:=\sqrt{\dim{D}}=2.$  Thus Poisson summation leads to longer dual sums over $\gamma$ than in the abelian case $d=1$.
Usual heuristics suggest that the dual sum (i.e.~the sum after applying Poisson summation) is of the same complexity as the original.
However, we discover new, nonabelian \emph{vanishing phenomena} in \S~\ref{SEC:local-exp-sums}, which allow for progress over $\gamma\ne 0$.
The situation is still complicated by the fact that our complete exponential sums fail square-root cancellation, asymptotically in $n$ \cite[Appendix~B]{Appendix}.
Nonetheless, by identifying various geometric sources of vanishing, cancellation, and sparsity, we eventually obtain
an error term that has an exponent $\frac34$ times the exponent of the trivial bound $X^{4n}$.\footnote{This ``slope'' $\frac34$ also appears in work on cubics, such as \cite{davenport2005analytic,heath1983cubic,vaughan1986waring,hooley1986Lfunctions,hooley_greaves_harman_huxley_1997,HB:circle:diagonal:cubic}.}
The tools involved are diverse, featuring
Cartan decomposition, matrix identities, iterative Gauss sum calculations, the geometry of quadric fibrations, and lattice-point methods (geometry of numbers).

Our work on $\gamma\ne 0$ appears in \S\S~\ref{SEC:general-hessian-based-analysis}--\ref{sec:geo:num}.
General Hessian-based analysis is done in \S~\ref{SEC:general-hessian-based-analysis}, using a rank estimate proven via the Cayley--Hamilton theorem.
Deeper local integral estimates are proven in \S~\ref{SEC:local-exp-sums}.
A final global application of the geometry of numbers appears in \S~\ref{sec:geo:num}, based on some fortuitous bounds on the successive minima of a relevant family of lattices.

We leave Kloosterman-type averaging as an interesting open question.
By this, we mean
the use of nontrivial averaging over the Fourier inversion variable
in our local integrals.
This would be analogous to Kloosterman's idea of averaging over numerators in the abelian case, as is done in \cite{HBnewcircle,GetzQuad}, explicitly or implicitly.
We expect the prime-case formulas in \cite[Appendix~B]{Appendix} to be useful for this purpose.
Alternatively, one could improve on Theorems~\ref{thm:simplest-case} and~\ref{thm:general-case} by sharpening our geometry of numbers estimates in \S~\ref{sec:geo:num} on average.
We also suspect that a secondary term of order $X^{3n-2}$ may exist, and may be detectable by methods of \cite{HBnewcircle,HB:circle:diagonal:cubic,GetzQuad,tran_thesis,wang2023_isolating_special_solutions}.

We also leave the case $d =\sqrt{\dim{D}} \ge 3$ open; it is unclear what bounds to expect.
Extrapolating from $d \in \{1,2\}$, a tentative guess is that after accounting for nonabelian vanishing phenomena,
we might morally have sums of length $\ll (X^{d^2(d-1)})^{(d-1)/(d^2-d)}$
and local cancellations of quality $\gg (X^d)^{d-1} (X^d)^{1/2}$,
all raised to the $n$th power,
for a total saving of roughly $X^{dn/2}$ over the trivial bound $X^{d^2n}$.
If so, then we might hope for an asymptotic point count with main term $\asymp X^{(n-2)d^2}$ when $n>4d$, even without Kloosterman-type averaging.
Yet by the methods of \cite[Appendix~A]{Appendix}, we have $\sigma_{\RR} \ge ((d-1)^2+1)n$ for $d\ge 2$, at least if $D_\infty$ is split (which is always the case if $2\nmid d$).
If $\sigma_{\RR} = (d^2-2d+c)n$, where $c\ge 2$, then \cite{myerson} requires $n>\frac{8d^2}{2d-c}\ge 4d+4+\frac{4}{d-1}$,
and \cite{birch} requires strictly more.

Finally, we mention that our nonabelian delta method
might have a natural interpretation in terms of ``nonabelian Dirichlet arcs'' like $q\theta - a \ll 1/Q$, for various elements $a,q\in D$ and parameters $Q$.
It is also possible that ``nonabelian exponential sums'' on these arcs could be bounded without using Poisson summation as we do.
It may also be worth comparing our setting to the \emph{orbital circle method} and \emph{orbital exponential sums} studied by Bourgain, Kontorovich, et al.~\cite{kontorovich2014orbital}.
However, a key difference is that the equations discussed in \cite{kontorovich2014orbital} are valued in $F$, whereas our equations are valued in $D$.

These are all interesting questions, which we leave open.

\subsection*{Acknowledgements}

This project is the outcome of the 2023-2024 Duke Research Scholars program, organized by J.~R.~G. and funded by the RTG DMS-2231514.

N.~A. thanks Simon Rydin Myerson for helpful conversations, and was funded through the Engineering and Physical Sciences Research Council Doctoral Training Partnership at the University of Warwick. Part of this work was performed while N.~A. was in residence at the Mittag-Leffler Institute in 2024, the hospitality and financial support of which are gratefully acknowledged.

J.~R.~G.~ is partially supported by DMS--2400550.  He thanks  Oscar Marmon and Lillian Pierce for inviting him to the Hausdorff School on the circle method in 2021. This is where he first announced Lemma \ref{lem:delta}, which ultimately led to the present collaboration.  He also acknowledges John Voight's help with notational conventions. 

J.~H. is partially supported by DMS--1902173.

V.~Y.~W. thanks Tim Browning, Jakob Glas, and Damaris Schindler
for encouragement and for interesting discussions on closely related topics,
and was supported in part by the European Union's Horizon~2020 research and innovation program under the Marie Sk\l{}odowska-Curie Grant Agreement No.~101034413,
and briefly by the Stanford Math Department.
Thanks are also due to Tim Browning for several helpful comments.

\section{Conventions}
\label{SEC:conventions}

\subsection{Local fields}

Let $F$ be a global field and let $v$ be a place of $F.$
We let $|\cdot|_v$ be the number theorist's norm on $F_v$.
If $F_v \cong \RR$ then $|\cdot|_v$ is the absolute value, and if $F \cong \CC$ then $|z|_v=z\bar{z}$.
If $F_v$ is non-Archimedean, let $\OO_{F_v}$ be its ring of integers and $\varpi_v$ be a uniformizer. Then $|\varpi^{-1}_v|_v$ is the cardinality $q_v$ of the residue field $\OO_{F_v}/\varpi_v\OO_{F_v}$. The valuation associated to the place $v$ is denoted by the same symbol: $v:F_v\to \ZZ\cup \{\infty\}.$  Thus $v(\varpi_v)=1$.

The usual norm $z \mapsto (z\overline{z})^{1/2}$ on $\CC$ is denoted by $|\cdot|.$ 
Fix a place $v$ of $F$ and drop it from notation, writing $F:=F_v,$ $|\cdot|:=|\cdot|_v,$ etc.  This creates the possibility of confusion when $F \cong \CC.$  To alleviate this, we do not identify $F$ and $\CC$ and we take the standard convention that $|x|=x\overline{x}$ if $x \in F$ and $|y|=(y\overline{y})^{1/2}$ if $y \in \CC.$  In other words, the meaning of $|\cdot|$ changes depending on whether we use the symbol $F$ or $\CC$ for the domain.    


If $F$ is a finite extension of $\QQ_p$,
the \emph{standard additive character} $\psi_F: F \to \CC^\times$ 
is defined by the formula $x\mapsto e^{2\pi i(\operatorname{Tr}_{F/\QQ_p}(x)\bmod{\ZZ_p})}$, where $\operatorname{Tr}_{F/\QQ_p}(x)\bmod{\ZZ_p}$ denotes the image of $\operatorname{Tr}_{F/\QQ_p}(x)$ in $\QQ_p/\ZZ_p \cong \ZZ[1/p]/\ZZ \subset \RR/\ZZ$.  Here the first isomorphism is induced by the inclusion of $\ZZ[1/p]$ into $\QQ_p.$  If $F$ is Archimedean, then the standard additive character is defined by $x\mapsto e^{-2\pi i\operatorname{Tr}_{F/\RR}(x)}$. In general, a nontrivial additive character on $F$ is given by $x\mapsto \psi_{F}(\alpha x)$ for some $\alpha\in F^\times$. 

Let $F$ be a finite extension of $\QQ_p$. For a nontrivial additive character $\psi$ on $F$, its \textit{conductor} $c(\psi)$ is the integer such that $\psi$ is trivial on $\varpi^{c(\psi)}\OO_F$ but not on $\varpi^{c(\psi)-1}\OO_F$. We say $\psi$ is \textit{unramified} if $c(\psi)=0$.

\subsection{Central simple algebras}
For a central simple algebra $D$ over a field $F,$ let $d:=\sqrt{\dim_F D}$.
We let $\mathrm{trd}:D \to F$ and $\mathrm{nrd}:D \to F$ be the reduced trace and the reduced norm.
When $d=2$,
we let $\dagger$ denote the unique standard involution on $D$ \cite[Definition~3.2.4 and Corollary~3.4.4]{voight2021quaternion}.
Then $\trd(A) = A+A^\dagger$ and $\nrd(A) = A^\dagger A = A A^\dagger$.
When $F$ is a number field, we use the same notation for the adelic and local versions of these maps.
At all split places $v$, we have $D_v:=D \otimes_FF_v \cong M_d(F_v)$ by definition.
By the Skolem-Noether theorem, for any such isomorphism 
$\mathrm{trd}$ is sent to $\mathrm{tr}$ and $\mathrm{nrd}$ is sent to $\det.$
Moreover, if $d=2$, then $\dagger$ is sent to
the map sending an element of $M_2(F_v)$ to its adjugate,
because the adjugate is the unique standard involution on $M_2(F_v)$ \cite[Example~3.2.8 and Corollary~3.4.4]{voight2021quaternion}.
Any maximal order of $M_d(F_v)$ over $\OO_{F_v}$ is isomorphic to $M_d(\OO_{F_v})$ by an inner automorphism of $M_d(F_v)$ \cite[Chapter~10]{voight2021quaternion}.

Let $F$ be a local field.
For any integer $m\ge 0$, the \textit{box norm} on $F^m$ is given by
\begin{align}
\label{box-vector-norm}
\norm{(v_1,\ldots, v_m)}:=\max_{1\le i\le m} |v_i|.
\end{align}
Often we will make an identification of $D$ with $F^{d^2}$ as $F$-vector spaces and let $\norm{\cdot}$ be the induced box norm on $D$.
We sometimes call $\norm{x}$ the \emph{magnitude} of $x$, to avoid linguistic confusion with the reduced norm $\nrd(x)$ for $x\in D$.

\begin{rem}
If $F\cong \CC$, then $\abs{\cdot}$ and $\norm{\cdot}$
are not norms in the traditional sense
(they do not satisfy the triangle inequality).
\end{rem}

\subsection{Measures and Fourier transforms}
Let $F$ be a local field. If $F$ is non-Archimedean, the Schwarz space $\mathcal{S}(D)=C_c^\infty(D)$ is the space of compactly supported smooth functions.  If $F$ is Archimedean $\mathcal{S}(D)$ is the usual Schwartz space on $D$ (viewed as an $F$-vector space).  Fix a nontrivial additive character $\psi:F\to \CC^\times$. For a central simple $F$-algebra $D$, let 
\begin{align*}
    \mathcal{F}_D(f)(y):=\int_{D} f(x)\psi(\trd(xy))\,dx
\end{align*}
be the Fourier transform on the Schwartz space $\mathcal{S}(D).$ Here the measure $dx$ on $D$ is always normalized to be self-dual with respect to $\psi.$  In other words, it is normalized so that the Fourier inversion holds:
\begin{align}\label{fourier-inversion}
    f(y)=\int_{D} \mathcal{F}_D(f)(x)\psi(-\trd(xy))\,dx.
\end{align}
We will often drop the subscript and write $\mathcal{F}=\mathcal{F}_D$ whenever the context is clear. Any Haar measure on $D^\times$ is a $\RR_{>0}$ multiple of $|\nrd|^{- d} dx.$  We will specify the constant when necessary.

Suppose $F$ is non-Archimedean. Let $\OO_D$ be a maximal order of $D$ over $\OO_F$ and 
\begin{align*}
    \OO_D^\#:=\big\{ y\in D: \trd(xy)\in \OO_F\, \forall x\in \OO_D\big\}
\end{align*}
be its dual. Note that $\OO_D^{\#\#}=\OO_D$ and 
\begin{align*}
    \mathcal{F}(\one_{\OO_D})=dx(\OO_D)\one_{\varpi^{c(\psi)}\OO_D^\#}.
\end{align*}
If $D$ is split and $\psi$ is unramified, then our normalizations ensure that the function $\one_{\OO_{D}}$ is self-dual, i.e., $\mathcal{F}(\one_{\OO_{D}})=\one_{\OO_{D}}$ and $\vol(\OO_{D})=1$.


\subsection{Gauss sum}
Generalized Gauss sums will be repeatedly used in our local estimates throughout the paper.
Assume $F$ is non-Archimedean
and define for $(a,t,\xi)\in F\times F^\times \times F$, the \emph{generalized Gauss sum}
\begin{align*}
\mathcal{G}(a,t,\xi):=\int_{\OO_F} \psi\left(\frac{ay^2+\xi y}{t}\right) \,dy.
\end{align*}
We emphasize that our Gauss sums are normalized so that the trivial bound is $\abs{\mathcal{G}} \le \vol(\OO_F)$; in classical settings, Gauss sums are often normalized differently.

\begin{lem}\label{LEM:classical-gauss-sum-properties}
Suppose $F$ is a finite extension of $\QQ_p$
and $\psi$ is unramified.
In general
there is a constant $c_F>0$, with $c_F=1$ if $p\ne 2$,
such that
\begin{align}\label{eq:gaussbd}
c_F\, |\mathcal{G}(a,t,\xi)|_\infty
\leq \min(1,|t/a|)^{1/2}
= q^{\min(v(a)-v(t),0)/2},
\end{align}
which is an equality if
$v(\xi) \ge \min(v(a),v(t))$
and $p\neq 2$.
If $|t/a|\ge 1$, i.e.~$v(a) \ge v(t)$, then
\begin{align*}
\mathcal{G}(a,t,\xi)
= \mathcal{G}(t,t,\xi)
= \mathcal{G}(0,t,\xi)
= \one_{\OO_F}(\xi/t).
\end{align*}
Suppose $v(a) \le v(t)$, i.e.~$|t/a|\le 1$.
Then $\mathcal{G}(a,t,\xi)=0$ unless
\begin{equation}
\label{gauss-sum-non-vanishing-condition}
\xi/a\in \OO_F,
\end{equation}
i.e.~$|\xi/a|\le 1$.
In addition, if $\xi/a\in 2\OO_F$, then
\begin{align}
\label{evaluate-off-center-gauss-sum}
\mathcal{G}(a,t,\xi)=\psi\left(-\frac{\xi^2}{4at}\right)\mathcal{G}(a,t,0).
\end{align}
\end{lem}

\begin{proof}
Everything is clear if $v(a)\ge v(t)$, so we may assume $v(a) \le v(t)$.
Dividing $(a,t,\xi)$ by $a\ne 0$, we may assume $a=1$.
All but the bound for $p=2$ follow from \cite[Lemma 6.1]{GetzQuad} and its proof.
For $p=2$, it follows from the Van der Corput Lemma \cite[Proposition 3.3]{Cluckers:vdC} 
\end{proof}

\subsection{Adelic measure and Fourier transform }

Let $D$ be a central simple algebra over a number field $F$. Fix a maximal order $\OO_D \subset D$ over $\OO_F$.
At every finite place $v$, the localization $\OO_{D_v} \subset D_{F_v}$ is a maximal order over $\OO_{F_v}$ by \cite[Lemma~10.4.3]{voight2021quaternion}.

The Schwartz space of $D_{\A_F}$ is defined as
\begin{align*}
  \mathcal{S}(D_{\A_F}):=  \mathcal{S}(D_{F_\infty}) \otimes \mathcal{S}(D_{\A_F^\infty}).
\end{align*}
Here
\begin{align*}
    \mathcal{S}(D_{\A_F^\infty}):=\sideset{}{'}\bigotimes_{v\nmid\infty} \mathcal{S}(D_{F_v})=\sideset{}{'}\bigotimes_{v\nmid\infty} C^\infty_c(D_{F_v})
\end{align*}
where the restricted tensor product is taken with respect to the basic functions $\one_{\OO_{D_v}}$, and
\begin{align*}
    \mathcal{S}(D_{F_\infty}):=\widehat{\bigotimes}_{v|\infty}\mathcal{S}(D_{F_v})
\end{align*}
where $\mathcal{S}(D_{F_\infty})$ and $\mathcal{S}(D_{F_v})$ are the usual Schwartz spaces on the $\RR$-vector spaces $D_{F_\infty}$ (resp.~$D_{F_v}$) with the natural Fr\'echet topology, and $\widehat{\bigotimes}$ is the completion of the (projective) tensor product.
Concretely, this means that any element of $\mathcal{S}(D_{\A_F})$ is a finite sum of functions of the form $f_\infty \otimes \left(\otimes_{v\nmid \infty} f_v\right)$ with $f_\infty \in \mathcal{S}(D_{F_\infty})$ and $f_v \in \mathcal{S}(D_{F_v})$ where $f_v=\one_{\OO_{D_v}}$ for all but finitely many $v.$  

Let $\psi:F \backslash \A_F \to \CC^\times$ be a nontrivial additive character. Then $\psi=\otimes_v \psi_v$ and we have an adelic Fourier transform $\mathcal{F}_{D_{\mathbb{A}_F}}:=\otimes_v \mathcal{F}_{D_v}$. Explicitly, it is given by
\begin{align*}
\mathcal{F}_{D}:\mathcal{S}(D_{\A_F}) & \tilde{\lto} \mathcal{S}(D_{\A_F})\\
 f&\longmapsto \left(y \mapsto   \int_{D_{\A_F}}f(x)\psi\left(\mathrm{trd}(yx) \right)dx \right)
\end{align*}
where $dx=\prod_v dx_v$. 

On $\mathcal{S}(D_{\mathbb{A}_F}^m)$ we let $\mathcal{F}_i$ be the Fourier transform in the $i$th entry. For instance if $\Phi\in \mathcal{S}(D_{\mathbb{A}_F}^2)$, then
\begin{align*}
    \mathcal{F}_{2}(\Phi)(x,y)=\int_{D_{\A_F}} \Phi(x,z)\psi\left(\mathrm{trd}(yz)\right) dz.
\end{align*}
is the Fourier transform in the second variable. We use the same notation over local fields.

\subsection{Analytic number theory conventions}  We use the phrase ``dominated by'' as a synonym for ``is bounded by a constant times.''

\section{The \texorpdfstring{$\delta$}{delta}-symbol method}
\label{SEC:delta}

Let $D$ be a central division algebra over a number field $F$.
For $\alpha \in D$ let 
\begin{align*}
\delta_{\alpha}:=\one_{\alpha=0}:=\begin{cases} 1 & \textrm{ if } \alpha=0,\\
    0 & \textrm{ otherwise.} \end{cases}
\end{align*}
The following is a nonabelian expansion of the $\delta$-symbol:
\begin{lem}
\label{lem:delta}
Let $\Phi\in \mathcal{S}(D^2_{\mathbb{A}_F})$. Assume $\Phi(t,0)=0$ for all $t\in D_{\mathbb{A}_F}$
and $\mathcal{F}_2(\Phi)(0,0)\neq 0$.
For all $Q\in A_{\GG_m}$,
\begin{equation}
\label{EQN:delta-nitpick}
b_{\Phi,Q}\, \delta_{\alpha}
= \frac{1}{|Q|^{\dim _FD}}
\sum_{\delta \in D^\times}\left(\Phi\left(\frac{\alpha\delta^{-1}}{Q},\frac{\delta}{Q}\right)
- \Phi\left(\frac{\delta}{Q},\frac{\delta^{-1}\alpha}{Q}\right)\right),
\end{equation}
where
\begin{equation*}
b_{\Phi,Q} \defeq \sum_{\delta \in D} \mathcal{F}_2(\Phi)(0,Q\delta)
= \mathcal{F}_2(\Phi)(0,0)+O_{N,\Phi}(|Q|^{-N})
\end{equation*}
for all $N \geq 0$.
Moreover, for $|Q|$ sufficiently large
\begin{equation}
\label{EQN:delta-main}
\delta_{\alpha}
= \frac{c_{\Phi,Q}}{|Q|^{\dim _FD}}
\sum_{\delta \in D^\times}\left(\Phi\left(\frac{\alpha\delta^{-1}}{Q},\frac{\delta}{Q}\right)
- \Phi\left(\frac{\delta}{Q},\frac{\delta^{-1}\alpha}{Q}\right)\right),
\end{equation}
where $c_{\Phi,Q} \defeq b_{\Phi,Q}^{-1} = \frac{1}{\mathcal{F}_2(\Phi)(0,0)}+O_{N,\Phi}(|Q|^{-N})$.
\end{lem}

\begin{proof}
If $\alpha=0$ then by our vanishing assumption $\Phi(t,0)=0$, we have
\begin{align*}
\sum_{\delta \in D^\times}
\left(\Phi\left(\frac{\alpha\delta^{-1}}{Q},\frac{\delta}{Q}\right)
- \Phi\left(\frac{\delta}{Q},\frac{\delta^{-1}\alpha}{Q}\right)\right)
= \sum_{\delta \in D^\times}\Phi\left(0,\frac{\delta}{Q}\right)
= \sum_{\delta \in D}\Phi\left(0,\frac{\delta}{Q}\right).
\end{align*}
By Poisson summation this is
\begin{align*}
|Q|^{\dim_FD} \sum_{\delta \in D} \mathcal{F}_2(\Phi)(0,Q\delta)
= |Q|^{\dim_FD}\mathcal{F}_2(\Phi)(0,0)+O_{N,\Phi}(|Q|^{-N}).
\end{align*}
This establishes the equality \eqref{EQN:delta-nitpick} in the lemma when $\alpha=0$.

If $\alpha \neq 0,$ by changing variables
$\delta \mapsto  \delta^{-1} \alpha$
we see that 
\begin{align*}
\sum_{\delta \in D^\times}\Phi\left(\frac{\alpha\delta^{-1}}{Q},\frac{\delta}{Q}\right)
=\sum_{\delta \in D^\times}\Phi\left(\frac{\delta}{Q},\frac{\delta^{-1}\alpha}{Q}\right)
\end{align*}
and hence the sum on the right in \eqref{EQN:delta-nitpick} vanishes.
Of course, \eqref{EQN:delta-main} is immediate from \eqref{EQN:delta-nitpick}.
\end{proof}


\begin{rem}
Lemma~\ref{lem:delta} can be shown to hold even if $D$ is an \emph{alternative} non-associative division algebra, e.g.~an octonion algebra.  This is because in such an algebra unique two-sided inverses of nonzero elements exist, and if $x$ and $y$ are both nonzero one has $(xy)^{-1}=y^{-1}x^{-1}$ and $y=x^{-1}(xy).$
\end{rem}


Since $D$ is a division algebra, for $\alpha \in D,$  $\alpha = 0$ holds if and only if $\nrd(\alpha) = 0$.  
If we replace $D$ by the split central simple algebra $M_n(F),$ then the idea behind Lemma~\ref{lem:delta} naturally produces an identity where the condition $\alpha = 0$ is replaced by the weaker condition $\det(\alpha)=0.$  This observation leads to the generalization we explore in the next subsection.

\subsection{A remark on generalizations}
\label{SEC:generalize-delta}

Let $G$ be a linear algebraic group over $F$ acting on a quasi-affine scheme $X$ over $F:$
\begin{align} \label{action:map}
X \times G \lto X.
\end{align}
Assume that $X$ admits a unique open $G$-orbit $O,$ that $b_0 \in O(F),$ and that the stabilizer of $b_0$ in $G$ is trivial.  Thus the action map induces a bijection $G(F) \tilde{\to} O(F).$
Finally assume that there are automorphisms (of schemes, not group schemes) $\iota:X \to X$ and $\iota:G \to G$ of order $2$ such that $\iota(gh)=\iota(h)\iota(g)$ and $\iota(b_0.g)=b_0.\iota(g).$

Let $\Phi:X(\A_F) \to \CC$ be a function and let $x \in X(F).$  Consider
\begin{align} \label{deltaPhi}
f_{\Phi}(x):=\sum_{g \in G(F)} \Phi\left( x.g^{-1}, b_0.\iota(g) \right)-\sum_{g \in G(F)}\Phi\left(b_0.g,\iota(x).\iota(g^{-1})\right).
\end{align}
Here we assume that $\Phi$ is sufficiently nice that both sums converge absolutely.  For example, if $X$ is smooth we can take $\Phi \in C_c^\infty(X(\A_F))=C_c^\infty(X(F_\infty)) \otimes C^{\infty}_c(X(\A_F^\infty)).$  

If $x \in O(F)$ we can choose an $h \in G(F)$ such that $b_0.h=x.$  We then change variables $g \mapsto g^{-1}h$ in the first sum in the definition of $f_{\Phi}(x)$ to see that $f_{\Phi}(x)=0$ for $x \in O(F).$  In other words, $f_{\Phi}(x)$ vanishes on the open orbit $O(F).$ We hasten to point out that $f_{\Phi}$ is not necessarily constant away from $O(F)$, unless $X(F)-O(F)$ consists of a single point.

Expressions of the type \eqref{deltaPhi} could be useful for counting points of schemes that do not lie in an open orbit under a group action.
The $\delta$-symbol in Lemma \ref{lem:delta} is essentially the special case where $X$ and $G$ are the scheme and group with points in an $F$-algebra $R$ given by $X(R):=D \otimes R,$ $G(R):=(D \otimes R)^\times.$  In this case $X(F)-O(F)=\{0\},$ we can take $\iota$ to be induced by 
the standard involution $\dagger$ on $D,$ and $b_0=I_D.$  In order to simplify the presentation we did not incorporate the involution into Lemma \ref{lem:delta}.

Another natural case is when $G=\GL_n$ and $X=M_{n \times n}$ ($n \times n$ matrices) equipped with the usual right action.  In this case $X(F)-O(F)=\{Z \in M_{n \times n}(F):\det (Z)=0\}$ and we can take $\iota$ to be the transpose.  In unpublished work Myerson and Vishe considered this setting when $n=2.$

One could also work more generally with reductive monoids, and make use of the Poisson summation formula for these objects conjectured to exist by Braverman and Kazhdan \cite{BK-lifting,NgoSums}.  The conjectural Poisson summation formula has been established in certain cases, see \cite{BK:basic:affine,BK:normalized,Choie:Getz,Getz:Hsu,Getz:Hsu:Leslie,Getz:Liu:BK,Getz:Liu:Triple}.

\section{Poisson summation}
\label{SEC:poisson}

Let $f \in C_c^\infty(D_{\A_F}^n)$ and $X \in A_{\GG_m}.$  In this section we study the analytic behavior as $|X|\to \infty$ of 
\begin{align*}
    \Sigma(X):=\sum_{\substack{\gamma' \in D^n\\P(\gamma')=0}}f\left(\frac{\gamma'}{X} \right).
\end{align*}
Assume henceforth that $P$ is homogeneous of degree $2$
and that $D$ is a quaternion division algebra over $F.$  This is the setting relevant for Theorem \ref{thm:general-case} and making this assumption now simplifies the notation required below.

Choose a pure tensor $\Phi=\otimes_v\Phi_v \in C^\infty_c(D_{\A_F}^2)$ such that the assumption of Lemma \ref{lem:delta} is satisfied, i.e. $\Phi(t,0)=0$ for all $t\in D_{\A_F}$ and $\int_{D_{\A_F}}\Phi(0,t)dt \neq 0.$  Applying the $\delta$-symbol expansion of Lemma \ref{lem:delta} (with $Q=X\gg_{\Phi} 1$), we have 
\begin{align}\label{eq:afterdelta}
    \Sigma(X)=\frac{c_{\Phi,X}}{|X|^4} \sum_{\gamma \in D^n}f\left(\frac{\gamma'}{X}\right)\sum_{\delta \in D^\times} \diff\left(\frac{\gamma'}{X},\frac{\delta}{X}\right),
\end{align}
where for convenience
\begin{align}\label{eq:defdiff}
\diff(\gamma',\delta)
:= \Phi\left(P(\gamma')\delta^{-1},\delta\right)
    - \Phi\left(\delta,\delta^{-1}P(\gamma')\right)
\end{align}
Here we have used the fact that $P(\gamma'/X) = P(\gamma')/X^2$.

Writing $\gamma'=(\gamma_1',\dots,\gamma_n')$, we define pairings
\begin{align}\label{eq:pairing}
\begin{split}
   \cdot: D_{\A_F}^n \times D_{\A_F}^n &\lto D_{\A_F}\\
   (\gamma,\gamma') &\longmapsto \sum_{i=1}^n \gamma_i\gamma_i',
\\\langle\,,\,\rangle:D_{\A_F}^{n} \times D_{\A_F}^{n} &\lto \A_F\\
    (\gamma,\gamma')
    &\longmapsto \mathrm{trd}(\gamma\cdot\gamma')
    = \sum_{i=1}^{n}\mathrm{trd}\,\gamma_i\gamma_i'.
\end{split}
\end{align}
Applying Poisson summation in $\gamma' \in D^n$ to \eqref{eq:afterdelta},  we have
\begin{align*}
 \Sigma(X)=   \frac{c_{\Phi,X}}{|X|^4}\sum_{\gamma \in D^{n}}\int_{D_{\A_F}^n}f\left(\frac{t}{X}\right) \sum_{\delta \in D^\times} \diff\left(\frac{t}{X},\frac{\delta}{X}\right) \psi\left(\langle \gamma,t\rangle \right) \, dt.
\end{align*}
By the change of variables $\gamma\mapsto \frac{\gamma}{\mathrm{nrd}(\delta)}$, we arrive at 
\begin{align*}
 \Sigma(X) = \frac{c_{\Phi,X}}{|X|^4} \sum_{\gamma \in D^{n}}\int_{D_{\A_F}^n}f\left(\frac{t}{X}\right)\sum_{\delta \in D^\times} \diff\left(\frac{t}{X},\frac{\delta}{X}\right) \psi\left(\frac{\langle \gamma,t\rangle}{\mathrm{nrd}(\delta)} \right) \, dt.
\end{align*}
Changing variables $t \mapsto Xt$ and using the fact that $\nrd(X) = X^2$, we obtain
\begin{align}\label{eq:after-Poisson}
 \Sigma(X) = c_{\Phi,X}|X|^{4(n-1)} \sum_{\gamma \in D^{n}}\int_{D_{\A_F}^n}f(t)\sum_{\delta \in D^\times} \diff\left(t,\frac{\delta}{X}\right) \psi\left(\frac{\langle \gamma/X,t\rangle}{\mathrm{nrd}(\delta/X)} \right) \, dt.
\end{align}

For $(\delta,\gamma) \in D_{\A_F}^\times \times D_{\A_F}^n,$ define
\begin{align}\label{eq:I0I1adelic}
\begin{split}
    I_{0}(\delta,\gamma):=I_{0}(f,\Phi,\delta,\gamma)
     &:= \int_{D_{\A_F}^n} f(t) \Phi(\delta, \delta^{-1} P(t))\psi\left(\frac{\langle \gamma,t\rangle}{\mathrm{nrd}(\delta)} \right) \, dt, \\
I_{1}(\delta,\gamma):=I_1(f,\Phi,\delta,\gamma) &:= \int_{D^n_{\A_F}} f(t) \Phi(P(t) \delta^{-1}, \delta) \psi\left(\frac{\langle \gamma,t\rangle}{\mathrm{nrd}(\delta)} \right)  \, dt. 
\end{split}
\end{align}
Then \eqref{eq:after-Poisson} can be written as
\begin{equation}\label{renormalize-after-poisson}
\begin{split}
\Sigma(X) &= c_{\Phi,X}X^{4(n-1)}
\sum_{\gamma \in D^{n}} \sum_{\delta \in D^\times}
\left(I_{1}\left(\frac{\delta}{X},\frac{\gamma}{X}\right) -  I_0\left(\frac{\delta}{X},\frac{\gamma}{X} \right)\right)
\end{split}
\end{equation}
Defining \eqref{eq:pairing} and  \eqref{eq:I0I1adelic} analogously over local fields, one has for pure tensors $f=\otimes_v f$
\begin{align*}
    I_i(f,\Phi,\delta,\gamma)=\prod_v I_{i,v}(f_v,\Phi_v,\delta,\gamma).
\end{align*}
To avoid potential confusion, we point out that since $X \in A_{\GG_m} \cong \RR_{>0},$ if $f$ and $\Phi$ are pure tensors then $I_{i}\left(\frac{\delta}{X},\frac{\gamma}{X}\right)=I_{i}\left(\delta,\gamma\right)$ for all $v \nmid \infty.$  In other words only the Archimedean factors $I_{i,\infty}$ (i.e.~those with $v\mid \infty$) in \eqref{renormalize-after-poisson} depend on the parameter $X.$

Moreover, by the change of variables $t\mapsto t^\dagger$ on $D_{\A_F}^n$, and the cyclicity of $\trd$, we have
\begin{equation}
\label{EQN:general-I0-I1-symmetry}
I_{0}(f,\Phi,\delta,\gamma)
= I_{1}(f\circ\dagger,\sw\circ\dagger,\delta^\dagger,\gamma^\dagger),
\end{equation}
where $f\circ\dagger$ is the map $t\mapsto f(t^\dagger)$,
and $\sw\circ\dagger$ is the map $(x,y)\mapsto \Phi(y^\dagger,x^\dagger)$.
So we may concentrate on either $I_{0}$ or $I_1$ for the most part.

Write $\Sigma(X)=\Sigma_0(X)+E_1(X),$ where
\begin{equation}\label{eq:term_gamma=0}
\Sigma_{0}(X):=c_{\Phi,X}X^{4(n-1)}\sum_{\delta \in D^\times}\left(I_{1}\left(\frac{\delta}{X},0\right) -  I_0\left(\frac{\delta}{X},0 \right)\right), 
\end{equation}
\begin{equation}\label{eq:term_gammaneq0}
E_1(X):=c_{\Phi,X}X^{4(n-1)}\sum_{\gamma \in D^n-\{0\}} \sum_{\delta \in D^\times} \left(I_{1}\left(\frac{\delta}{X},\frac{\gamma}{X}\right) -  I_0\left(\frac{\delta}{X},\frac{\gamma}{X} \right)\right). 
\end{equation}
We expect that $E_1(X) \ll_{f,\Phi,\epsilon} X^{3n+\epsilon}$ for all $\epsilon>0.$  In Theorem \ref{THM:off-center-total-bound}  we prove this when $F=\QQ,$ $D$ is nonsplit at $\infty$ and $2,$ and $f=\otimes_{v}f_v$, $\Phi:=\otimes_v\Phi_v$ with $f_v:=\one_{\OO_{D_v}^n}$ and  $\Phi_v =\one_{\OO_{D_v}^2}$
for all $v$ where $D_v$ is split.
This suffices for the proof of Theorem \ref{thm:general-case}.

\subsection{Contribution of $\Sigma_0$}
For the remainder of this section we focus on $\Sigma_0(X).$  For this contribution we do not assume that $F=\QQ$ because we do not require this assumption.
For functions $\phi:D^\times \backslash D_{\A_F}^\times \to \CC^\times$ and $s \in \CC$ let $\phi_s(x):=\phi(x)|\nrd(x)|^s.$

In order to state a lemma let us discuss measures on $A_{\GG_m}$ and $D_{\A_{F}}^\times.$
Let
 $$
 (D_{\A_F}^\times)^1:=\ker(|\nrd|:D_{\A_F}^\times \lto \RR_{>0}).
 $$
  Endow $\RR_{>0}$ with the usual measure $\frac{dt}{t}.$  Choose a Haar measure on $D_{\A_F}^\times.$  We assume that the map $|\nrd|:A_{\GG_m} \tilde{\to} \RR_{>0}$ is measure preserving and that the natural isomorphisms $D_{\A_F}^\times =A_{\GG_m} \times (D_{\A_F}^\times)^1 $ and $(D_{\A_F}^\times)^1 \tilde{\to} A_{\GG_m} \backslash D_{\A_F}^\times$ are measure preserving.  We then obtain a measure on $A_{\GG_m} D^\times \backslash D^\times_{\A_F}$ in the usual manner since $D^\times <A_{\GG_m} \backslash D^\times_{\A_F}$ is discrete. We define orthonormal bases on subspaces of $L^2(A_{\GG_m}D^\times \backslash D_{\A_F}^\times)$ using this measure.

\begin{lem} \label{lem:ker}
Let $h \in C_c^\infty(D_{\A_F}).$  Then
\begin{align}\label{eq:Mellin-inv}
\sum_{\delta\in D^\times} h(\delta g)=\frac{1}{2\pi i} \int_{\mathrm{Re}(s)=\sigma} \sum_{\pi} \sum_{\phi\in B_\pi} \left(\int_{D_{\A_F}^\times} h(x) \bar\phi_s(x) \, d^\times x\right) \phi_{-s}(g) \, ds,
\end{align}
for all $g\in D_{\A_F}^\times$ and $\sigma \gg 1$, where the sum over $\pi$ is taken over isomorphic classes of automorphic representations of $A_{\GG_m} D^\times \backslash D^\times_{\A_F}$ and $B_\pi$ is an orthonormal basis of $\pi$. 
\end{lem}

\begin{proof}
By Mellin inversion we have
\begin{align*}
\sum_{\delta \in D^\times} h(\delta g)=\frac{1}{2\pi i}\int_{\mathrm{Re}(s)=\sigma}\left(\sum_{\delta \in D^\times}\int_{A_{\GG_m}}h(a\delta g)|a|^sd^\times a\right)ds
\end{align*}
for $\sigma  \gg 1.$  
On the other hand, for $\mathrm{Re}(s)>2$ the function $g \mapsto \int_{A_{\GG_m}}h(ag)|a|^sda$ lies in $L^1(( D_{\A_{F}}^\times)^1)$ \cite[Lemma 12.5]{GodementJacquetBook}.
Moreover since $D$ is a division algebra, $A_{\GG_m} D^\times \backslash D_{\A_F}^\times$
is compact, and hence $L^2(A_{\GG_m}D^\times \backslash D_{\A_F}^\times)$ decomposes discretely as a representation of $D_{\A_F}^\times.$  Hence
\begin{align*}
    \sum_{\delta \in D^\times}\int_{A_{\GG_m}}h(a\delta g)|a|^sd^\times a&=\sum_{\pi}\sum_{\phi \in B_\pi}\int_{(D_{\A_F}^\times)^1}\int_{A_{\GG_m}}h(ax)|a|^sd a\overline{\phi}(x) d^\times x \phi_{-s}(g)\\
    &=\sum_{\pi}\sum_{\phi \in B_\pi}\int_{D_{\A_F}^\times}h(x)\overline{\phi}_s(x) d^\times x \phi_{-s}(g).
\end{align*}
\end{proof}

Now apply \eqref{eq:Mellin-inv} with $g=I_D$ to \eqref{eq:term_gamma=0}. We obtain 
\begin{align*}
    \begin{split}
        \Sigma_0(X)  &= \frac{c_{\Phi,X}|X|^{4(n-1)}}{2\pi i}\int_{\mathrm{Re}(s)=\sigma}  \sum_{\pi} \sum_{\phi\in B_\pi} \int_{D_{\A_F}^\times \times D_{\A_F}^n} f\left(t\right) \\& \times \left(\Phi\left(P(t)Xx^{-1},x\right)
    - \Phi\left(x,x^{-1}XP(t)\right)\right)    \, dt \bar\phi_s(x) \phi(I_D) \, d^\times x ds.
    \end{split}
\end{align*}
By the change of variables $x\mapsto Xx$, we obtain 
\begin{align*}
    \begin{split}
         \Sigma_0(X)&=  \frac{c_{\Phi,X}|X|^{4n-4}}{2\pi i }   \int_{\mathrm{Re}(s)=\sigma} |X|^{2s}\sum_{\pi}  \sum_{\phi\in B_\pi}   \\ 
         &\times \int_{D_{\A_F}^\times} \int_{D_{\A_F}^n} f\left(t\right) \left(\Phi\left(P(t)x^{-1},x\right)
    - \Phi\left(x,x^{-1}P(t)\right)\right)   \, dt \bar\phi_s(x) \phi(I_D) \, d^\times x ds \\ 
      &   =  \frac{c_{\Phi,X}|X|^{4n-4}}{2\pi i } \int_{\mathrm{Re}(s)=\sigma}  |X|^{2s} \sum_{\pi} \sum_{\phi\in B_\pi}  \left(I^1(f,\Phi,\bar\phi_s)-I^0(f,\Phi,\bar\phi_s)\right) \, ds, 
    \end{split}
\end{align*}
where 
\begin{align*}
\begin{split}
    I^0(f,\Phi,\bar\phi_s)
&:= \int_{D_{\A_F}^\times} \int_{D_{\A_F}^n}  f\left(t\right) \Phi \left(x,x^{-1}P(t)\right)  \, dt \bar\phi_s(x) \phi(I_D) \, d^\times x \\
&= \int_{D_{\A_F}^\times} I_0(f,\Phi,x,0) \bar\phi_s(x) \phi(I_D) \, d^\times x, 
\end{split}
\end{align*}
\begin{align*}
\begin{split}
    I^1(f,\Phi,\bar\phi_s)
&:= \int_{D_{\A_F}^\times} \int_{D_{\A_F}^n}  f\left(t\right) \Phi \left(P(t)x^{-1},x\right)  \, dt \bar\phi_s(x) \phi(I_D) \, d^\times x \\ 
&= \int_{D_{\A_F}^\times} I_1(f,\Phi,x,0) \bar\phi_s(x) \phi(I_D) \, d^\times x . 
\end{split}
\end{align*}
for cusp forms $\phi:A_{\GG_m}D^\times \backslash D_{\A_F}^\times \to \CC.$ 

Let $K_\infty < D_{\infty}^\times$ be a maximal compact subgroup and let $K:=K_\infty \widehat{\OO}_D^\times=\prod_v K_v$.  
We assume henceforth that $\Phi$ is bi-$K$-invariant in both entries. We also assume that 
\begin{align} \label{conj:inv}
f(kt_1k^{-1},\dots,kt_nk^{-1})=f(t_1,\dots,t_n)
\end{align}
for $(t_1,\dots,t_n) \in D^n_{\A_F}$ and $k \in K_v$ for $v$ such that $D_v$ is split.
Then only $\pi$ with $\pi^K \neq 0$ contribute a nonzero summand to $\Sigma_0(X)$ (this follows from the same argument proving Lemma \ref{lem:unram}).  For $\pi^K \neq 0,$ fix a vector $\phi_\pi\in \pi^K$ of $L^2$-norm $1.$  The vector $\phi_{\pi}$ is unique up to a constant of norm $1.$  It follows that  
\begin{align}
 \Sigma_0(X)&=\frac{c_{\Phi,X}|X|^{4n-4}}{2\pi i }\int_{\mathrm{Re}(s)=\sigma}\sum_{\pi:\pi^K \neq 0}|X|^{2s}\left(I^1(f,\Phi,\overline{\phi}_{\pi s})-I^0(f,\Phi,\overline{\phi}_{\pi s})\right)ds.
\end{align}

\quash{

This should be a consequence of the local estimations below. These estimations should also allow us to shift the contour. \textcolor{red}{The following two classical results may be preciser/unnecessary if one assumes that $D$ is quaternion.} It is known (see \cite[Theorem 19.4.3]{Getz:Hahn}) that the image of Jacquet-Langlands correspondence consists of discrete automorphic representations that are $D$-compatible. By Moeglin-Waldspurger's theorem (see \cite[Theorem 10.8.1]{Getz:Hahn}), any discrete automorphic representation of $\GL_d$ is isomorphic to a Speh representation $(\sigma, m)$ where $m|d$ and $\sigma$ is a cuspidal automorphic representation of $\GL_{d/m}$. }

\begin{lem}\label{lem:matrix coeff decomp of I}
Assume that $\pi^K \neq 0.$ One has
    \begin{align*}
        I^1(f,\Phi,\bar{\phi_{\pi s}})&=|\phi_{\pi}(I_D)|^2 \int_{D_{\A_F}^\times} \int_{D_{\A_F}^n} f(t) \Phi \left(P(t)x^{-1},x\right)   \bar{m_\pi(x)}|\nrd(x)|^s\, dt  d^\times x\\
        &=|\phi_\pi(I_D)|^2 \prod_v \int_{D_{v}^\times} \int_{D_{v}^n} f_v(t_v) \Phi_v \left(P(t_v)x_v^{-1},x_v\right)   \bar{m_{\pi_v}(x_v)}|\nrd(x_v)|_v^s\,  dt_v d^\times x_v.
    \end{align*}
    Here $m_\pi(x)$ is the matrix coefficient $\langle\pi(x)\phi_\pi,\phi_{\pi}\rangle_\pi$, and $m_{\pi_v}$ is the zonal spherical function of $\pi_v$. A similar formula holds for $I^0$.
\end{lem}
\noindent In the lemma we realize $\pi$ as a subrepresentation of $L^2(A_{\GG_m} D^\times \backslash D_{\A_F}^\times),$ and  $\langle\,,\, \rangle_{\pi}$ is the pairing on the space of $\pi$ given by restriction of the pairing on $L^2(A_{\GG_m} D^\times \backslash D_{\A_F}^\times).$

\begin{proof}
    By $K$-invariance,
    \begin{align*}
        I^1(f,\Phi,\bar\phi_{\pi s})&=\frac{1}{\vol(K)}\phi_{\pi}(I_D)\int_{D_{\A_F}^\times} \int_{D_{\A_F}^n} f\left(t\right) \Phi \left(P(t)x^{-1},x\right)   \left(\int_K\bar\phi_\pi(kx)dk\right)|\nrd(x)|^s \, dt  \, d^\times x.
    \end{align*}
We have
\begin{align*}
\int_K\phi_\pi(kx)dk=\mathrm{vol}(K)\phi_{\pi}(I_D)\langle \pi(x)\phi_{\pi},\phi_{\pi} \rangle_{\pi}
\end{align*}
since spherical functions are unique up to a constant multiple.   Hence,
    \begin{align*}
        I^1(f,\Phi,\bar{\phi_{\pi s}})
        &=|\phi_{\pi}(I_D)|^2\int_{D_{\A_F}^\times} \int_{D_{\A_F}^n} f\left(t\right) \Phi \left(P(t)x^{-1},x\right)   \bar{\langle\pi(x)\phi_{\pi},\phi_{\pi}\rangle_\pi}|\nrd(x)|^s \, dt  \, d^\times x.
    \end{align*}
    The last identity in the lemma follows from the fact that the global matrix coefficient attached to a cusp form is factorizable if the cusp form is factorizable.
\end{proof}

Using notation from \eqref{loc:sp:int} below, we have shown
\begin{align} \label{eq:factor}
I^i(f,\Phi,\overline{\phi}_{\pi s})=|\phi_\pi(I_D)|^2 \prod_vI^i ( f_v,\Phi_v,\overline{m}_{\pi_v s}).
\end{align}
Thus we can use the work in \S \ref{sec:nonarch:sp} and \S \ref{sec:arch:bounds} to understand these functions. 

At this point it is convenient to recall some basic automorphic representation theory.  Unitary automorphic representations of $D_{\A_F}^\times$ are either $1$-dimensional or infinite dimensional.  The one-dimensional representations are precisely of the form $\chi \circ \nrd$ for a character $\chi:F^\times \backslash \A_F^\times \to \CC^\times.$  
 If $\pi$ is the one-dimensional representation $\chi \circ \mathrm{nrd},$ then its standard $L$-function is 
$$
L(s,\pi)=L(s+\tfrac{1}{2},\chi)L(s-\tfrac{1}{2},\chi).
$$   
The infinite dimensional representations correspond under the Jacquet-Langlands correspondence to cuspidal automorphic representations of $\GL_2(\A_F)$ (although not all representations of $\GL_2(\A_F)$ are in the image of the Jacquet-Langlands correspondence).  The (standard) $L$-function $L(s,\pi)$ of an infinite-dimensional automorphic representation of $D^\times_{\A_F}$ is by definition the standard $L$-function of its Jacquet-Langlands transfer. 

There is only one setting in which we have to understand the local Jacquet-Langlands transfer of an irreducible admissible representation $\pi_v$ of $D_v^\times$ when $D_v$ is nonsplit.  This is the case where $\pi_v \cong |\mathrm{nrd}|^{it}.$  In this case the Jacquet-Langlands transfer is $\mathrm{St}|\det|^{it},$ where $\mathrm{St}$ is the Steinberg (also known as special) representation \cite[Above Proposition 15.5]{JL:GL2}.  We point out that $L(s,\mathrm{St}|\cdot|^{it})=\zeta_v(s+it+\tfrac{1}{2})$ \cite[Theorem 6.15]{Gelbart:JL}.
    
We will also require the completed zeta function $\zeta_D(s):=\prod_{v}\zeta_{D_v}(s)$ of $D.$  Let $\Gamma_{\RR}(s):=\pi^{-s/2}\Gamma(s/2).$  We then define 
\begin{align}\label{eq:trivL}
\begin{split}
    \zeta_{D_v}(s):=\begin{cases}
   \Gamma_{\RR}(s+\tfrac{1}{2})\Gamma_{\RR}(s+\tfrac{3}{2})=2(2\pi)^{-(s+1/2)}\Gamma(s+\tfrac{1}{2})
        & \textrm{ if } F_v=\RR, D_v \textrm{ nonsplit, } \\
\zeta_{v}\left(s+\frac{1}{2}\right)& \textrm{ if }  v \nmid \infty, D_v \textrm{ nonsplit, } \\
       \zeta_v\left(s+\frac{1}{2}\right)\zeta_v\left(s-\frac{1}{2}\right) & \textrm{ if } D_v \textrm{ split. }
    \end{cases}
\end{split}
\end{align}
Up to harmless normalizations this is the same as the definition in \cite[\S 29.5-29.6]{voight2021quaternion}.
We note that $\zeta_{D}(s) \neq L(s,1 \circ \mathrm{nrd}),$ although the corresponding local identity is true for all places $v$ where $D_v$ is split.
The function $\zeta_D(s)$ is meromorphic.  Because we assumed $D$ is nonsplit it is holomorphic apart from  simple poles at $s=-\frac{1}{2},\frac{3}{2}.$

Let $Z_{D_\infty^\times}$ be the center of $D_\infty^\times.$  We point out that if $F=\QQ,$ then $Z_{D_\infty^\times} \backslash D_\infty^\times$ is compact if and only if $A_{\GG_m} \backslash D_\infty^\times$ is compact.

\begin{thm}\label{thm:globalspectral}
Suppose $n\ge 5$ and $\pi^K\neq 0$. Assume $Z_{D_\infty^\times} \backslash D_\infty^\times$ is compact, and that  $(f,\Phi)=(f_\infty\one_{\widehat{\OO}_D^n},\Phi_\infty\one_{\widehat{\OO}_D^2})$ where $\Phi_\infty$ is bi-invariant under $K_\infty$ in both entries.  The quotient
\begin{align*}
    \frac{I^i(f,\Phi,\overline{\phi}_{\pi s})}{L\left(s+\frac{3}{2},\pi\right)}
\end{align*}
is holomorphic on $\sigma=\mathrm{Re}(s)>-\frac{n-1}{2},$
and $I^i(f,\Phi,\overline{\phi}_{\pi s})$ is rapidly decreasing on vertical strips (away from its poles).  If $\pi$ is the trivial representation, then the same is true if we replace $L(s+\tfrac{3}{2},\pi)$ by $\zeta_D(s+\tfrac{3}{2}).$
\end{thm}
\begin{proof}
     By \eqref{eq:factor}, the assertion follows from Propositions \ref{prop:nonarchsplittrivconv}, \ref{prop:nonarchnonsplittrivconv}, and \ref{prop:archconv}, our comments on the Jacquet Langlands correspondence above, and standard bounds on $L$-functions \cite[Lemma 5.2]{Iwaniec:Kowalski}.
\end{proof}

Let $1:D^\times \backslash D_{\A_F}^\times \to \CC$ be the function that is identically $1.$  It spans the trivial subrepresentation of $L^2(A_{\GG_m}D^\times \backslash D_{\A_F}^\times).$ We assume $1$ on $A_{\GG_m}D^\times \backslash D^\times_{\mathbb{A}_{F}}$ has $L^2$-norm $1,$ and thus the measure $d^\times x$ in \eqref{eq:Mellin-inv} is chosen so that $\vol(A_{\GG_m}D^\times \backslash D^\times_{\mathbb{A}_{F}})=1$. 

\begin{lem}\label{lem:nontrivspec} Suppose $n\ge 5$. Assume $A_{\GG_m} \backslash D_\infty^\times$ is compact, $F=\QQ$ and $(f,\Phi)$ is as above.  For any $\epsilon>0$ one has
    \begin{align*}
    \Sigma_0(X)=\frac{c_{\Phi,X}|X|^{4n-4}}{2\pi i} \int_{\mathrm{Re}(s)=\sigma} |X|^{2s}(I^{1}(f,\Phi,\bar{1}_s)-I^{0}(f,\Phi,\bar{1}_s)) \, ds+O_{\epsilon,f,\Phi}(|X|^{3n-3+\epsilon}).
    \end{align*}
\end{lem}
\begin{proof}
Consider 
\begin{align} \label{negl}
\frac{c_{\Phi,X}|X|^{4n-4}}{2\pi i} \int_{\mathrm{Re}(s)=\sigma} \sum_{\pi: \pi \neq \mathrm{triv}} |X|^{2s}(I^{1}(f,\Phi,\overline{\phi}_{\pi s})-I^{0}(f,\Phi,\overline{\phi}_{\pi s}))\, ds
\end{align}
where the sum is over nontrivial automorphic representations of $A_{\GG_m} \backslash D^\times_{\A_F}$ such that $\pi^K$ is nonzero.  The lemma is equivalent to the assertion that \eqref{negl} is $O_{\epsilon,f,\Phi}(|X|^{3n-3+\epsilon}).$  

Since $A_{\GG_m} \backslash D_\infty^\times$ is compact, the quotient
$A_{\GG_m}D^\times \backslash  D_{\A_F}^\times/K$
is finite. Thus $L^2(A_{\GG_m}D^\times \backslash D_{\A_F}^\times/K)$ is finite dimensional and we deduce that the sum over $\pi$ in \eqref{negl} is finite.  Hence it suffices to prove that the contribution of any given $\pi$ is $O_{\epsilon,f,\Phi,\pi}(|X|^{3n-3+\epsilon})$. 

 Since the complete $L$-function $L(s,\pi)$ (i.e. with the Archimedean factor) is entire, by Theorem \ref{thm:globalspectral} we can perform a contour shift
\begin{align*}
&\frac{c_{\Phi,X}|X|^{4n-4}}{2\pi i} \int_{\mathrm{Re}(s)=\sigma} |X|^{2s}(I^{1}(f,\Phi,\overline{\phi}_{\pi s})-I^{0}(f,\Phi,\overline{\phi}_{\pi s}))\, ds\\&=\frac{c_{\Phi,X}|X|^{4n-4}}{2\pi i} \int_{\mathrm{Re}(s)=-\tfrac{n-1}{2}+\epsilon} |X|^{2s}(I^{1}(f,\Phi,\overline{\phi}_{\pi s})-I^{0}(f,\Phi,\overline{\phi}_{\pi s})) \, ds\\
&=O_{\epsilon,f,\Phi,\pi}(|X|^{3n-3+2\epsilon}).
\end{align*}    
\end{proof}

\begin{rem}
    With more local work at the Archimedean places one could remove the assumption that $A_{\GG_m} \backslash D_\infty^\times$ is compact and $F=\QQ.$
\end{rem}

\begin{lem}
\label{LEM:formal-cancellation}
For $(f,\Phi)$ as in Theorem \ref{thm:globalspectral}, one has 
\begin{align*}
    \mathrm{Res}_{s=0}\left( I^1(f,\Phi,\bar1_s)-I^0(f,\Phi,\bar1_s)\right) =0. 
\end{align*}
\end{lem}

\begin{proof}
    From \eqref{EQN:general-I0-I1-symmetry}, we deduce that 
    \begin{align*}
        I^0(f_v, \Phi_v, \bar1_s)=I^1(f_v\circ\dagger, \sw_v\circ\dagger, \bar1_s). 
    \end{align*}
    But for all finite $v,$ we have $\one_{\OO_{D_v}^n}=\one_{\OO_{D_v}^n}\circ\dagger$ and $\one_{\OO_{D_v}^2}=\one_{\OO_{D_v}^2}^{\mathrm{sw}}\circ\dagger$. Therefore, by Theorem \ref{thm:globalspectral} and \eqref{eq:trivL}, it suffices to verify that  
    \begin{align*}
        I^1(f_\infty,\Phi_\infty,\bar1_0)-I^0(f_\infty,\Phi_\infty,\bar1_0)=0. 
    \end{align*}    
By Fourier inversion we have
    \begin{align*}
        &I^1(f_\infty,\Phi_\infty,\bar1_0)\\
         &=\int_{D_\infty^\times } \int_{D_\infty^n}f_\infty(Y)\left(\int_{D_\infty}\psi\left(\langle P(Y),Z\rangle\right) \mathcal{F}_1(\Phi_\infty)(-xZ,x)|\nrd(x)|^{2} \, dZ\right) dYdx
\end{align*}
By the bounds in Lemma \ref{lem:Archzeta} all of the integrals here converge absolutely.  Thus by the Fubini-Tonelli theorem  we can rearrange the integrals to see that the above is
\begin{align*}
        &\int_{D_\infty^\times \times D_\infty} \left(\int_{D_\infty^n}f_\infty(Y)\psi\left(\langle P(Y),Z\rangle\right)\, dY \right)\mathcal{F}_1(\Phi_\infty)(-xZ,x)|\nrd(x)|^{2} \, dZ dx.
\end{align*}
Here and for the rest of the proof we use capital letters to denote Haar measures on the additive group $D_{\infty}$ and lowercase letters to denote Haar measures on the multiplicative group $D^\times_\infty.$  Thus $|\nrd(x)|^2dx=dX.$  The above is
\begin{align*}
        &\int_{D_\infty \times D_\infty} \left(\int_{D_\infty^n}f_\infty(Y)\psi\left(\langle P(Y),Z\rangle\right)\, dY \right)\left(\int_{D_\infty}\Phi_\infty(H,X)\psi(\langle H, -XZ\rangle) dH dX \right)dZ \\
        &=\int_{D_\infty \times D_\infty} \left(\int_{D_\infty^n}f_\infty(Y)\psi\left(\langle P(Y),Z\rangle\right)\, dY \right)\left(\int_{D_\infty}\Phi_\infty(H,X)\psi(\langle -ZH, X\rangle) dX dH\right) dZ \\
        &=\int_{D_\infty^\times \times D_\infty} \left(\int_{D_\infty^n}f_\infty(Y)\psi\left(\langle P(Y),Z\rangle\right)\, dY \right)\mathcal{F}_2(\Phi_\infty)(h,-Zh)|\nrd(h)|^2  dZ dh \\
        &=I^0(f_\infty,\Phi_\infty,\bar1_0).
    \end{align*}  
\end{proof}
    

\begin{proof}[Proof of Theorem \ref{thm:general-case}]
We have $\Sigma(X)=\Sigma_0(X)+E_1(X)$ as explained above \eqref{eq:term_gamma=0} and \eqref{eq:term_gammaneq0}.  As mentioned earlier, $E_1(X) \ll_{\epsilon,f,\Phi} X^{3n+\epsilon}$ for all $\epsilon>0$ by Theorem \ref{THM:off-center-total-bound}.  Thus 
by Lemma \ref{lem:nontrivspec}
we have 
\begin{align*}
    \Sigma(X)=\frac{c_{\Phi,X}|X|^{4n-4}}{2\pi i} \int_{\mathrm{Re}(s)=\sigma} |X|^{2s}(I^{1}(f,\Phi,\bar{1}_s)-I^{0}(f,\Phi,\bar{1}_s)) \, ds+O_{\epsilon,f,\Phi}(|X|^{3n+\epsilon}).
    \end{align*}
We perform a countour shift using 
 Theorem \ref{thm:globalspectral} and Lemma \ref{LEM:formal-cancellation} to obtain
\begin{align} \label{for:leading}
    \Sigma(X)&=c_{\Phi,X}\bigg(\mathrm{Res}_{s=-2} I^{1}(f,\Phi,\bar{1}_s)-I^{0}(f,\Phi,\bar{1}_s)\bigg)|X|^{4n-8}+O_{\epsilon,f,\Phi}(|X|^{3n+\epsilon}).
\end{align}
\end{proof}

Thus we have reduced Theorem \ref{thm:general-case} to a study of the Eulerian integrals $I^i(f,\Phi,\phi_s)$ and the error term $E_1(X).$ The error term $E_1(X)$ will be treated in \S\S~\ref{SEC:general-hessian-based-analysis}--\ref{sec:geo:num}. We treat the integrals $I^i(f,\Phi,\phi_s)$ first, starting with the non-Archimedean case in the following section. Before this we pause to explain the proof of Theorem \ref{thm:constant}, modulo some local results contained below.
In contrast to the rest of this section, for the proof of Theorem \ref{thm:constant} we will make use of our assumption that $F=\QQ,$ although with more effort this assumption could be removed.
\begin{proof}[Proof of Theorem \ref{thm:constant}]
By Lemma \ref{lem:delta} and \eqref{for:leading} the constant $c(f_\infty \one_{\widehat{\OO}_D^\infty})$ is
\begin{align*}
    \frac{1}{\mathcal{F}_2(\Phi)(0,0)}\bigg(\mathrm{Res}_{s=-2} I^{1}(f,\Phi,\bar{1}_s)-I^{0}(f,\Phi,\bar{1}_s)\bigg).
\end{align*}
Using Theorem \ref{thm:globalspectral} this is 
\begin{align*}
    &\frac{\mathrm{Res}_{s=-2}\zeta_D(s+3/2)}{\mathcal{F}_2(\Phi_\infty)(0,0)\vol(\widehat{\OO}_{D})}\frac{I^{1}(f,\Phi,\bar{1}_s)-I^{0}(f,\Phi,\bar{1}_s)}{\zeta_D(s+3/2)}\Bigg|_{s=-2}\\
    &=\frac{\mathrm{Res}_{s=-1/2}\zeta_D(s)}{\mathcal{F}_2(\Phi_\infty)(0,0)}\frac{I^{1}(f_\infty,\Phi_\infty,\bar{1}_s)-I^{0}(f_\infty,\Phi_\infty,\bar{1}_s)}{\zeta_{D_\infty}(s+3/2)}\Bigg|_{s=-2}\vol(\widehat{\OO}_D^\times)\prod_{p} c(\one_{\OO_{D_p}^n}).
\end{align*}
Here we have used the fact that $I^1(\one_{\OO_{D_p}^n},\one_{\OO_{D_p}^2},\bar{1}_s)=I^0(\one_{\OO_{D_p}^n},\one_{\OO_{D_p}^2},\bar{1}_s)$ for all $p,$ \eqref{loc:dens}, \eqref{eq:c(f):bis}, and \eqref{eq:nonsplitlocaldensity}. 

Since $\Sigma(X)$ is independent of $\Phi_\infty,$  Theorem \ref{thm:general-case} implies that the factor at infinity above is also independent of $\Phi_\infty.$
We claim, moreover, that it remains unchanged if we replace $\Phi_\infty$ by any function $\Phi_\infty' \in \mathcal{S}(D_\infty^2)$ that is  bi-$K_\infty$-invariant in both variables. To see this, choose a sequence $\{\Phi_i\}\subset C^\infty_c(D_\infty^2)$ consisting of functions that are bi-$K_\infty$-invariant in both variables such that $\Phi_i \to \Phi_\infty'$ in $\mathcal{S}(D_\infty^2).$  Then 
\begin{align*}
    \frac{\mathrm{Res}_{s=-1/2}\zeta_D(s)}{\mathcal{F}_2(\Phi_i)(0,0)}\frac{I^{1}(f_\infty,\Phi_i,\bar{1}_s)-I^{0}(f_\infty,\Phi_i,\bar{1}_s)}{\zeta_{D_\infty}(s+3/2)}\Bigg|_{s=-2}
\end{align*}
is independent of $\Phi_i.$  On the other hand, it converges to 
\begin{align}\label{eq:Archcomp}
    \frac{\mathrm{Res}_{s=-1/2}\zeta_D(s)}{\mathcal{F}_2(\Phi_\infty')(0,0)}\frac{I^{1}(f_\infty,\Phi_\infty',\bar{1}_s)-I^{0}(f_\infty,\Phi_\infty',\bar{1}_s)}{\zeta_{D_\infty}(s+3/2)}\Bigg|_{s=-2}.
\end{align}
This justifies the claim. 

Using the notation in the proof of Lemma \ref{lem:Archnonvanish}, take $\Phi_\infty':=\Phi_{2,2}-\Phi_{2,1}$. Then using \eqref{eq:I1specific} and \eqref{eq:I0specific} we see that \eqref{eq:Archcomp} equals
\begin{align*}
    -\mathrm{Res}_{s=-\tfrac{1}{2}}\zeta_D(s)\frac{\vol((D^\times_\infty)^1)}{2}\int_{D} \left(\int_{D^n} f(Y)\psi\left(\left\langle P(Y),Z\right\rangle\right)\, dY\right)\, dZ.
\end{align*}
By  \cite[Theorem 29.10.23, Lemma 29.5.18, Lemma 29.8.24]{voight2021quaternion}  we have
\begin{align*}
    \mathrm{Res}_{s=-\tfrac{1}{2}}\zeta_D(s)=-\frac{\mathrm{Res}_{s= 1}\zeta_\QQ(s)}{2\pi^2\widehat{\tau}^\times(\widehat{\OO}_D^\times)}=-\frac{1}{2\pi^2 \widehat{\tau}^\times(\widehat{\OO}_D^\times)}.
\end{align*}
The definition of $\widehat{\tau}^\times(\widehat{\OO}_D^\times)$ will be recalled in a moment.
We conclude that
\begin{align} \label{main:term}
    c(f_\infty \one_{\widehat{\OO}_D^\infty})=\frac{\vol(K)}{4\pi^2 \widehat{\tau}^\times(\widehat{\OO}_D^\times)}c(f_\infty)\prod_{p} c(\one_{\OO_{D_p}^n}).  
\end{align}

Let us be a little more explicit about the measures in this expression.    Above Lemma \ref{lem:nontrivspec} we assumed that the volume of $A_{\GG_m}D^\times \backslash D_{\A_F}^\times$ is $1.$  There is a unique Haar measure on  $A_{\GG_m} \backslash D_{\A_F}^\times$ that induces this measure.  By \cite[Theorem 29.11.3]{voight2021quaternion}, this measure is the Tamagawa measure.  It decomposes as a product of a Haar measure on $A_{\GG_m} \backslash D_{F_\infty}^\times$ and a Haar measure $\widehat{\tau}^\times$ on $D^\times_{\A_F^\infty}.$
Since $A_{\GG_m} \backslash D^\times_{\infty}$ is compact, the restriction of the Tamagawa measure on $A_{\GG_m} \backslash D^\times_{\A_F}$ to $K$ induces a Haar measure on $K.$  This measure is used to compute the volume in \eqref{main:term}.  

By \cite[(29.8.15), Lemma 29.5.9]{voight2021quaternion}, we have
$    \vol(K)=4\pi^2\widehat{\tau}^\times(\widehat{\OO}_D^\times)$.
Thus $c(f_\infty \one_{\widehat{\OO}_{D}^n})=c(f_\infty)\prod_{p} c(\one_{\OO_{D_p}^n}).$  The nonvanishing statement follows from Corollary \ref{cor:nonvanish split local density} and Lemmas \ref{lem: nonsplit local density} and \ref{lem:Archop}.
\end{proof}

\section{Spectral bounds for split cases} \label{sec:nonarch:sp}

Let $v$ be a place of $F$ which we drop from notation, writing $F:=F_v, D:=D_v,$ etc. We do not yet make any assumption on $F$ or $D$.  Let $\pi$ be an irreducible admissible unitary representation of $D^\times.$  We realize $\pi$ in the category of smooth representations.   Let $m$ be a matrix coefficient of $\pi$ and let 
\begin{align*}
    m_s(g):=m(g)|\mathrm{nrd}(g)|^s.
\end{align*}
For $(f,\Phi)\in \mathcal{S}(D^n)\times\mathcal{S}(D^2)$  and $s\in \CC$, define
\begin{align} \label{loc:sp:int}
I(f,\Phi,m_s):=\int_{D^\times}\int_{D^n} f(Y)\Phi\left(P(Y)g^{-1},g\right)m_s(g)\, dY dg.
\end{align}

Since $\pi$ is unitary, the function $g \mapsto \int_{D^n} f(Y)\Phi\left(P(Y)g^{-1},g\right)m(g)\, dY$ is bounded by the restriction of a Schwartz function in $\mathcal{S}(D)$ to $D^\times.$  Hence the integral above converges absolutely for $\mathrm{Re}(s)>1$ \cite[Proposition 1.1]{GodementJacquetBook}.
 By Fourier inversion, for $\mathrm{Re}(s) >1$ the integral \eqref{loc:sp:int} is
\begin{align*}
&\int_{D^n}f(Y)\int_{D^\times}\left(\int_D\mathcal{F}_1(\Phi)(Z,g)\psi(-\langle P(Y)g^{-1},Z\rangle )m_s(g) \, dZ\right)\, dgdY\\
    &=\int_{D^n}f(Y)\left(\int_D\psi(\langle P(Y),Z\rangle )\int_{D^\times} \mathcal{F}_1(\Phi)(-gZ,g)m_{s+2}(g)\, dg dZ\right)\, dY
\end{align*}
For a given $Z\in D$, define the zeta integral
\begin{align*}
    \mathcal{Z}_Z(\Phi,s,m):=\int_{D^\times} \mathcal{F}_1(\Phi)(-gZ,g)m_{s+2}(g)\, dg,
\end{align*}
and define
\begin{align*}
    H_Z(\Phi,s,m):=\frac{\mathcal{Z}_Z(\Phi,s,m)}{L(s+3/2,\pi)}.
\end{align*}
When $\Phi$ is fixed, we often write $\mathcal{Z}_Z(s,m):=\mathcal{Z}_Z(\Phi,s,m)$ and $H_Z(s,m):=H_Z(\Phi,s,m)$. 

\begin{lem} \label{lem:rationality}
If $F$ is non-Archimedean, $\mathcal{Z}_Z(\Phi,s,m)$ lies in $\CC(q^{-s})$ and $H_Z(\Phi,s,m)$ lies in $\CC[q^{-s},q^s].$
\end{lem}

\begin{proof}
    This follows from \cite[Theorem 3.3]{GodementJacquetBook}. 
\end{proof}

Assume now $F$ is non-Archimedean and $D$ is split, so we identify $D$ with $M_2(F)$. Assume $\psi$ is unramified and fix $\Phi=\one_{M_2(\OO_F)^2}$. We normalize the measure on $\GL_2(F)$ so that $dg(\GL_2(\OO_F))=1$. In the rest of the section, we make the following assumption:
\begin{align*}
f \textrm{ is invariant under conjugation by }\GL_2(\OO_F) \textrm{ (acting diagonally)}.
\end{align*}
Then $I(f,\Phi,m_s)$ is nonzero only if $\pi$ is unramified and 
$m$ is left and right $\GL_2(\OO_F)$-invariant. 
 We henceforth assume this. Moreover, we shall normalize $m$ so that $m(e)=1$. This is compatible with Lemma \ref{lem:matrix coeff decomp of I}.

We have $\mathcal{F}(\one_{M_2(\OO_F)})=\one_{M_2(\OO_F)},$ so $H_Z(s,m)$ is invariant under $Z \mapsto Z+W$ for $W \in M_2(\OO_F).$ We make an additional assumption:
\begin{align*} 
\textrm{ If } f(Y)\neq 0, \textrm{ then } P(Y)\in M_2(\OO_F).
\end{align*}
Under this assumption, we can write
\begin{align}\label{eq:simplified}
    \frac{I(f,\one_{M_2(\OO_F)^2},m_s)}{L(s+3/2,\pi)}&=\int_{M_2(F)^n}f(Y)\left(\sum_{Z\in M_2(F)/M_2(\OO_F)} \psi(\langle P(Y),Z\rangle )H_Z(s,m)\right) \, dY.
\end{align}

Recall the box norm defined in \eqref{box-vector-norm}, which is the matrix norm in this case.

\begin{lem} \label{lem:reps}
    Let $Z\in M_2(F)$. There exists $W\in M_2(\OO_F)$ so that $Z+W \neq 0,$
     $|\det(Z+W)| \ge \norm{Z+W}$ (or equivalently $(Z+W)^{-1}\in M_2(\OO_F)$) and $|\mathrm{tr}(Z+W)|\ge |2|.$ 
\end{lem}
\begin{proof}
    We can and do assume that $Z\in \GL_2(F)$. By the Cartan decomposition, there exist $k_1,k_2\in\GL_2(\OO_F)$ and integers $n_1\ge n_2$ such that $Z=k_1\begin{psmatrix}
        \varpi^{n_1}&\\&\varpi^{n_2}
    \end{psmatrix}k_2$. Let $W=k_1 \begin{psmatrix}
        \lambda_1&\\&\lambda_2
    \end{psmatrix} k_2$, where
    \begin{align*}
        \lambda_i=\begin{cases}
            1-\varpi^{n_i}&\text{ if }n_i>0,\\
            0&\text{ otherwise,}
        \end{cases}\quad\text{ for }i=1,2.
    \end{align*}
    It is clear that $W\in M_2(\OO_F)$  and $|\det(Z+W)| \ge \norm{Z+W}$.
    
    We may therefore assume $ Z=(z_{ij})$ satisfies $|\det(Z)|\ge\norm{Z},$ so $|\det(Z)|\ge 1$. Suppose $|\mathrm{tr}(Z)|< |2|\le1$. If $|z_{11}|<\norm{Z}$, choose $W=\begin{psmatrix}
    0&0\\0&1
\end{psmatrix},$ so
\begin{align*}
    |\mathrm{tr}(Z+W)|=|\mathrm{tr}(Z)+1|=1.
\end{align*}
Notice that $\det(Z+W)=\det(Z)+z_{11}$ and we have $|\det(Z)|\ge\norm{Z}>|z_{11}|$ and $|\det(Z)|\ge1$, so
\begin{align*}
    |\det(Z+W)|=|\det(Z)|\ge\norm{Z}\ge\norm{Z+W}.
\end{align*}

Assume $|z_{11}|=\norm{Z}.$  In this case we claim that we can take $W=\begin{psmatrix}
    0&0\\0&\alpha
\end{psmatrix}$
for an appropriately chosen $\alpha.$
In more detail, we choose $\alpha=2$ if $2|q.$  If $q$ is odd and $p$ is the prime dividing $q,$ we choose $\alpha \in \{1,\dots,p-1\}$ so that 
$|\det (Z+W)|=|\det(Z)|.$  This is possible because 
$\det(Z+W)=\det(Z)+\alpha z_{11}$
and $|\det(Z)|\ge|z_{11}|$. Since
\begin{align*}
    \norm{Z+W}=\max(|z_{11}|,|z_{22}+\alpha|)\le\max(\norm{Z},1)\le|\det(Z)|,
\end{align*}
and $|\mathrm{tr}(Z+W)|=|\mathrm{tr}(Z)+\alpha|=|\alpha|\ge|2|$, the assertion is justified.
\end{proof}

By the above lemma, for each coset in $M_2(F)/M_2(\OO_F)$, we can and do choose a representative $Z$ such that $Z^{-1}\in M_2(\OO_F)$. Then
\begin{align*}
    \one_{M_2(\OO_F)}(gZ)=1 \Rightarrow  \one_{M_2(\OO_F)}(g)=1,
\end{align*}
and
\begin{align}\label{eq:normalizedzeta}
\nonumber \mathcal{Z}_Z(s,m)&=\int_{\GL_2(F)}\one_{M_2(\OO_F)}(gZ)m_{s+2}(g)\, dg\\
    &=|\det(Z)|^{-(s+2)} \int_{\GL_2(F)}\one_{M_2(\OO_F)}(g)m(gZ^{-1})|\det(g)|^{s+2}\, dg.
\end{align}

\begin{lem}\label{lem:HZ}
We have
\begin{align*}
    H_{Z}(s,m)= m(Z^{-1})|\det(Z)|^{-(s+2)} .
\end{align*}
\end{lem}
\begin{proof}
 According to \eqref{eq:normalizedzeta}, we need to compute 
\begin{align*}
    \mathcal{Z}_{Z}(s,m)&=|\det(Z)|^{-(s+2)} \int_{\GL_2(F)}\one_{M_2(\OO_F)}(g)m(gZ^{-1} )|\det g|^{s+2}\, dg\\
    &=|\det(Z)|^{-(s+2)} \int_{\GL_2(F)}\one_{M_2(\OO_F)}(g)|\det g|^{s+2}\int_{\GL_2(\OO_F)}m(gk Z^{-1})\, dk dg.
\end{align*}
Observe that the function
    \begin{align*}
        g \longmapsto \int_{\GL_2(\OO_F)}m(gk Z^{-1})\, dk
    \end{align*}
is  bi-$\GL_2(\OO_F)$-invariant and hence, by the uniqueness of the spherical function, equals $m(Z^{-1})m(g)$. The assertion then follows from the identity (see e.g. \cite[Lemma 6.10]{GodementJacquetBook})
    \begin{align*}
L\left(s+3/2,\pi\right)=\int_{\GL_2(F)}\one_{M_2(\OO_F)}(g)m(g )|\det g|^{s+2}\, dg.
    \end{align*}
\end{proof}
 
\begin{lem} \label{lem:eq:normalizedinitial}
    One has 
    \begin{align*}
    \frac{I(f,\one_{M_2(\OO_F)^2},m_s)}{L(s+3/2,\pi)}
    &=\sum_{Z} m(Z^{-1})|\det(Z)|^{-(s+2)} \int_{D^n}f(Y)\psi(\langle P(Y),Z\rangle )\, dY
\end{align*} 
for $\mathrm{Re}(s)>2$.  Here the sum in $Z$ is over a set of representatives of cosets in
$M_2(F)/M_2(\OO_F)$ such that $Z^{-1} \in M_2(\OO_F).$
\end{lem}
\begin{proof}
By Lemma \ref{lem:HZ} we can rewrite \eqref{eq:simplified} as
\begin{align*}
      \frac{I(f,\one_{M_2(\OO_F)^2},m_s)}{L(s+3/2,\pi)}&=\int_{M_2(F)^n}f(Y)\left(\sum_{Z} m(Z^{-1})|\det(Z)|^{-(s+2)}\psi(\langle P(Y),Z\rangle )\right) \, dY.
\end{align*}
Since $f$ is Schwartz and $|m(Z^{-1})|\le 1,$ by the Fubini-Tonelli theorem to prove the assertion it suffices to show $\sum_{Z} |\det(Z)|^{-(\sigma+2)}$ converges for $\sigma=\mathrm{Re}(s)>2$. 

We define for each $a\in \ZZ_{>0}$ the set
\begin{align*}
    E_{a}:=\{ Z\in M_2(F)/M_2(\OO_F) : \norm{Z}=q^{a}\}.
\end{align*}
Note that $|E_{a}|\le q^{4a}$. Since $|\det(Z)|\ge\norm{Z}$, for $\sigma>2$
\begin{align*}
\sum_{Z} |\det(Z)|^{-(\sigma+2)}\le 1+\sum_{a=1}^\infty |E_a| q^{-a(\sigma+2)}\le 1+\sum_{a=1}^\infty q^{-a(\sigma-2)}<\infty
\end{align*}
\end{proof}

\subsection{Bounding nonabelian Gauss sums}
In view of Lemma \ref{lem:eq:normalizedinitial}, we need to understand the analytic behavior of the integral
\begin{align*}
  \int_{M_2(F)^n}f(Y) \psi(\langle P(Y),Z\rangle )\, dY,
\end{align*}
which is a nonabelian analogue of a Gauss sum.
For this we will use the classical one-dimensional estimates of Lemma~\ref{LEM:classical-gauss-sum-properties}. We continue to assume $\psi$ is unramified.

\begin{lem}\label{lem:unrint}
Suppose $2\nmid q$. Let $Z\in M_2(F)$ with $|\mathrm{tr}(Z)|\ge 1$. Then
\begin{align*}
     \bigg|\int_{M_2(\OO_F)} \psi(\mathrm{tr}(Y^2Z)) \, dY\bigg|&=|\mathrm{tr}(Z)|^{-1/2}\max(|\mathrm{tr}(Z)\det(Z)|, \norm{Z}^2)^{-1/2}\\
     &\le |\mathrm{tr}(Z)|^{-1/2}\norm{Z}^{-1}.
\end{align*}
\end{lem}

\begin{proof} 
 Write $Y=(y_{ij})$ and $Z=(z_{ij})$. Below we identify $M_2(F)$ with $F^4$, viewed as a space of column vectors. Then $\mathrm{tr}(Y^2Z)$ is a quadratic form in $Y$ whose associated matrix is given by
\begin{align*}
    J=\left(\begin{matrix}
    2z_{11} & z_{21} & z_{12} & 0\\
    z_{21} & 0 & z_{11}+z_{22}& z_{21} \\ 
   z_{12} & z_{11}+z_{22}& 0 & z_{12}\\
    0 & z_{21}& z_{12} & 2z_{22}
    \end{matrix}\right).
\end{align*}
That is $\mathrm{tr}(Y^2 Z)=\frac{1}{2}Y^t JY$. Consider the matrix 
\begin{align*}
    R=\left(\begin{matrix}
    0 & 0 & 1 &-z_{22}\\
    1 & 1 & 0 & z_{12}\\
    1 & -1 & 0 & z_{21}\\
    0 & 0 & -1 & -z_{11}
    \end{matrix}\right).
\end{align*}
Let $r:=\mathrm{tr}(Z), r_+:=z_{12}+z_{21}, r_-:=z_{12}-z_{21}.$
One has
\begin{align*}
    R^tJR&=2r\left(\begin{matrix}
    1 & & & \\
    & -1 & & \\
    & & 1 & \\ 
    & & & \det(Z)
    \end{matrix}\right),\\
    \det(R)&=2r, \\
    R^{-1}&=\frac{-1}{2r}\left(\begin{matrix}
   -r_+ & -r & -r&  -r_+ \\
    -r_- & -r & r&  -r_-\\
    -2z_{11}& 0& 0 & 2z_{22}\\ 
    2 & 0& 0& 2
    \end{matrix}\right).
\end{align*}

Since $2\nmid q$, 
\begin{align*}
    R^{-1}\OO_F^4&= \OO_Fe_1+\OO_Fe_2+\OO_Fe_3+r^{-1}\OO_F\begin{psmatrix}
     r_+\\
     r_-\\
      2z_{11}\\
      -2
    \end{psmatrix},
\end{align*}
Changing variables $Y\mapsto RY$, we have
\begin{align*}
    &\int_{M_2(\OO_F)} \psi(\mathrm{tr}(Y^2Z))\, dY\\
    &=|r|\int_{r^{-1}\OO_F}\int_{\OO_F^3} \psi(r(w+r_+z)^2-r(x+r_-z)^2+r(y+2z_{11}z)^2+r\det( Z)4z^2)\, dwdxdydz\\
    &=|r|\int_{r^{-1}\OO_F} \psi(r(r_+^2-r_-^2+(2z_{11})^2+4\det(Z))z^2)\\
    &\times \int_{\OO_F^3} \psi(r(w^2+2r_+wz)-r(x^2+2r_-xz)+r(y^2+4z_{11}yz))\, dwdxdydz
\end{align*}
Since $|r|\ge 1$, by Lemma \ref{LEM:classical-gauss-sum-properties} we have
\begin{align*}
    &\left|\int_{M_2(\OO_F)} \psi(\mathrm{tr}(Y^2Z)) \, dY\right|\\
    &=|r|\left|\int_{r^{-1}\OO_F}\one_{\OO}(r_+z)\one_{\OO}(r_-z)\one_{\OO}(2z_{11}z)\psi(4r\det(Z)z^2)\, dz\right|\left|\int_{\OO_F}\psi(rw^2-rx^2+ry^2) \, dwdxdy\right|\\
    &=|r|^{-\frac{1}{2}}\left|\int_{\varpi^m\OO_F} \psi(4r\det(Z)z^2)\, dz\right|,
\end{align*}
where $m\ge 0$ is the integer such that
\begin{align*}
    q^{-m}=\max(|r|,|r_+|,|r_-|,|2z_{11}|)^{-1}=\norm{Z}^{-1}.
\end{align*}
Changing variables, we have again by Lemma \ref{LEM:classical-gauss-sum-properties}
\begin{align*}
    \left|\int_{\varpi^m\OO_F} \psi(4r\det(Z)z^2)\, dz\right|&=q^{-m}\left|\int_{\OO_F} \psi(4\varpi^{2m}r\det(Z)z^2)\, dz\right|\\
    &=\max(|\mathrm{tr}(Z)\det(Z)|,\norm{Z}^2)^{-1/2}
\end{align*}
and the assertion follows.
\end{proof}

\begin{lem}\label{lem:unrint2}
    Suppose $2|q$. Let $Z\in M_2(F)$ with $|\mathrm{tr}(Z)|\ge |2|$. Then
    \begin{align*}
    &\bigg|\int_{M_2(\OO_F)} \psi(\mathrm{tr}(Y^2Z)) \, dY\bigg|\ll |\mathrm{tr}(Z)|^{-1/2}\norm{Z}^{-1}.
    \end{align*}
\end{lem}

\begin{proof}
    We retain the notation in the previous lemma. In this case 
    \begin{align*}
        R^{-1}\OO_F^4:=\OO_Fe_1+2^{-1}\OO_F\begin{psmatrix}
            1 \\
            1 \\
            0 \\
            0 \\
        \end{psmatrix}+\OO_Fe_3+(2r)^{-1}\OO_F\begin{psmatrix}
     r_+\\
     r_-\\
     2z_{11}\\
      -2
    \end{psmatrix}.
    \end{align*} 
Thus
    \begin{align*}
      & \int_{M_2(\OO_F)} \psi(\mathrm{tr}(Y^2Z))dY\\
&=|r|\int_{r^{-1}\OO_F}\int_{\OO_F^3} \psi(r(w+2^{-1}(x+r_+z))^2-4^{-1}r(x+r_-z)^2+r(y+z_{11}z)^2+r(\det(Z)z^2))\, dwdxdydz\\
&|r|\int_{r^{-1}\OO_F}\psi(r(4^{-1}r_+^2-4^{-1}r_-^2+z_{11}^2+\det (Z))z^2)\\
     &\times \int_{\OO_F^3} \psi(r(w^2+w(x+r_+z)+2^{-1}x(r_+-r_-)z)+r(y^2+2yz_{11}z))\, dwdxdydz.
    \end{align*}

If $|r|< 1$, we can and do assume $\norm{Z}\ge 1$. Then the integral equals
\begin{align*}
    &|r|\int_{r^{-1}\OO_F}\psi(r(4^{-1}r_+^2-4^{-1}r_-^2+z_{11}^2+\det (Z))z^2)\one_{\OO_F}(2rz_{11}z)\\
     &\times \int_{\OO_F} \one_{\OO_F}(r(x+r_+z))\psi(2^{-1}xr(r_+-r_-)z))\, dxdz\\
     &=|r|\int_{r^{-1}\OO_F}\psi(r(4^{-1}r_+^2-4^{-1}r_-^2+z_{11}^2+\det (Z))z^2)\one_{\OO_F}(2rz_{11}z)\one_{\OO_F}(rr_+z)\one_{\OO_F}(2^{-1}r(r_+-r_-)z) \, dz,
\end{align*}
which is dominated by $\norm{Z}^{-1}$. Now suppose $|r|\ge 1.$  By Lemma \ref{LEM:classical-gauss-sum-properties} the integral equals
\begin{align*}
    &|r|\int_{r^{-1}\OO_F}\psi(r(4^{-1}r_+^2-4^{-1}r_-^2+z_{11}^2+\det (Z))z^2)\one_{\OO_F}(2yz_{11})\\
     &\times \int_{\OO_F^3} \psi(r(w^2+w(x+r_+z)+2^{-1}x(r_+-r_-)z)+r(y^2+2yz_{11}z))\one_{\OO_F}(x+r_+z)\, dwdxdydz\\
     &=|r|\int_{r^{-1}\OO_F}\psi(r(-4^{-1}r_+^2+2^{-1}r_+r_- -4^{-1}r_-^2+z_{11}^2+\det (Z))z^2)\one_{\OO_F}(2yz_{11})\one_{\OO_F}(r_+z)\\
     &\times \int_{\OO_F^3} \psi(r(w^2+wx+2^{-1}x(r_+-r_-)z)+r(y^2+2yz_{11}z))\, dwdxdydz
\end{align*}
Observe that since $|r|\ge 1$
\begin{align*}
    &\int_{\OO_F^2}\psi(r(w^2+wx+2^{-1}x(r_+-r_-)z)\, dwdx\\
    &=\int_{\OO_F}\psi(rw^2)\one_{\OO_F}(r(w+2^{-1}(r_+-r_-)z))\, dw\\
    &=\one_{\OO_F}(2^{-1}(r_+-r_-)z)\int_{\OO_F}\psi(rw^2)\one_{\OO_F}(r(w+2^{-1}(r_+-r_-)z))\, dw.
\end{align*}
Thus the original integral equals
\begin{align*}
&|r|\int_{r^{-1}\OO_F}\psi(r(-4^{-1}r_+^2+2^{-1}r_+r_- -4^{-1}r_-^2+z_{11}^2+\det (Z))z^2)\one_{\OO_F}(2yz_{11})\one_{\OO_F}(r_+z)\one_{\OO_F}(2^{-1}(r_+-r_-)z)\\
     &\times \int_{\OO_F^2} \psi(rw^2)\one_{\OO_F}(r(w+2^{-1}(r_+-r_-)z))\psi(r(y^2+2yz_{11}z))\, dwdydz
\end{align*}
The assertion follows by applying Lemma \ref{LEM:classical-gauss-sum-properties} to the integral over $y$ and then applying the trivial bound.
\end{proof}

\subsection{Analytic continuation}
    
\begin{lem}\label{lem:error term in density}
    Suppose $n\ge 4$. For $\sigma\in\RR$, consider the series
    \begin{align*}
       \sum_{Z}|\det(Z)|^{-(\sigma+2)}|\mathrm{tr}(Z)|^{-n/2}\norm{Z}^{-n},
    \end{align*}
where the sum in $Z$ is taken over a set of representatives of cosets in $(M_2(F)-M_2(\OO_F))/M_2(\OO_F)$ such that $\det(Z) \ge \norm{Z}$ and $|\mathrm{tr}(Z)|\ge |2|.$ The series converges when $\sigma>-\frac{n+1}{2}$ and is bounded above by 
        \begin{align*}
            c_Fq^{-(\sigma+\frac{n+1}{2})}\zeta(n-3)\zeta(n/2-1)\zeta\left(\sigma+\frac{n+1}{2}\right),
        \end{align*}
for some positive constant $c_F,$ which is $1$ if $2\nmid q$.
\end{lem}

\begin{proof}
We prove only for the case we need, i.e. $2\nmid q$, and leave the case $2|q$ to the reader. We rewrite the series as
 \begin{align}\label{eq:series}
     \sum_{m=1}^{\infty} q^{-m(\sigma+2)}\sum_{Z: |\det (Z)|=q^m} |\mathrm{tr} (Z)|^{-\frac{n}{2}}\norm{Z}^{-n}.
 \end{align}
 For a fixed $m\ge 1$, we first study the inner sum. We define for integers $0\leq r\leq a$ the set
\begin{align*}
    E_{a,r}:=\{ Z : \norm{Z}=q^{a}, |\mathrm{tr}(Z)|=q^{r}\}.
\end{align*}
Note that $|E_{a,r}|\le q^{3a+r}$.  Since $\norm{Z} \le |\det(Z)|\le \norm{Z}^2,$ the inner sum in \eqref{eq:series} is bounded by
\begin{align*}
    \sum_{\substack{0 \leq a \leq m \leq 2a\\ 0\le r\le a}} q^{3a+r}q^{-\frac{rn}{2}}q^{-na}&\le \zeta(n/2-1)\sum_{0\le \lceil \frac{m}{2}\rceil\le a\le m} q^{(3-n)a}\\&\le q^{-m(n-3)/2}\zeta(n-3)\zeta(n/2-1).
\end{align*}
Thus \eqref{eq:series} is convergent when $\sigma>-\frac{n+1}{2}$ and is bounded above by
\begin{align*}
   \zeta(n-3)\zeta(n/2-1)\sum_{m=1}^\infty q^{-m(\sigma+\frac{n+1}{2})}=  q^{-(\sigma+\frac{n+1}{2})}\zeta(n-3)\zeta(n/2-1)\zeta\left(\sigma+\frac{n+1}{2}\right).
\end{align*}
\end{proof}

\begin{prop}\label{prop:nonarchsplittrivconv}
   Suppose $n\ge 4$ and coefficients of $P$ lie in $\OO_F^\times.$  Then the function
   \begin{align*}
       \frac{I(\one_{M_2(\OO_F)^n},\one_{M_2(\OO_F)^2},m_s)}{L\left(s+\frac{3}{2},\pi\right)}
   \end{align*}
   extends to a holomorphic function in $s$ on $\sigma=\mathrm{Re}(s)>-\frac{n+1}{2}$. Moreover, 
\quash{\begin{align*}
    \left| \frac{I(\one_{M_2(\OO_F)^n},\one_{M_2(\OO_F)^2},m_s)}{L\left(s+\frac{3}{2},\pi\right)}\right|\le  q^{c(\psi)(2(\sigma+2)+3n/2)}|J(s)|\left(1+q^{-(\sigma+\frac{n+1}{2})}\zeta(n-3)\zeta(n/2-1)\zeta(\sigma+\frac{n+1}{2})\right).
\end{align*}}
\begin{align*}
    \left| \frac{I(\one_{M_2(\OO_F)^n},\one_{M_2(\OO_F)^2},m_s)}{L\left(s+\frac{3}{2},\pi\right)}-1\right|\le c_Fq^{-(\sigma+\frac{n+1}{2})}\zeta(n-3)\zeta(n/2-1)\zeta\left(\sigma+\frac{n+1}{2}\right),
\end{align*}
where $c_F>0$ is a constant that is $1$ if $2\nmid q$.
\end{prop}

\begin{proof}
By Lemma \ref{lem:eq:normalizedinitial} it suffices to bound
\begin{align}\label{eq:trivsum2}
    \sum_{Z} |\det(Z)|^{-(s+2)}\int_{M_2(\OO_F)^n}\psi(\langle P(Y),Z\rangle ) \, dY.
\end{align}
Here, using Lemma \ref{lem:reps}, the sum in $Z$ can be taken over a set of representatives of cosets in $(M_2(F)-M_2(\OO_F))/M_2(\OO_F)$ such that  $|\det(Z)| \ge \norm{Z}$ and $|\mathrm{tr}(Z)|\ge |2|.$ By Lemmas \ref{lem:unrint} and \ref{lem:unrint2}, the sum \eqref{eq:trivsum2} is bounded by
\begin{align*}
    C\sum_{Z} |\det(Z)|^{-(\sigma+2)}|\mathrm{tr}(Z)|^{-\frac{n}{2}}\norm{Z}^{-n}.
\end{align*}
for some constant $C>0$ that is $1$ if $2\nmid q$. The assertion then follows from Lemma \ref{lem:error term in density}.
\end{proof}

For the remainder of this section we continue to assume that the coefficients of $P$ lie in $\OO_{F}^\times.$
Recall for the trivial representation, $L(s+3/2,1\circ \nrd)=\zeta_D(s+3/2)=\zeta(s+2)\zeta(s+1).$ Motivated by global considerations, we define the local density to be
\begin{align} \label{loc:dens}
    c\left(\one_{M_2(\OO_F)^n}\right):=\left.\frac{I(\one_{M_2(\OO_F)^n},\one_{M_2(\OO_F)^2},1_s)}{\zeta(s+2)\zeta(s+1)}\right|_{s=-2},
\end{align}
which is well defined for $n\ge 4$ by Proposition \ref{prop:nonarchsplittrivconv}.   This definition is consistent with the definition in Theorem \ref{thm:constant}, as we now explain.
By the proof of Proposition \ref{prop:nonarchsplittrivconv},
    \begin{align*}
        c\left(\one_{M_2(\OO_F)^n}\right) = \sum_{Z} \int_{M_2(\OO_F)^n}\psi(\langle P(Y),Z\rangle ) \, dY
    \end{align*}
is an absolutely convergent series. As the integral on the right  is invariant under $Z \mapsto Z+W$ for $W \in M_2(\OO_F),$ we have
\begin{align}\label{eq:c(f):bis}
        c\left(\one_{M_2(\OO_F)^n}\right) = \int_{M_2(F)} \int_{M_2(\OO_F)^n}\psi(\langle P(Y),Z\rangle ) \, dY.
    \end{align}

\begin{cor}\label{cor:nonvanish split local density}
 Suppose $2 \nmid q$ and $n\ge 5$. Then $c\left(\one_{M_2(\OO_F)^n}\right)\neq 0$.
\end{cor}

\begin{proof}
   By Proposition \ref{prop:nonarchsplittrivconv}, it suffices to check that
   \begin{align*}
       1> q^{-\frac{n-3}{2}}\zeta(n-3)\zeta(n/2-1)\zeta\left(\frac{n-3}{2}\right)
   \end{align*}
   when $q=3$ and $n=5,$ which is straightforward.
\end{proof}

We can alternatively write
\begin{align*}
    c\left(\one_{M_2(\OO_F)^n}\right)&=\int_{M_2(F)}\left(\int_{M_2(\OO_F)^n}\psi(\langle P(Y),Z\rangle ) dY \right)\, dZ\\
    &=\lim_{m\to\infty} \int_{M_2(F)} \one_{\varpi^{-m}M_2(\OO_F)}(Z)\left(\int_{M_2(\OO_F)^n}\psi(\langle P(Y),Z\rangle ) dY \right)\, dZ\\
    &=\lim_{m\to\infty} q^{4m}\int_{M_2(\OO_F)^n}\one_{\varpi^{m}M_2(\OO_F)}(P(Y))\, dY\\
    &=\lim_{m\to\infty} q^{4m(1-n)}\sum_{Y\in (M_2(\OO_F)/\varpi^{m}M_2(\OO_F))^n}\one_{\varpi^{m}M_2(\OO_F)}(P(Y))>0
\end{align*}
when $2\nmid q$ and $n\ge 5$.

\begin{cor}\label{cor:split local density formula}
    For $2\nmid q$ and $n\ge 5,$ 
    \begin{align*}
       \#\{Y\in (M_2(
        \OO_F
        )/\varpi^{m} M_2(\OO_F))^n:P(Y)\equiv 0\mod{\varpi^{m}}\}\sim c\left(\one_{M_2(\OO_F)^n}\right)q^{4m(n-1)}.
    \end{align*}
   \qed
\end{cor}

\section{Spectral bounds for nonsplit cases} \label{sec:arch:bounds}

Let $F$ be any local field of characteristic zero and assume $D$ is nonsplit (hence a division algebra). Then we have an exact sequence
\begin{align*}
    1\lto (D^\times)^1\lto D^\times\lto \nrd(D^\times)\lto 1,
\end{align*}
where $(D^\times)^1$ is the subgroup of norm $1,$ which is also the maximal compact subgroup of $D^\times$. 
\begin{lem} \label{lem:unram}
    Suppose $\Phi$ is (bi-)$(D^\times)^1$-invariant in both entries. Then $I(f,\Phi,m_s)=0$ unless $\pi$ is unramified.
\end{lem}

\begin{proof}
For $\mathrm{Re}(s)>0,$ we have for $k_1,k_2\in (D^\times)^1$
\begin{align*}
    I(f,\Phi,m_s)&=\int_{D^\times}\int_{D^n} f(Y)\Phi\left(P(Y)k_2^{-1}g^{-1}k_1^{-1},k_1gk_2\right)m_s(k_1gk_2)\, dY dg\\
    &=\int_{D^\times}\int_{D^n} f(Y)\Phi\left(P(Y)g^{-1}gk_2^{-1}g^{-1},g\right)m_s(k_1gk_2)\, dY dg.
\end{align*}
Since $(D^\times)^1$ is a normal subgroup of $D^\times,$ by our assumptions on $\Phi$ the integral equals
\begin{align*}
    \int_{D^\times}\int_{D^n} f(Y)\Phi\left(P(Y)g^{-1},g\right)m_s(k_1gk_2)\, dY dg.
\end{align*}
\end{proof}

As in the split case we assume $\Phi$ is (bi-)$(D^\times)^1$-invariant in both entries, so $I(f,\Phi,m_s)$ is zero unless $\pi$ is unramified by Lemma \ref{lem:unram}.
The only unitary unramified representations $\pi$ of $D^\times$ are characters $|\mathrm{nrd}|^{it}$ where $t \in \RR.$ Thus for our purposes in this section it is no loss of generality to assume that $\pi$ is the trivial representation.

In view of this, we henceforth assume $\pi$ is the trivial representation and the matrix coefficient $m$ is identically $1.$ For $\mathrm{Re}(s)>0$, we have by Fourier inversion
\begin{align*}
I(f,\Phi,1_s)=\int_D\left(\int_{D^n}f(Y)\psi(\langle P(Y),Z\rangle )\, dY\right)\left(\int_{D^\times } \mathcal{F}_1(\Phi)(-gZ,g)|\nrd(g)|^{s+2}\, dg\right)\, dZ.
\end{align*}

We will follow the strategy employed in the split case to prove the analytic continuation of $I(f,\Phi,1_s).$ Recall the zeta integral
\begin{align*}
    \mathcal{Z}_{Z}(\Phi,s):=\mathcal{Z}_{Z}(\Phi,s,1)=\int_{D^\times } \mathcal{F}_1(\Phi)(-gZ,g)|\nrd(g)|^{s+2}\, dg,
\end{align*}
and 
\begin{align*}
    H_{Z}(\Phi,s):=H_{Z}(\Phi,s,1)=\frac{\mathcal{Z}_{Z}(\Phi,s,1)}{\zeta_D(s+3/2)}.
\end{align*}
Thus
\begin{align}\label{eq:Inonsplit}
    \frac{I(f,\Phi,1_s)}{\zeta_D(s+3/2)}=\int_{D} \left(\int_{D^n} f(Y)\psi\left(\left\langle P(Y),Z\right\rangle\right)\, dY\right) H_Z(\Phi,s)\, dZ.
\end{align}

\subsection{The non-Archimedean case}

Suppose $F$ is non-Archimedean and $\psi$ is unramified. We begin by recalling standard facts on quaternion  division algebras to fix notations. We refer to \cite[\S 6.1.5, \S 13.3.10]{voight2021quaternion} for details. 

 Let $\varpi_D\in\OO_D$ be a uniformizer in $D$, chosen so that $\varpi_D^2=\varpi$. Then $\varpi_D\OO_D$ is a two-sided ideal in $\OO_D$ and $\OO_D/\varpi_D\OO_D$ is a quadratic extension over $\mathbb{F}_q$. Let $K$ be the unique quadratic unramified extension of $F$. We may identify $D$ with $K\oplus K\varpi_D$ equipped with the multiplication rule
\begin{align*}
    \varpi_D\alpha=\bar{\alpha}\varpi_D
\end{align*}
for $\alpha\in K$, where $\bar{\,\cdot\,}$ is the Galois conjugation for $K/F$. Under this identification, the maximal order $\OO_D$ can be identified with $\OO_K\oplus\OO_K\varpi_D$, where $\OO_K$ is the ring of integers of $K$.

Let $u\in\OO_F^\times-(\OO_F^\times)^2$ so that $K=F[x]/(x^2-u)$. We may write $K=F\oplus F\sqrt{u}$, and 
\begin{align*}
    D=K\oplus K\varpi_D=F\oplus F\sqrt{u}\oplus F\varpi_D\oplus F\sqrt{u}\varpi_D,
\end{align*}
with the multiplication rule $\varpi_D\sqrt{u}=-\sqrt{u}\varpi_D$. We will often identify $D$ with $F^4$ as a vector space via the $F$-linear isomorphism
\begin{align*}
    F^4&\tilde{\lto} D\\
    (z_1,z_2,z_3,z_4)&\longmapsto z_1+z_2\sqrt{u}+z_3\varpi_D+z_4\sqrt{u}\varpi_D.
\end{align*}
and define a vector norm on $D$ by
\begin{align*}
    \|Z\|:=\max(|z_1|,|z_2|,|z_3|,|z_4|).
\end{align*}
Note that $\trd(Z)=2z_1, \nrd(Z)=z_1^2-z_2^2u-\varpi z_3^2+\varpi u z_4^2$. Therefore, if $Z \neq 0$ then
\begin{align}\label{eq:nrd=norm}
    1 \ll \frac{|\nrd(Z)|^{\frac{1}{2}}}{\norm{Z}} \le 1.
\end{align}
Indeed, the right inequality is obvious.  To prove the left inequality we may assume that $Z \in \OO_D$ and some $z_i \in \OO_F^\times;$   hence $\norm{Z}=1.$  The set of all $Z \in \OO_D$ such that $z_i \in \OO_F^\times$ is a compact subset of $D.$  Since $\mathrm{nrd}$ is anisotropic, $|\nrd|$ does not vanish on this compact set and hence is bounded below.

\begin{lem}\label{lem :nonsplitcoeff}
Let $\ell,e\in \ZZ$. Let $Z\in D$ with $|\trd(Z)|\ge q^{e}$. Then
\begin{align*}
    \left|\int_{\varpi_D^\ell\OO_D} \psi(\trd(Y^2 Z)) \, dY\right|\ll_{e,\ell} |\trd(Z)|^{-\frac{1}{2}}\norm{Z}^{-1}.
\end{align*}
\end{lem}
\begin{proof}
By changing variables $Y\mapsto \varpi^{-2\lfloor\frac{\ell}{2}\rfloor}Y$ and replacing $Z$ and $e$, we can and do assume $\ell\in\{0,1\}.$ Write $Y=(y_i)$ and $Z=(z_i)$. Consider first $\ell=0$. We have
\begin{align*}
    \mathrm{trd}(Y^2Z)&=2z_1y_1^2+2uz_1y_2^2+2\varpi z_1y_3^2-2 u\varpi z_1y_4^2 +4uz_2y_1y_2+4\varpi z_3y_1y_3-4u\varpi z_4y_1y_4.
\end{align*}
Thus, up to a positive constant $\int_{\OO_D} \psi(\trd(Y^2 Z))\, dY$ is
\begin{align*}
    &\int_{\OO_F^4}\psi(2(z_1y_1^2+uz_1y_2^2+\varpi z_1y_3^2- u\varpi z_1y_4^2 +2uz_2y_1y_2+2\varpi z_3y_1y_3-2u\varpi z_4y_1y_4))\, dy_1dy_2dy_3dy_4\\
    &=\int_{\OO_F}\psi(2z_1y_1^2 )\left(\int_{\OO_F}\psi(2uz_1(y_2^2+2z_1^{-1}z_2y_1y_2))\, dy_2\right)\\
    &\times\left(\int_{\OO_F}\psi(2\varpi z_1( y_3^2+2z_1^{-1}z_3y_1y_3))\, dy_3\right)\left(\int_{\OO_F}\psi(-2u\varpi z_1(y_4^2 +2z_1^{-1}z_4y_1y_4))\, dy_4\right)\, dy_1.
\end{align*}
By Lemma \ref{LEM:classical-gauss-sum-properties}, for $2\le i\le 4$ the integral over $y_i$ vanishes unless $|z_1^{-1}z_iy_1|\ll \min(1,|z_1|)^{-1}.$  When $z_i$ is nonzero this bound is equivalent to $|y_1|\ll |z_i|^{-1}\max(|z_1|,1)$. By applying \eqref{eq:gaussbd} to the integrals over $y_i$ for $2\le i\le 4$ and the trivial bound to the integral over $y_1$, one has
\begin{align*}
    \left|\int_{\OO_D} \psi(\trd(Y^2 Z)) \, dY\right|\ll |z_1|^{-\frac{3}{2}} \min(1,\max(|z_2|,|z_3|,|z_4|)^{-1}\max(|z_1|,1))
\end{align*}
Here by convention if $z_2=z_3=z_4=0$ then $\min(1,\max(|z_2|,|z_3|,|z_4|)^{-1}\max(|z_1|,1))=1.$
One easily checks this is dominated by $|\trd(Z)|^{-\frac{1}{2}}\norm{Z}^{-1},$ and thus the assertion follows.

For $\ell=1$, by changing variables $Y\mapsto \varpi_D Y$, it suffices to bound the integral
\begin{align*}
\int_{\OO_D}\psi(\trd((\varpi_DY)^2Z))\, dY.
\end{align*}
By a direct calculation, we have
\begin{align*}
    \trd((\varpi_DY)^2Z)&=2\varpi\left( z_1y_1^2 -uz_1y_2^2+\varpi z_1y_3^2+u\varpi z_1y_4^2+2\varpi z_3y_1y_3+2u\varpi z_4y_2y_3-2u\varpi z_2y_3y_4  \right)
\end{align*}
The assertion can be proved similarly as in the case $\ell=0$ and we leave it to the reader.
\end{proof}

At this point it is helpful to recall that
\begin{align}
\int_{\OO_D \cap D^\times}|\mathrm{nrd}(g)|^{s}d^\times g=d^\times g(\OO_D^\times)\zeta(s).
\end{align}
We refer to \cite[Lemma 29.7.17]{voight2021quaternion} for the proof.

\begin{lem}\label{lem:HZnonarchnonsplit}
    We have
    \begin{align*}
        |H_Z(\Phi,s)|\ll \max(|\nrd(Z)|,1)^{-(\mathrm{Re}(s)+2)}J(q^{-\mathrm{Re}(s)}).
    \end{align*}
    for some polynomial $J\in \RR[s]$ with positive coefficients that is independent of $Z$.
\end{lem}

\begin{proof}
    Assume first $|\nrd(Z)|\le 1$. Then there exists
    a finite set $A\subset \mathcal{S}(D)$ independent of $Z$ such that the function $g\mapsto \mathcal{F}_1(\Phi)(-gZ,g)$ can be written as $\sum_{\widetilde{\Phi}\in A_Z} \widetilde{\Phi}$ for some subset $A_Z\subset A$. Thus there exists a finite set $A'\subset \CC[q^{-s}]$ such that $H_Z(\Phi,s)=\sum_{J\in A'_Z} J(s)$ for some subset $A'_Z\subset A'$. The assertion then follows easily. For $|\nrd(Z)|\ge 1$, a similar argument applies after changing variables $g\mapsto gZ^{-1}$.
\end{proof}

\begin{prop}\label{prop:nonarchnonsplittrivconv}
   Suppose $n\ge 3$. Then the function
   \begin{align*}
       \frac{I(\one_{\OO_D^n},\Phi, 1_s)}{\zeta_D\left(s+3/2\right)}
   \end{align*}
   extends to a holomorphic function in $s$ on $\sigma=\mathrm{Re}(s)>-\frac{n+1}{2}$.
\end{prop}

\begin{proof}
By Lemma \ref{lem:HZnonarchnonsplit}, it suffices to study the convergence of
\begin{align*}
    &\int_{D}\max(|\nrd(Z)|,1)^{-(\sigma+2)}\left|\int_{\OO_D^n} \psi(\langle P(Y),Z\rangle)\, dY\right| \, dZ\\
    &=\vol(\OO_D)\left(1+\sum_{Z\in(D-\OO_D)/\OO_D} |\nrd(Z)|^{-(\sigma+2)}\left|\int_{\OO_D^n} \psi(\langle P(Y),Z\rangle)\, dY\right| \right)
\end{align*}
For each coset in $(D-\OO_D)/\OO_D,$ we choose a representative $Z$ with $|z_{11}|\ge 1$ so $|\trd(Z)|\ge |2|$. By Lemma \ref{lem :nonsplitcoeff} and \eqref{eq:nrd=norm}, it suffices to study the convergence of the series
\begin{align*}
    \sum_{Z} |\mathrm{trd}(Z)|^{-\frac{n}{2}}\norm{Z}^{-(2\sigma+4+n)}.
\end{align*}
It equals
\begin{align*}
    \sum_{a=1}^\infty q^{-a(2\sigma+4+n)} \sum_{Z:\norm{Z}=q^a} |\mathrm{trd}(Z)|^{-\frac{n}{2}}&\ll \sum_{a=1}^\infty q^{-a(2\sigma+4+n)} \sum_{r=0}^{a} q^{3a+r}q^{-r\frac{n}{2}}\\
    &\le \zeta(n/2-1)\sum_{a=1}^\infty q^{-a(2\sigma+1+n)} 
\end{align*}
which converges when $\sigma>-\frac{n+1}{2}$.
\end{proof}

Recall that $\zeta_D(s+3/2)=\zeta(s+2)$. Define the local density to be
\begin{align} \label{loc:dens2}
c\left(\one_{\OO_D^n}\right):=\frac{1}{dg(\OO_D^\times)\mathrm{vol}(\OO_D)}\frac{I(\one_{\OO_D^n},\one_{\OO_D^2},1_s)}{\zeta(s+2)}\Bigg|_{s=-2}. 
\end{align}
Since $\OO_D^\#=\varpi^{-1}_D \OO_D,$ we have
\begin{align*}
    \mathcal{Z}_{Z}(\one_{\OO_D^2},s)&=\vol(\OO_D)\int_{D^\times\cap \OO_D} \one_{\OO_D}(\varpi_DgZ)|\nrd(g)|^{s+2}\,  dg\\
    &=dg(\OO_D^\times)\vol(\OO_D)\max(|\nrd(\varpi_DZ)|,1)^{-(s+2)}\zeta(s+2).
\end{align*}
Therefore, 
\begin{align}\label{eq:nonsplitlocaldensity}
    c\left(\one_{\OO_D^n}\right)=\int_{D}\left(\int_{\OO_D^n} \psi(\langle P(Y),Z\rangle)\, dY\right)\, dZ.
\end{align}
Thus the definition \eqref{loc:dens2} is consistent with the definition in Theorem \ref{thm:constant}.

\begin{lem}\label{lem: nonsplit local density}
    Suppose $n\ge 5$ and coefficients of $P$ lie in $\OO_F^\times$. Then $c\left(\one_{\OO_D^n}\right)>0$.
\end{lem}

\begin{proof}
     We have that
    \begin{align*}
        &\int_{D}\left(\int_{\OO_D^n} \psi(\langle P(Y),Z\rangle)\, dY\right)\, dZ\\
        &=\lim_{m\to \infty}\int_{D} \one_{\varpi^{-m}\OO_D}(Z)\left(\int_{\OO_D^n} \psi(\langle P(Y),Z\rangle)\, dY\right)\, dZ\\
        &=\vol(\OO_D)\lim_{m\to \infty} q^{4m}\int_{\OO_D^n} \one_{\varpi_D^{2m-1}\OO_D}(P(Y))\, dY
    \end{align*}
    is nonnegative. To prove the assertion, it suffices to show
    \begin{align}\label{eq:nonarchasym}
        \int_{\OO_D^n} \one_{\varpi_D^{2m-1}\OO_D}(P(Y))\, dY\gg q^{-4m}.
    \end{align}

Since $n\ge 5,$ by \cite[Chapter VI Theorem 2.12]{Lam:Quad} there exists nonzero $Y_0=(y_1,\ldots, y_n)\in \OO_F^n$ such that $P(Y_0)=0$. We can and do assume $y_1,y_2\in \OO_F^\times$. By the discussion above Theorem 4.5 in \cite{Conrad_adelic_points}, the map $\varphi:D^\times\to D^\times$ given by $\varphi(Y):=Y^2$ is open.  Thus there exists $N\in \ZZ_{\ge 1}$ such that $1+\varpi^{2N}\OO_D\le \varphi(\OO_D^\times).$ Let $U:=\varphi^{-1}(1+\varpi^{2N}\OO_D).$ The integral on the left in \eqref{eq:nonarchasym} is bounded below by
    \begin{align} \label{bbl}
         \int_{y_1U}\left(\int_{y_2U\times \prod_{i=3}^n(y_i+\varpi^{2N}\OO_D)} \one_{\varpi_D^{2m-1}\OO_D}(P(Y))\, dY_n\ldots dY_{2} \right)\, dY_1.
    \end{align}
Since
\begin{align*}
    \{\coeff_2Y_2^2+\ldots+\coeff_nY_n^2 : Y_2\in y_2U, Y_i\in y_i+\varpi^{2N}\OO_D, i\ge 3\}= -\coeff_1y_1^2+\varpi^{2N}\OO_D,
\end{align*}
we have \eqref{bbl} is  $\gg_{q,n} q^{-4m}$ for $m>2N+1$. 

\end{proof}

\subsection{The Archimedean case}
In this subsection $F \cong \RR$ and $D$ is the Hamiltonian quaternions. Note that $D^\times\cong \RR_{>0}\times \mathrm{SU}(2)\cong \RR_{>0}\times S^3$. Let $\varepsilon>0.$ For $f$ in a Schwartz space we use the symbol $\ll_{f,\varepsilon}$ to indicate that the implicit constant depends on $\varepsilon$ and for each $\varepsilon$ it is continuous as a function of $f$ with respect to the usual Fr\'echet topology on the Schwartz space.

\begin{lem}\label{lem:Archimedeancoeff}
    For $f\in \mathcal{S}(D)$ and $Z\in D^\times$ with $\nrd(Z)\ge 1$, we have
    \begin{align*}
        \left|\int_{D} f(Y)\psi(\trd(Y^2Z))\, dY\right|\ll_{f,\varepsilon} \max(|\trd(Z)|,1)^{-1} \nrd(Z)^{\varepsilon-\frac{1}{2}}.
    \end{align*}
\end{lem}

\begin{proof}
 Write $Y=y_1+y_2i+y_3j+y_4k$ and $Z=z_1+z_2i+z_3j+z_4k$.  We use these coordinates to identify $D$ with $\RR^4.$ Then $\mathrm{trd}(Z)=2z_1$ and $\mathrm{nrd}(Z)=z_1^2+z_2^2+z_3^2+z_4^2.$  We have
\begin{align*}
    2^{-1}\mathrm{trd}(Y^2Z)= (y_1^2-y_2^2-y_3^2-y_4^2)z_1-2y_1y_2z_2-2y_1y_3z_3-2y_1y_4z_4.
\end{align*}
After replacing $\psi(x)$ by $\psi(2x)$ the integral in the lemma is
\begin{align}\label{eq:archoriginal}
    \int_{\RR^4} f(y_1,y_2,y_3,y_4)\psi((y_1^2-y_2^2-y_3^2-y_4^2)z_1-2y_1y_2z_2-2y_1y_3z_3-2y_1y_4z_4) \, dy_1dy_2dy_3dy_4.
\end{align}

Assume first that $|z_1| \leq \tfrac{1}{2}.$  Choose $2 \leq i \leq 4$ such that $\norm{Z}=|z_i|.$  We then repeatedly apply integration by parts to the expression
\begin{align*}
    \int_{\RR^4} \left(f(y_1,y_2,y_3,y_4)\psi((y_1^2-y_2^2-y_3^2-y_4^2)z_1)\right)\psi(-2y_1y_2z_2-2y_1y_3z_3-2y_1y_4z_4) \, dy_1dy_2dy_3dy_4
\end{align*}
where we differentiate the expression in parenthesis with respect to $y_i.$  For any 
$N \in \ZZ_{\geq 1}$ the contribution of $|y_1| \geq |z_i|^{\varepsilon-1}$ is $O_{f,N}(|z_i|^{-N\epsilon}).$  Bounding the integral over $|y_1|<|z_i|^{\varepsilon-1}$ trivially we obtain the lemma in this case when $|z_1| \leq \frac{1}{2}.$ 

Now assume that $|z_1| \geq \tfrac{1}{2}.$  In this case we apply the method of stationary phase.  The quadratic form $(y_1^2-y_2^2-y_3^2-y_4^2)z_1-2y_1y_2z_2-2y_1y_3z_3-2y_1y_4z_4$ in $(y_1,y_2,y_3,y_4)$  is nondegenerate. Thus the gradient vanishes only at $(y_1,y_2,y_3,y_4)=0.$  The Hessian is 
\begin{align}
  \det  \begin{psmatrix} 2z_1 & -2z_2 & -2z_3 & -2z_2\\ -2z_2 & -2z_1& 0 & 0 \\ -2z_3 & 0 & -2z_1 & 0\\ -2z_4 & 0 & 0 & -2z_1\end{psmatrix}=-8z_1^2(\mathrm{nrd}(Z))
\end{align}
Thus the standard stationary phase bound \cite[Lemma 7.7.3]{Hormander:PDI}
suffices in this case.
\end{proof}

\begin{lem}\label{lem:Archzeta}
 The function $H_Z(\Phi,s)$ is entire. Moreover, given any real numbers $\sigma_1<\sigma_2,$ $0<\epsilon,$ and polynomial $J\in \CC[s],$ one has
 \begin{align*}
     |J(s)\mathcal{Z}_Z(\Phi, s)|\ll_{\sigma_1,\sigma_2, J,\epsilon,\Phi} \max(\nrd(Z),1)^{-(\mathrm{Re}(s)+2)} 
 \end{align*} 
 for $\sigma_1<\mathrm{Re}(s)<\sigma_2$ and $\epsilon<\mathrm{Im}(s).$
\end{lem}

\begin{proof}
    Note that $\mathcal{F}_1(\Phi)$ is invariant under $S^3$ in both entries as $\Phi$ is. Using spherical coordinates, $\mathcal{Z}_{Z}(\Phi,s)$ is up to a positive constant
\begin{align}\label{eq:spherical}
    \int_{0}^\infty \mathcal{F}_1(\Phi)(ar,r)r^{2s+3}\, dr.
\end{align}
where $a:=\nrd(Z)^{\frac{1}{2}}$. By Tate's thesis, $\Gamma(s+2)^{-1}\mathcal{Z}_{Z}(\Phi,s)$ is entire and the first assertion follows. The second assertion is a consequence of \cite[Chapter 1 Theorem 3.1]{Igusa:forms} if $a\le 1.$ For $a>1,$ the same argument applies after a change of variables $r\mapsto a^{-1}r$.
\end{proof}

\begin{prop}\label{prop:archconv}
Suppose $n\ge 3$. For any $f\in \mathcal{S}(D^n),$ the function
\begin{align*}
    \frac{I(f,\Phi,1_s)}{\zeta_D(s+3/2)}
\end{align*}
extends to a holomorphic function on $\sigma=\mathrm{Re}(s)>-\frac{n+1}{2}$. Moreover, $I(f,\Phi,1_s)$ is rapidly decreasing on vertical strips (away from its  poles).
\end{prop}

\begin{proof}
By \eqref{eq:Inonsplit} and Lemmas \ref{lem:Archimedeancoeff} and \ref{lem:Archzeta}, it suffices to prove that
\begin{align*}
    \int_{\nrd(Z)\le 1} dZ+
    \int_{\nrd(Z)>1} \max(|\trd(Z)|,1)^{-n}\nrd(Z)^{-(\sigma+2+\frac{n}{2}-\epsilon)}\, dZ.
\end{align*}
converges for $\sigma>-\frac{n+1}{2}+\epsilon$. Using spherical coordinates, the integral above is equal to a positive constant times
    \begin{align*}
        &\int_{1}^\infty\int_{0}^\pi \int_{0}^\pi\int_{0}^{2\pi}\max(2r|\cos(\theta_1)|,1)^{-n}r^{-(n+1+2\sigma-2\epsilon)}\sin(\theta_1)^2\sin(\theta_2)\, d\theta_3d\theta_2d\theta_1dr\\
        &\ll \int_{1}^\infty \int_{\cos^{-1}(r^{-1})}^{\pi/2} r^{-(n+1+2\sigma-2\epsilon)}\sin(\theta)^2\, d\theta dr\\
        &+\int_{1}^\infty \int_0^{\cos^{-1}(r^{-1})}r^{-(2n+1+2\sigma-2\epsilon)}(\sec(\theta)^n-\sec(\theta)^{n-2})\, d\theta dr.
    \end{align*}
    We have 
    \begin{align*}
        \int_{\cos^{-1}(r^{-1})}^{\pi/2} \sin^2(\theta)\, d\theta=\frac{\theta-\sin(\theta)\cos(\theta)}{2}\bigg|_{\cos^{-1}(r^{-1})}^{\pi/2}\ll r^{-1}
    \end{align*}
    and by the reduction formula for antiderivaties of $\sec^n$
    \begin{align*}
        \int_0^{\cos^{-1}(r^{-1})}\left(\sec(\theta)^n-\sec(\theta)^{n-2}\right)\, d\theta\ll r^{n-1}.
    \end{align*}
    Therefore, we are to examine the convergence of
    \begin{align*}
        \int_{1}^\infty r^{-(n+2+2\sigma-2\epsilon)}\, dr.
    \end{align*}
    It converges when $n+2+2\sigma-2\epsilon>1$, i.e., $\sigma>-\frac{n+1}{2}+\epsilon$. 
\end{proof}

\begin{lem}\label{lem:Archop}
    Suppose $n\ge 4$. The linear operator on $\mathcal{S}(D^n)$ given by
    \begin{align*}
        f\mapsto c(f):=\int_{D} \left(\int_{D^n} f(Y)\psi\left(\left\langle P(Y),Z\right\rangle\right)\, dY\right)\, dZ
    \end{align*}
    is nonzero.
\end{lem}

\begin{proof}
Since we are working over $F\cong\RR,$ upon changing variables if necessary we may assume  $P(Y)=Y_1^2+\ldots+Y_n^2.$ Consider the closed $F$-subscheme $U$ of $\mathbb{A}^{4n}$ cut out by the system of equations \eqref{EQN:explicit-system}. The Jacobian of \eqref{EQN:explicit-system} has full rank if $w_i\neq 0$ for some $i.$
Therefore, the $F$-points of the smooth locus $U^{\mathrm{sm}}(F)$ is an open dense subset of $U(F)$. If $f\in C^\infty_c(D^n)$ is a function whose support does not intersect the subset
$$
\{(w,x,y,z) \in (F^4)^n:w=0\}
$$
then $f|_{U(F)}$ is compactly supported on $U^{\mathrm{sm}}(F)$ and
\begin{align*}
    \int_{D} \left(\int_{D^n} f(Y)\psi\left(\left\langle P(Y),Z\right\rangle\right)\, dY\right) \, dZ=\int_{U^{\mathrm{sm}}(F)} f(u)\, du
\end{align*}
for an appropriate measure $du$ on $U^{\mathrm{sm}}(F)$ (see \cite[\S 7]{Getz:Hsu:Leslie}). Therefore, the operator is easily seen to be nonzero.
\end{proof}

\begin{lem}\label{lem:Archnonvanish}
    Suppose $n\ge 4$. Let $c_1,c_0\in \CC^\times$ be given. There exist $(f,\Phi) \in C^\infty_c(D^n) \times C_c^\infty(D^2)$ such that $\Phi$ is bi-invariant under $S^3$ in both entries, satisfies $\mathcal{F}_2(\Phi)(0,0)\neq 0$ and $\Phi(t,0)=0$ for all $t\in D,$ and
    \begin{align*}
        \frac{c_1 I^1(f,\Phi,1_s)-c_0 I^0(f,\Phi,1_s)}{\zeta_D(s+3/2)}\Bigg|_{s=-2}\neq 0.
    \end{align*}
\end{lem}
\begin{proof}
    It suffices to show there exists $f\in C^\infty_c(D^n)$ such that the $\CC$-span of   
    \begin{align*}
       \zeta_D(s+3/2)^{-1}(I^1(f,\Phi,1_s), I^0(f,\Phi,1_s))\Bigg|_{s=-2}
    \end{align*}
    with $\Phi\in C^\infty_c(D^2)$ satisfying the stated conditions is $\CC^2$. Furthermore, since both coordinates are continuous in $\Phi,$ it suffices to show that  $\CC$-span as $\Phi$ ranges over all functions in $\mathcal{S}(D^2)$ satisfying the stated conditions is $\CC^2.$

    For $a\in \RR_{>0},$ consider $\Psi_a\in \mathcal{S}(D)$ given by $\Psi_a(X)=e^{-a\nrd(X)}$. One has
    \begin{align*}
        \int_{D^\times} \Psi_a(g)\nrd(g)^{s} dg=\vol((D^\times)^1)\int_{0}^\infty e^{-ar^2}|r|^{2s}d^\times (r^2)=\frac{\vol((D^\times)^1)}{2}a^{-s}(2\pi)^s\zeta_D(s-1/2)
    \end{align*}
     For $a_1,a_2\in \RR_{>0}$, let $\Phi_{a_1,a_2}:=\Psi_{a_1}\otimes \Psi_{a_2}$. Then one has
    \begin{align*}
        \mathcal{Z}_{Z}(\Phi_{a_1,a_2},s)=\mathcal{F}_1(\Phi_{a_1,a_2})(0,0)\frac{\vol((D^\times)^1)}{2}Q(\nrd(Z))^{-s-2}\zeta_D(s+3/2)
    \end{align*}
    for some linear function $Q$. 
    Thus
    \begin{align}\label{eq:I1specific}
            \frac{I^1(f,\Phi_{a_1,a_2},1_s)}{\zeta_D(s+3/2)}\Bigg|_{s=-2}=\mathcal{F}_1(\Phi_{a_1,a_2})(0,0)\frac{\vol((D^\times)^1)}{2}\int_{D} \left(\int_{D^n} f(Y)\psi\left(\left\langle P(Y),Z\right\rangle\right)\, dY\right)\, dZ.
    \end{align}
    Similarly, 
    \begin{align}\label{eq:I0specific}
         \frac{I^0(f,\Phi_{a_1,a_2},1_s)}{\zeta_D(s+3/2)}\Bigg|_{s=-2}=\mathcal{F}_2(\Phi_{a_1,a_2})(0,0)\frac{\vol((D^\times)^1)}{2}\int_{D} \left(\int_{D^n} f(Y)\psi\left(\left\langle P(Y),Z\right\rangle\right)\, dY\right)\, dZ. 
    \end{align}

        By Lemma \ref{lem:Archop}, we can choose $f$ so that the integral in \eqref{eq:I1specific} and \eqref{eq:I0specific} does not vanish. Then 
        \begin{align*}
             \zeta_D(s+3/2)^{-1}(I^1(f,\Phi_{2,2}-\Phi_{2,1},1_s), I^0(f,\Phi_{2,2}-\Phi_{2,1},1_s))\Bigg|_{s=-2}
        \end{align*}
        and 
        \begin{align*}
           \zeta_D(s+3/2)^{-1}(I^1(f,\Phi_{2,2}-\Phi_{1,1},1_s), I^0(f,\Phi_{2,2}-\Phi_{1,1},1_s))\Bigg|_{s=-2}
        \end{align*}
        are two linear independent vectors in $\CC^2$. This completes the proof.
\end{proof}

\section{General Hessian estimates}
\label{SEC:general-hessian-based-analysis}

In this section, we assume $P(Y)=\sum_{1\le i\le n} \coeff_i Y_i^2$, where $\coeff_i\in F^\times$ are fixed.  We fix a place $v$ of $F$ and omit it from notation, writing $F:=F_v,$ $D:=D_{F_v}$, etc.
For any integer $N\ge 0$, we fix vector norms $\norm{\cdot}$ on $F^N$ and $M_N(F)\cong F^{N^2}$,
by the formula \eqref{box-vector-norm}.

We need the following quadric fibration result.
We emphasize that $D-\{0\}$ may be strictly larger than $D^\times$,
since $D$ has nontrivial zero divisors if $D$ is split.

\begin{lem}
\label{LEM:general-rank-bounds}
Fix a choice of coordinates on $D$,
i.e.~an $F$-linear isomorphism $D \cong F^4$.
Let
\begin{align*}
    m:M_{4n}(F) &\lto F^N\\
    A &\longmapsto (m_1(A),\dots,m_N(A))
\end{align*}
be the map sending $A\in M_{4n}(F)$ to its $N:=\binom{4n}{2n}^2$ size $2n\times 2n$ minors.
Let $W\in D$.
\begin{enumerate}
\item
If $W\ne 0$, then $\trd(W Y^2)$ is a quadratic form over $F$ of Hessian rank $\ge 2$.

\item
If $W\ne 0$, then $\trd(W P(Y))$ is a quadratic form over $F$ of Hessian rank $\ge 2n$.

\item Let $H_W\in M_{4n}(F)$ be the Hessian matrix of $\trd(W P(Y))$.
The map $$D \to F^N,\quad W\mapsto m(H_W)$$
sends any compact subset of $D-\{0\}$ to a compact subset of $F^N - \{0\}$.

\item
Let $H,A_1,A_2\in M_{4n}(F)$.
Suppose $\norm{H} \ll 1$ and $\norm{m(H)} \gg 1$.
If $\norm{A_i}\ll 1$ and $\abs{\det(A_i)} \asymp 1$,
then $\norm{A_1HA_2} \ll 1$
and $\norm{m(A_1HA_2)} \gg 1$.
\end{enumerate}
\end{lem}

\begin{proof}

We first remark that if $Q$ is a nonzero quadratic form on $F^r$ and $L$ is a nonzero linear form on $F^r$, then
\begin{equation*}
\rank(Q\vert_{L=0}) - \rank(Q) \in \{0,-1,-2\}.
\end{equation*}
Indeed, after linearly changing variables so that $L$ is the projection to the first coordinate, the Hessian matrix of $Q\vert_{L=0}$ is a sub-matrix of the Hessian matrix of $Q$, with the first row and column deleted.

(1):
If $\trd(Y) = 0$, then $Y^2 = -\nrd(Y)$ by Cayley--Hamilton,
so $\trd(WY^2) = -\trd(W)\nrd(Y)$.
If $\trd(W) \ne 0$, then
\begin{equation*}
\rank(\trd(WY^2) \vert_{\trd(Y)=0})
= \rank(\nrd(Y) \vert_{\trd(Y)=0}).
\end{equation*}
After passing to an algebraic closure we see that this rank is $\rank(ad-bc\vert_{a+d=0})
= 3,$
so $\rank(\trd(WY^2)) \in \{3,4\}$.
Now suppose $\trd(W) = 0$.
Then by Cayley--Hamilton,
\begin{equation*}
\trd(WY^2)
= \trd(Y)\trd(WY) - \nrd(Y)\trd(W)
= \trd(Y)\trd(WY).
\end{equation*}
But the linear forms $\trd(Y)$ and $\trd(WY)$ are linearly independent, because $\trd(W)=0$ and $W\ne 0$.\footnote{If we had $\lambda\, \trd(Y) = \mu\, \trd(WY)$, then taking $Y=1$ would force $\lambda=0$.
But $\trd(WY)$ is a nonzero linear form, since $W\ne 0$; so we would then have $\mu=0$ as well.}
Therefore, $\rank(\trd(WY^2)) = 2$ in the case $\trd(W) = 0$.

(2):
Immediate from (1),
since $P(Y) = \sum_{1\le i\le n} \coeff_i Y_i^2$.

(3):
Immediate from (2), the determinantal definition of rank, and compactness.

(4):
Immediate from the determinantal definition of rank, and compactness.
%
\end{proof}

In addition, we need a \emph{uniform} diagonalization result for quadratic forms over local fields.

\begin{lem}
\label{LEM:uniform-diagonalization}
Let $c\defeq 1$ unless $F$ has residue characteristic $2$, in which case let $c\defeq \abs{2}_F^{-1}$.
Let $H=H^t\in M_N(F)$ be a symmetric matrix.
Then there exists a matrix $A\in M_N(F)$ 
with entries of absolute value less than or equal to $c^N$
and with $\abs{\det{A}} = 1$,
such that $A^tHA$ is diagonal.
\end{lem}

\begin{proof}
We may assume $H\ne 0$.
If $F=\RR$ (resp.~$F=\CC$), the spectral theorem (resp.~Autonne--Takagi factorization) for symmetric matrices over $F$ lets us take $A\in \GL_N(F)\cap U(N)$ unitary.
Now assume $F$ non-Archimedean.
By compactness, there exists a point $x=x_0\in \OO_F^N$ maximizing the quantity $\abs{x^tHx}$.
Because $H$ is nonzero $x_0\ne 0$.
Since dividing $x_0$ by $\varpi_F$ cannot decrease $\abs{x_0^tHx_0}$, we may further assume $x_0$ primitive.
So by an $\OO_F$-linear change of variables, we may assume $x_0=(0,\dots,0,1)$.
By completing the square in the last variable of the quadratic form $x^tHx$, there exists a matrix $A_0\in 2^{-1}M_N(\OO_F)$ with $\det(A_0)=1$ such that $A_0^tHA_0 = \begin{psmatrix} H_0 & 0 \\ 0 & x_0^tHx_0 \end{psmatrix}$.
By induction on $N$, there exists a matrix $A_1\in 2^{1-N}M_{N-1}(\OO_F)$ with $\det(A_1)=1$ such that $A_1^tH_0A_1$ is diagonal.
To complete the proof, take $A=A_0\begin{psmatrix} A_1 & 0 \\ 0 & 1 \end{psmatrix}$.
%
\end{proof}

We also need stationary phase in the clean form suggested by \cite[Lemma~2.7]{browning2015rational}.

\begin{lem}
\label{LEM:reduce-to-near-stat-phase}
Assume $F\cong \RR$.
Fix an integer $N\ge 0$.
Let $w,\Psi\in C^\infty(F^N)$.
Suppose $w$ is supported on $\norm{\bm{x}}\le R \ll 1$,
and that its $k$th order derivatives are $\ll_k 1$ for $\bm{x}\in F^N$, for $k\ge 0$.
Suppose the $k$th order derivatives of $\Psi$ are $\ll_k 1$ for $\norm{\bm{x}}\le R+1$, for $k\ge 0$.
If $T\ge 1$, then
\begin{equation*}
\int_{F^N} w(\bm{x}) \psi(T\Psi(\bm{x})) \,d\bm{x}
\ll_{A,\epsilon} \vol(\Omega_{T,\epsilon}) + T^{-A}
\end{equation*}
for all $A\ge 0$ and $\epsilon\in (0, \frac12)$,
where $\Omega_{T,\epsilon} \defeq \{\bm{x}\in \operatorname{Supp}(w): \norm{\nabla{\Psi}(\bm{x})} \le T^{-1/2+\epsilon}\}$.
\end{lem}

\begin{proof}
Fix a smooth partition of unity $1=\nu_0+\nu_1+\dots+\nu_N$ of $F^N$
such that $\nu_0$ is supported on $\norm{\bm{x}} \le 1$,
and $\nu_i$ for $1\le i\le N$ is supported on $\abs{x_i} \ge \frac12$.
Trivially
\begin{equation*}
\int_{F^N} \nu_0(T^{1/2-\epsilon}\nabla{\Psi}(\bm{x}))
w(\bm{x}) \psi(T\Psi(\bm{x})) \,d\bm{x}
\ll \vol(\Omega_{T,\epsilon}).
\end{equation*}
Also, if $1\le i\le N$ then by non-stationary phase in $\bm{x}$ \cite[Theorem 7.7.1]{Hormander:PDI} 
\begin{equation*}
\int_{F^N} \nu_i(T^{1/2-\epsilon}\nabla{\Psi}(\bm{x}))
w(\bm{x}) \psi(T\Psi(\bm{x})) \,d\bm{x}
\ll_{A,\epsilon} T^{-A}.
\end{equation*}
\quash{because the derivative bounds
\begin{equation*}
\frac{1 + \norm{\frac{\partial}{\partial{x_i}}(T^{1/2-\epsilon}\nabla{\Psi}(\bm{x}))}}
{\abs{T\frac{\partial\Psi}{\partial{x_i}}(\bm{x})}}
\ll \frac{T^{1/2-\epsilon}}
{T^{1/2+\epsilon}}
= T^{-2\epsilon},
\quad \frac{\abs{\frac{\partial}{\partial{x_i}}(T\frac{\partial\Psi}{\partial{x_i}}(\bm{x}))}}
{\abs{T\frac{\partial\Psi}{\partial{x_i}}(\bm{x})}^2}
\ll \frac{T}
{(T^{1/2+\epsilon})^2}
= T^{-2\epsilon},
\end{equation*}
hold whenever $\nu_i(T^{1/2-\epsilon}\nabla{\Psi}(\bm{x})) w(\bm{x}) \ne 0$.
The former derivative bound accounts for the product and chain rules on the $\nu_i w$ factor,
while the latter derivative bound accounts for the quotient rule on any factors of $T\frac{\partial\Psi}{\partial{x_i}}$ in the denominator.}
\end{proof}

Recall for $(f,\Phi,\delta,\gamma) \in \mathcal{S}(D^n) \times \mathcal{S}(D^2) \times D^\times \times D^n$ we define
\begin{equation}
\begin{split}
\label{EQN:define-local-integrals-I_0,I_1}
I_{0}(\delta,\gamma):=I_{0}(f,\Phi,\delta,\gamma) &:= \int_{D^n} f(Y) \Phi(\delta, \delta^{-1} P(Y)) \psi(\trd(\gamma\cdot Y)/\nrd(\delta)) \, dY.
\end{split}
\end{equation}
We now have the following integral result at infinite places $v$.
A similar result holds when $F\cong \CC$, but we focus on $F\cong \RR$ for notational simplicity.
\begin{prop}
\label{PROP:infinite-place-integral-bounds}
Assume $F\cong \RR$.
Fix $f\in C^\infty_c(D^n)$ and $\Phi\in C^\infty_c(D^2)$.
Let $\delta\in D^\times$ and $\gamma\in D^n$.
Then $I_0\left(\frac{\delta}{X},\frac{\gamma}{X}\right) = 0$ unless $\norm{\delta/X} \ll 1$.
Moreover, for all $A,\epsilon>0$,
$$I_0\left(\frac{\delta}{X},\frac{\gamma}{X}\right)
\ll_{A,\epsilon}  (1 + \norm{\gamma}/\norm{\delta})^{-A}
(1 + \norm{X\gamma/\nrd(\delta)})^{\epsilon-n}.$$
\end{prop}

\begin{proof}
The first claim is clear since $\Phi$ is assumed to have compact support.
It remains to prove the second claim, assuming that $\norm{\delta/X} \ll 1$, so $1\ll \norm{X\delta^{-1}}.$ Fix $A,\epsilon>0$ and $B\gg_{A,\epsilon} 1$.
By non-stationary phase in $Y$ (i.e.~repeated integration by parts, integrating the $\psi$ factor and differentiating the complementary factor), we have
\begin{equation*}
I_0\left(\frac{\delta}{X},\frac{\gamma}{X}\right)
\ll_B \int_{D^n} \one_{\norm{Y}\ll 1}
\frac{\norm{X\delta^{-1}}^B}{\norm{X\gamma/\nrd(\delta)}^B} \, dY,
\end{equation*}
because differentiating $X \delta^{-1} P(Y)$ introduces a factor of $\ll \norm{X \delta^{-1}}.$
But $\delta^{-1} = \delta^\dagger/\nrd(\delta)$
and $\norm{\delta^\dagger} \asymp \norm{\delta}$, so this implies
\begin{equation*}
I_0\left(\frac{\delta}{X},\frac{\gamma}{X}\right)
\ll_B \frac{\norm{\delta}^B}{\norm{\gamma}^B}.
\end{equation*}
This estimate, combined with the trivial estimate $I_0 \ll 1$,
already suffices, except in the case where we have, say,
\begin{equation}
\label{key-gamma-delta-infty-range}
\norm{\gamma} \le (1+\norm{X\gamma/\nrd(\delta)})^\omega \norm{\delta},
\end{equation}
where $\omega \defeq \min(\frac{\epsilon}{2n+A},\frac{1}{100})$ is a parameter whose purpose will become clear later.

To go further, we will use Lemma~\ref{LEM:general-rank-bounds}.
First, by Fourier inversion, i.e.~\eqref{fourier-inversion}, we have
\begin{align*}
I_0\left(\frac{\delta}{X},\frac{\gamma}{X}\right)
= \int_{D} \mathcal{F}_2(\Phi)(\delta/X,Z)
\int_{D^n} f(Y) \overline{\psi}\left(\frac{\trd(ZX\delta^\dagger P(Y) - X\gamma\cdot Y)}{\nrd(\delta)}\right) \, dY \, dZ.
\end{align*}
By continuity, we may restrict to $Z\ne 0$.
After identifying $Y\in D^n$ with $\bm{x}\in F^{4n}$,
let $$\bm{u} \defeq X\gamma/\nrd(\delta) \in F^{4n},
\quad \lambda \defeq X\norm{Z\delta^\dagger}/\nrd(\delta) \in F,
\quad G(\bm{x}) \defeq \trd(Z\delta^\dagger P(Y))/\norm{Z\delta^\dagger}.$$
We now apply stationary phase, in the form of Lemma~\ref{LEM:reduce-to-near-stat-phase},
to the inner integral over $Y$,
with $T=\max(1,\abs{\lambda},\norm{\bm{u}})$ and $\Psi(\bm{x})=\frac{\lambda G(\bm{x})-\bm{u}\cdot\bm{x}}{T}$.
We may assume $\norm{\bm{u}} \gg 1,$ or there is nothing to prove.
 Lemma~\ref{LEM:uniform-diagonalization}, together with parts (3) and (4) of Lemma~\ref{LEM:general-rank-bounds}, implies that
up to a bounded change of variables over $F$, the quadratic form $G(\bm{x})$ is of the form $\sum_{1\le i\le 4n} a_i x_i^2$, where $\abs{a_1},\dots,\abs{a_{2n}} \asymp 1$ and $a_{2n+1},\dots,a_{4n} \ll 1$.
The measure of the ``near-stationary set''
\begin{equation*}
\{\bm{x}\in \operatorname{Supp}(f):
\norm{\lambda\nabla{G}(\bm{x}) - \bm{u}} \le T^{1/2+\omega}\}
\end{equation*}
is $0$ if $\lambda\ll \norm{\bm{u}}$ (with a small enough implied constant),
and is $\ll \prod_{1\le i\le 4n} \min(1, \frac{T^{1/2+\omega}}{\abs{\lambda a_i}}) \ll (T^{-1/2+\omega})^{2n}$ for all other $\lambda$.
Since $T\ge \norm{\bm{u}}$, we conclude that
\begin{equation*}
I_0\left(\frac{\delta}{X},\frac{\gamma}{X}\right)
\ll_{B,\omega} \int_{D} (1+\norm{Z})^{-B}
\left(\norm{X\gamma/\nrd(\delta)}^{-n+2n\omega}
+ \norm{X\gamma/\nrd(\delta)}^{-B}\right) \, dZ.
\end{equation*}
But $(1 + \norm{\gamma}/\norm{\delta})^A \ll_A (1+\norm{X\gamma/\nrd(\delta)})^{A\omega}$ by \eqref{key-gamma-delta-infty-range}.
So the last display is
$$
\ll_{A,\epsilon} (1 + \norm{\gamma}/\norm{\delta})^{-A} (1 + \norm{X\gamma/\nrd(\delta)})^{\epsilon-n},$$
because $2n\omega + A\omega \le \epsilon$ by definition of $\omega$.
This completes the proof of the proposition.
\end{proof}

Now assume $F$ is non-Archimedean for the rest of the section.
For $x\in D$, let
\begin{equation*}
v(x)\defeq \max\{e\in \ZZ\cup \{\infty\}: x\in \varpi^e\OO_D\}.
\end{equation*}
Let
\begin{equation*}
v(\gamma) := \min(v(\gamma_1),\dots,v(\gamma_n)),
\quad \norm{\gamma} := q^{-v(\gamma)}.
\end{equation*}
We have the following analogue of Proposition~\ref{PROP:infinite-place-integral-bounds}.

\begin{prop}
\label{PROP:nonsplit-finite-place-integral-bounds}
Let $(f,\Phi) \in C_c^\infty(D^n) \times C_c^\infty(D^2).$ 
\begin{enumerate}
\item \label{1} If $I_0(\delta,\gamma) \ne 0$,
then $v(\delta)\ge -O_{\Phi}(1)$ and $v(\gamma) \ge v(\delta) - O_{f,\Phi}(1)$.

\item \label{2} $I_0(\delta,\gamma) \ll_{f,\Phi} q^{\min(v(\gamma)-v(\nrd(\delta)),0)\,n}$.
\end{enumerate}
\end{prop}

\begin{proof}
By Fourier inversion, i.e.~\eqref{fourier-inversion}, we have
\begin{equation*}
I_0(\delta,\gamma) = \int_D \mathcal{F}_2(\Phi)(\delta,Z)
\int_{D^n} f(Y) \overline{\psi}\left(\frac{\trd(Z\delta^\dagger P(Y) - \gamma\cdot Y)}{\nrd(\delta)}\right) \, dY \, dZ.
\end{equation*}
Here $\mathcal{F}_2(\Phi)(\delta,Z)\, f(t) = 0$
unless $(\delta,Z)\in K_\Phi^2$ and $t\in K_f^n$, where $K_\Phi\subset D$ and $K_f\subset D^n$ are suitable compact sets depending on $\Phi$ and $f$, respectively.

\eqref{1}:
First, $\delta\in K_\Phi$ implies $v(\delta)\ge -O_\Phi(1)$.
Also, for any $W\in \OO_D$, the shift $Y\mapsto Y+\varpi^{v(\nrd(\delta))-v(\delta)+k}W$ leaves $f(Y)$ and the residue $Z\delta^\dagger P(Y) \bmod{\nrd(\delta)}$ invariant
provided $k$ in $\ZZ$ is sufficiently large in a sense depending on $f$, $K_\Phi$, and $P$.
Therefore, on averaging over $W\in \OO_D$ using the definition of $\mathcal{F}_D(\one_{\OO_D})(\varpi^{v(\nrd(\delta))-v(\delta)+k}\gamma/\nrd(\delta))$,
we get $I_0 = 0$
unless $\varpi^{-v(\delta)+k}\gamma \in \operatorname{Supp}(\mathcal{F}_D(\one_{\OO_D}))$.
It follows that $v(\gamma) \ge v(\delta)-k-1$.

\eqref{2}:
By linearity in $f$, we may assume $f$ is the indicator function of a small box.
Fix an identification $D\cong F^4$.
For each $W\in D$, let $H_W\in M_{4n}(F)$ be the Hessian matrix of the $4n$-variable quadratic form $WP(t)$ over $F$.
Given $(\delta,Z)\in K_\Phi^2$, let $W\defeq Z\delta^\dagger/\varpi^{v(Z\delta^\dagger)}\in \OO_D$ (assuming $Z\ne 0$).
This definition ensures that $v(W)=0$ and
\begin{equation*}
I(Z) \defeq
\int_{D^n} f(Y) \overline{\psi}\left(\frac{\trd(Z\delta^\dagger P(Y) - \gamma\cdot Y)}{\nrd(\delta)}\right) \, dY
= \int_{D^n} f(Y) \overline{\psi}\left(\frac{\trd(\varpi^{v(Z\delta^\dagger)} W P(Y) - \gamma\cdot Y)}{\nrd(\delta)}\right) \, dY.
\end{equation*}
By Lemma~\ref{LEM:uniform-diagonalization}, there exists a matrix $A\in 2^{-4n}M_{4n}(\OO_F)$, with $\abs{\det{A}}=1$, such that $A^tH_WA$ is diagonal.
By Lemma~\ref{LEM:general-rank-bounds}(3)--(4), the diagonal matrix $A^tH_WA$ has at least $2n$ entries of absolute value $\gg 1$, since $\norm{W}=1$.
Thus by \eqref{eq:gaussbd} we have
\begin{equation*}
I(Z) \ll_{q,f} (q^{\min(v(Z\delta^\dagger)-v(\nrd(\delta)),0)/2})^{2n}.
\end{equation*}
On the other hand, arguing as in (1), we have $I(Z) = 0$ unless $v(\gamma) \ge \min(v(Z\delta^\dagger),v(\nrd(\delta))) - O_{f,\Phi}(1)$.
Thus $I(Z) \ll_{f,\Phi} (q^{\min(v(\gamma)-v(\nrd(\delta)),0)})^{n}
$.
This extends to $Z=0$ by continuity.
\end{proof}

When $D$ is nonsplit,
we will ultimately not need much more than Proposition~\ref{PROP:nonsplit-finite-place-integral-bounds}.
In fact, for analytic purposes, nonsplit $D$ turn out to be simpler
than split $D$,
because fewer moduli appear.
The key observation is this:

\begin{lem}
\label{LEM:nonsplit-primitive-norms}
%
If $D$ is nonsplit,
then $ 2v(\delta)+1 \ge v(\mathrm{nrd}(\delta)) \ge 2v(\delta).$
\end{lem}

\begin{proof}
%
By scaling $\delta$ by an appropriate power of $\varpi$,
we may assume $v(\delta) = 0$, i.e.~$\delta\in \OO_D - \varpi\OO_D$.
%
%
If $D$ is nonsplit, then $\OO_D = \{x\in D: v(\nrd(x)) \ge 0\}$ by \cite[\S~13.3, Proposition~13.3.4]{voight2021quaternion}.
So $\varpi\OO_D = \{x\in D: v(\nrd(x)) \ge 2\}$.
Since $\delta\in \OO_D - \varpi\OO_D$, it follows that $0\le v(\nrd(\delta)) \le 1$.
\quash{

But $\eta \mapsto \nrd(\eta)$
defines a continuous map $\phi:
\OO_{D} - \varpi\OO_{D}
\to F_v^\times$,
since $D$ is nonsplit.
The domain of $\phi$ is compact in the $v$-adic topology, and therefore the image of $\phi$ is too.
In particular, we have $v(\nrd(\delta)) \in v(\im(\phi)) \ll 1$.}
\end{proof}

Split $D$ are much more intricate.
They are the focus of the next section.

\section{General exponential sums}
\label{SEC:local-exp-sums}


Adopt the setting of Proposition~\ref{PROP:nonsplit-finite-place-integral-bounds} from the last section,
with $F$ non-Archimedean, $\gamma\in D^n$, and $\delta\in D^\times$.
Our goal in this section is to obtain sharper integral estimates when $D$ is split. We assume $P(Y)=\sum_{1\le i\le n} \coeff_i Y_i^2$, where $\coeff_i\in \OO_{F}^\times.$ With more notation, we could likely handle general $\coeff_i\in F^\times$, but this would obscure the main ideas. Our main result in this setting is the estimate of Theorem \ref{thm:estimate} below.  

Throughout this section, we assume $D$ is split and we will often identify $D\cong M_2(F)$ over $F$ and $\OO_D\cong M_2(\OO_F)$ over $\OO_F$. We also assume $f=\one_{\OO_D^n}$, $\Phi=\one_{\OO_D^2},$ and that $\psi$ is unramified.  Then by \eqref{EQN:define-local-integrals-I_0,I_1}, we have
\begin{equation}
\begin{split}
\label{EQN:standard-I0-I1}
I_0(\delta,\gamma) &= \one_{\OO_{D}}(\delta) \int_{\OO_{D}^n}
\one_{\OO_{D}}(\delta^{-1} P(Y)) \psi(\trd(\gamma\cdot Y)/\nrd(\delta)) \, dY.
\end{split}
\end{equation}
In view of \eqref{EQN:standard-I0-I1}, we will henceforth assume $\delta \in \OO_D\cap D^\times$. By Fourier inversion, self-duality of $\one_{\OO_{D}}$, and the identity $\delta^{-1} = \delta^\dagger / \nrd(\delta)$, we have 
\begin{equation}
\label{standard-case-I0-fourier}
I_0(\delta,\gamma) = \int_{\OO_{D}} I_0(Z,\delta,\gamma) \,dZ,
\end{equation}
where
\begin{equation*}
I_0(Z,\delta,\gamma) := \int_{\OO_{D}^n}
\psi\left(\frac{\trd(Z\delta^\dagger P(Y) + \gamma\cdot Y)}{\nrd(\delta)}\right) \,dY.
\end{equation*}

Via the invertible affine map $Y_i\mapsto Y_i+\delta W_i$ on $\OO_D$ (where $W_i\in \OO_D$),
\begin{equation*}
I_0(Z,\delta,\gamma) = \int_{\OO_{D}^n} \int_{\OO_{D}^n}
\psi\left(\frac{\trd(Z\delta^\dagger P(Y+\delta W) + \gamma\cdot (Y+\delta W))}{\nrd(\delta)}\right) \,dY \,dW.
\end{equation*}
But $\delta^\dagger P(Y+\delta W) \equiv \delta^\dagger P(Y) + \sum_{1\le i\le n} \coeff_i \delta^\dagger Y_i \delta W_i\bmod{\nrd(\delta)}$, since $\delta^\dagger \delta \equiv 0 \bmod{\nrd(\delta)}$.
So on averaging over $W$, we conclude that
\begin{equation}
\label{EQN:prelim-vanishing}
I_0(Z,\delta,\gamma) = \int_{\OO_{D}^n}
\one_{\mathrm{nrd}(\delta)\OO_D^n}(\coeff_iZ\delta^\dagger Y_i\delta + \gamma_i\delta)\,
\psi\left(\frac{\trd(Z\delta^\dagger P(Y) + \gamma\cdot Y)}{\nrd(\delta)}\right) \,dY.
\end{equation}
The condition on $(\coeff_iZ\delta^\dagger Y_i\delta + \gamma_i\delta) \in \OO_D^n$ leads to new, nonabelian \emph{vanishing phenomena}, which we alluded to in the introduction of the paper, for the local integrals $I_0$.

This condition leads us to study the image of the set $\eta^\dagger \OO_D \eta := \{\eta^\dagger A \eta: A\in \OO_D\}$ in $\OO_D/\nrd(\eta)\OO_D,$ denoted  by $\eta^\dagger\frac{\OO_D}{\nrd(\eta)\OO_D}\eta$.

\begin{lem}\label{LEM:delta-module-arithmetic}
    Suppose  $\eta\in \OO_D-\varpi\OO_D$. We have an  isomorphism of $\OO_F$-modules
    $$\eta^\dagger\frac{\OO_D}{\nrd(\eta)\OO_D}\eta
\cong \OO_F/\nrd(\eta)\OO_F.$$
In particular, $\eta^\dagger\frac{\OO_D}{\nrd(\eta)\OO_D}\eta$ is a cyclic $\OO_F$-module, generated by some $N\in \eta^\dagger \OO_D \eta$.
\end{lem}

\begin{proof}
    
Let $m = \nrd(\eta)$.
If $\alpha,\beta\in \OO_D^\times$, then $(\alpha\eta\beta)^\dagger\frac{\OO_D}{m\OO_D}\alpha\eta\beta = \beta^\dagger\eta^\dagger\frac{\OO_D}{m\OO_D}\eta\beta \cong \eta^\dagger\frac{\OO_D}{m\OO_D}\eta$.
Thus, by the Cartan decomposition
we may assume $\eta = \begin{psmatrix}m & 0 \\ 0 & 1\end{psmatrix}$, since $\eta\in \OO_D-\varpi\OO_D$.
If $A = \begin{psmatrix}a & b \\ c & d\end{psmatrix}$,
then $\eta^\dagger A \eta
= \begin{psmatrix}1 & 0 \\ 0 & m\end{psmatrix} \begin{psmatrix}a & b \\ c & d\end{psmatrix} \begin{psmatrix}m & 0 \\ 0 & 1\end{psmatrix}
\equiv \begin{psmatrix}0 & b \\ 0 & 0\end{psmatrix} \bmod{m}$.
\end{proof}

Write $\delta = \varpi^{v(\delta)} \eta$, so $\eta\in \OO_D-\varpi\OO_D$ and $\nrd(\delta) = \varpi^{2v(\delta)} \nrd(\eta)$.
The restriction $\coeff_iZ\delta^\dagger Y_i\delta + \gamma_i\delta \equiv 0 \bmod{\nrd(\delta)}$ in \eqref{EQN:prelim-vanishing} implies first that
\begin{equation}
\label{EQN:e-gamma_i-divisibility}
\coeff_iZ\delta^\dagger Y_i + \gamma_i \in \OO_D\delta^\dagger,
\textnormal{ whence }
\gamma_i = \varpi^{v(\delta)} \gamma'_i
\end{equation}
for some $\gamma'_i\in \OO_D$,
and second that
\begin{equation}
\label{prelim-vanishing-condition}
\coeff_iZ\eta^\dagger Y_i\eta + \gamma'_i\eta \equiv 0 \bmod{\nrd(\eta)}.
\end{equation}
Letting $N = \eta^\dagger N_0 \eta \in \eta^\dagger \OO_D \eta$ generate $\eta^\dagger\frac{\OO_D}{\nrd(\eta)\OO_D}\eta$ as in Lemma~\ref{LEM:delta-module-arithmetic},
it follows that there exist $\lambda_1,\dots,\lambda_n\in \OO_F/\nrd(\eta)\OO_F$
such that
\begin{equation}
\label{Y_i-N-coeff}
\eta^\dagger Y_i\eta \equiv \lambda_i N \bmod{\nrd(\eta)},
\end{equation}
and therefore
\begin{equation}
\label{CONG:gamma'_i-proportional}
\gamma'_i\eta \equiv -\coeff_i \lambda_i Z N \bmod{\nrd(\eta)},
\textnormal{ or equivalently, }
\gamma'_i \in -\coeff_i \lambda_i Z\eta^\dagger N_0 + \OO_D\eta^\dagger.
\end{equation}
We will bound the measure of the set of $Z$ satisfying such a congruence on average in Lemma~\ref{LEM:Z-equivalence-classes} below. It strengthens the exponent of the bound in Theorem~\ref{THM:off-center-total-bound} for the sum over $\gamma \neq 0$ from $3n+2+\epsilon$ to  $3n+\epsilon.$

Given $\eta\in \OO_D-\varpi\OO_D$, we define the map
\begin{equation}
\label{EQN:define-local-map-Theta}
\Theta: \OO_F^n \times \OO_D \to (\OO_D/\OO_D\eta^\dagger)^n,
\quad (\lambda,Z) \mapsto \lambda Z\eta^\dagger N_0.
\end{equation}
For $M_0\in \OO_D/\OO_D\eta^\dagger$, let
\begin{equation}
\label{EQN:define-Z-measure-bound-W}
\mathcal{W}(M_0,\eta)
:= \sup_{\substack{\mu\in \OO_F^n-\varpi\OO_F^n}}
\int_{\OO_D}
\one_{\mu M_0\in \Theta(\OO_F^n,Z)}
\,dZ.
\end{equation}
We now have developed enough notation that we can state the main theorem of this section:

\begin{thm} \label{thm:estimate}
Assume $v\nmid 2$.
Let $\gamma \in D^n$.
Write $\delta=\varpi^{v(\delta)}\eta.$
Then $I_0(\delta,\gamma)=0$ unless $\gamma \in \varpi^{v(\delta)}\OO_D^n$ and there exists an $M_0 \in \OO_D$ and a $\mu \in \OO_F^n-\varpi\OO_F^n$ such that 
\begin{align} \label{cond}
\gamma\varpi^{-v(\delta)}\eta \in \mu M_0\eta+\nrd(\eta)\OO_D^n.
\end{align}
Assume $\gamma \in \varpi^{v(\delta)}\OO_D^n$ and write $\gamma':=\gamma \varpi^{-v(\delta)}.$  If \eqref{cond} holds then
\begin{align*}
\abs{I_0(\delta,\gamma)} \leq 
\frac{\mathcal{W}(M_0,\eta)\,|\nrd(\eta)|^{3n/2}\norm{\delta}^n}{\max(\norm{(\gamma'-\gamma'^{\dagger})\eta},|\nrd(\eta)|)^{n/2}\max(\norm{\gamma'},\norm{\delta}|\nrd(\eta)|)^n}.
\end{align*}
\end{thm}
\noindent The proof will occupy the remainder of this section. 

\begin{rem} We point out that Theorem \ref{thm:estimate}
implies that the integral $I_0(\delta,\gamma)$ vanishes unless the image of $\gamma'\in \OO_D^n$ in $(\OO_D/\OO_D\eta^\dagger)^n$ is in the image of the diagonal embedding $\OO_D/\OO_D\eta^\dagger \to (\OO_D/\OO_D\eta^\dagger)^n$, at least  up to $(\OO_F/\nrd(\eta)\OO_F)^n$-scaling.  Thus the exponential sums are supported on a relatively small set of residue classes modulo $\nrd(\eta).$
\end{rem}

In our eventual application of Theorem \ref{thm:estimate} to the proof of Theorem \ref{THM:off-center-total-bound} we will require the average estimate for $\mathcal{W}(M_0,\eta)$ given by the following lemma:

\begin{lem}
\label{LEM:Z-equivalence-classes}
We have
\begin{equation*}
\sum_{M_0\in \OO_F^\times \backslash (\OO_D/\OO_D\eta^\dagger)}
\mathcal{W}(M_0,\eta)
\leq v(\nrd(\eta)) +1.
\end{equation*}
\end{lem}

\begin{proof}
Given $M_0\in \OO_D/\OO_D\eta^\dagger$ and $\mu\in \OO_F^n-\varpi\OO_F^n$,
we have
\begin{equation*}
\mu M_0 \in \Theta(\OO_F^n,Z)
\Leftrightarrow
\mu_1M_0,\dots,\mu_nM_0 \in \OO_F Z\eta^\dagger N_0 + \OO_D\eta^\dagger
\Leftrightarrow
M_0 \in \OO_F Z\eta^\dagger N_0 + \OO_D\eta^\dagger.
\end{equation*}
The condition on the right is independent of $\mu$.
Therefore, by \eqref{EQN:define-Z-measure-bound-W}, we have
\begin{equation*}
\mathcal{W}(M_0,\eta)
= \int_{\OO_D}
\one_{M_0 \in \OO_F Z\eta^\dagger N_0 + \OO_D\eta^\dagger}
\,dZ
= \int_{\OO_D}
\one_{[M_0] \in [\OO_F Z\eta^\dagger N_0 + \OO_D\eta^\dagger]}
\,dZ,
\end{equation*}
where $[A]$ denotes the image of $A$ in the quotient set $\OO_F^\times \backslash (\OO_D/\OO_D\eta^\dagger)$.
Here, if $A$ is an element (resp.~subset) of either $\OO_D$ or $\OO_D/\OO_D\eta^\dagger$, then we interpret $[A]$ as an element (resp.~subset) of $\OO_F^\times \backslash (\OO_D/\OO_D\eta^\dagger)$.
By Fubini's theorem, then,
\begin{equation}
\label{EQN:allowable-Z-fubini}
\sum_{M_0\in \OO_F^\times \backslash (\OO_D/\OO_D\eta^\dagger)}
\mathcal{W}(M_0,\eta)
= \int_{\OO_D}
\sum_{M_0\in \OO_F^\times \backslash (\OO_D/\OO_D\eta^\dagger)}
\one_{M_0 \in [\OO_F Z\eta^\dagger N_0 + \OO_D\eta^\dagger]}
\,dZ.
\end{equation}
Given $Z$, the inner sum over $M_0$ is simply $\#[\OO_F Z\eta^\dagger N_0 + \OO_D\eta^\dagger]$.
However,
\begin{equation*}
[\OO_F Z\eta^\dagger N_0 + \OO_D\eta^\dagger]
= \bigcup_{0\le k\le v(\nrd(\eta)) } [\varpi^k Z\eta^\dagger N_0 + \OO_D\eta^\dagger],
\end{equation*}
which is a union of at most $v(\nrd(\eta)) +1$ distinct elements of $\OO_F^\times \backslash (\OO_D/\OO_D\eta^\dagger)$.
So the right-hand side of \eqref{EQN:allowable-Z-fubini} is less than or equal to $v(\nrd(\eta)) +1$.
\end{proof}

\quash{\begin{rem}
With more work, one can show that $\mathcal{W}(M_0,\eta) \ll q^{\tau-v(\nrd(\eta))}$, where $\tau \defeq \min(v(M_0\eta),v(\nrd(\eta)))$.
This would allow for a direct, but messier, proof of Lemma~\ref{LEM:Z-equivalence-classes}.
\end{rem}}

\subsection{Beginning of the proof  of Theorem \ref{thm:estimate}}

In the rest of the section, we assume $v\nmid 2$. We may assume $I_0(\delta,\gamma) \ne 0$, or there is nothing to prove.
Then by \eqref{standard-case-I0-fourier}, there exists an element $Z\in \OO_D$ with $I_0(Z,\delta,\gamma) \ne 0$.
Now $\gamma \in \varpi^{v(\delta)}\OO_D^n$ by \eqref{EQN:e-gamma_i-divisibility}.
In the congruence \eqref{CONG:gamma'_i-proportional}, replace each $\lambda_i$ with a lift $\lambda_i \in \OO_F \cap F^\times$ and let 
$\gcd(\coeff_1\lambda_1,\dots,\coeff_n\lambda_n) \in F^\times$ be any element such that 
$\gcd(\coeff_1\lambda_1,\dots,\coeff_n\lambda_n)\OO_F = \coeff_1\lambda_1 \OO_F+\dots+\coeff_n\lambda_n\OO_F.
$
Then condition \eqref{cond} holds with
\begin{equation*}
M_0 \defeq \gcd(\coeff_1\lambda_1,\dots,\coeff_n\lambda_n) Z\eta^\dagger N_0 \in \OO_D,
\quad \mu_i \defeq -\coeff_i\lambda_i / \gcd(\coeff_1\lambda_1,\dots,\coeff_n\lambda_n) \in \OO_F,
\end{equation*}
where the factor of $\gcd(\coeff_1\lambda_1,\dots,\coeff_n\lambda_n)\in \OO_F\cap F^\times$ ensures that $\mu = (\mu_1,\dots,\mu_n) \in \OO_F^n-\varpi\OO_F^n$.
We emphasize that $M_0$ and $\mu$ may depend on $(\delta,\gamma)$ and $Z$.
However, all that matters is that $M_0$ does not depend on $i$.


\begin{lem} \label{lem:Wsup}
Let $(\mu,M_0)$ be any pair satisfying \eqref{cond}. We have
$$
\abs{I_0(\delta,\gamma)} \leq \mathcal{W}(M_0,\eta)\,
\sup\limits_{Z \in \OO_D}|I_0(Z,\delta,\gamma)|.
$$
\end{lem}

\begin{proof}
By assumption we have $\gamma' \in \mu M_0 + \OO_D^n\eta^\dagger$.
If $I_0(Z,\delta,\gamma) \ne 0$,
then \eqref{CONG:gamma'_i-proportional} implies $\gamma' \in \Theta(\OO_F^n,Z)$ or equivalently $\mu M_0 \in \Theta(\OO_F^n,Z)$.
Therefore, by \eqref{EQN:define-Z-measure-bound-W} we have $\int_{\OO_D} \one_{I_0(Z,\delta,\gamma) \ne 0} \,dZ \le \mathcal{W}(M_0,\eta)$.
The lemma now follows from \eqref{standard-case-I0-fourier}.
\end{proof}

 By Lemma \ref{lem:Wsup}, to prove Theorem \ref{thm:estimate} it suffices to prove that for any $Z\in \OO_D$
\begin{align}\label{eq:I0Zbound}
    |I_0(Z,\delta,\gamma)|\le \frac{|\nrd(\eta)|^{3n/2}\norm{\delta}^n}{\max(\norm{(\gamma'-\gamma'^{\dagger})\eta},|\nrd(\eta)|)^{n/2}\max(\norm{\gamma'},\norm{\delta}|\nrd(\eta)|)^n}.
\end{align}
By definition $\eta=\delta\varpi^{-v(\delta)}$ and $\gamma'=\gamma \varpi^{-v(\delta)},$ so we can rewrite \eqref{EQN:prelim-vanishing} as
\begin{equation*}
I_0(Z,\delta,\gamma) = \int_{\OO_{D}^n}
\one_{\nrd(\eta)\OO_D}(\coeff_iZ\eta^\dagger Y_i\eta + \gamma'_i\eta)\,
\psi\left(\frac{\trd(Z\eta^\dagger P(Y) + \gamma'\cdot Y)}{\nrd(\eta)\varpi^{v(\delta)}}\right) \,dY.
\end{equation*}
Let 
\begin{align}\label{eq:b}
    b:=\min(v(Z\eta^\dagger), v(\delta)+v(\nrd(\eta))). 
\end{align}
By shifting $Y_i \mapsto Y_i + \nrd(\eta)\varpi^{v(\delta)-b} M_i$ for $M_i \in \OO_D$, we find that $I_0(Z,\delta,\gamma) = 0$ unless
\begin{equation}
\label{DIV:b-gamma'-divisibility}
\gamma'\in \varpi^b \OO_D^n.
\end{equation}
Therefore, we can write 
\begin{align}\label{I_0-convert-delta-to-eta-u-ell}
I_0(Z,\delta,\gamma) = \int_{\OO_{D}^n}
\one_{\nrd(\eta)\OO_D}(\coeff_iZ\eta^\dagger Y_i\eta + \gamma'_i\eta)\,
\psi\left(\frac{\trd(WP(Y) + \rho\cdot Y)}{\nrd(\eta)\varpi^{v(\delta)-b}}\right) \,dY,
\end{align}
where
\begin{align}\label{eq:rhoW}
    \rho &:= \gamma'/\varpi^b \in \OO_D^n,
\qquad W := (Z\eta^\dagger)/\varpi^b \in \OO_D.
\end{align}

From the quadratic form $\trd(WP(Y))$ in \eqref{I_0-convert-delta-to-eta-u-ell}, we proceed to
factor out generalized Gauss sums along various directions, in an appropriate iterative manner.
This is like a ``partial diagonalization'' process, but we will phrase it in the most conceptual way possible. We begin with the most significant direction of cancellation: translation by scalars.
The invertible affine map $Y_i\mapsto Y_i + k_i$, for $k_i\in \OO_{F}$,
leaves $Z\eta^\dagger Y_i\eta \bmod{\nrd(\eta)}$ invariant.
Therefore, changing variables $Y\mapsto Y+k$ in \eqref{I_0-convert-delta-to-eta-u-ell} gives
\begin{equation}
\label{scalar-shifted-I0}
I_0(Z,\delta,\gamma) = \int_{\OO_{D}^n}
\one_{\nrd(\eta)\OO_D}(\coeff_iZ\eta^\dagger Y_i\eta + \gamma'_i\eta)\,
\psi\left(\frac{\trd(W P(Y) + \rho\cdot Y)}{\nrd(\eta)\varpi^{v(\delta)-b}}\right)
J_1(Y,Z) \,dY,
\end{equation}
where
\begin{equation*}
J_1(Y,Z) := \int_{ \OO_{F}^n} \psi\left(\frac{\trd(W (P(Y+k)-P(Y)) + \rho\cdot k)}{\nrd(\eta)\varpi^{v(\delta)-b}}\right) \,dk.
\end{equation*}
Here $P(Y+k)-P(Y)
= \sum_{1\le i\le n} (\coeff_ik_i^2+2\coeff_ik_iY_i)$,
because $k_i$ commutes with $Y_i$.

Let
\begin{align}\label{eq:a}
a &:= \begin{cases}
    \trd(W)=\trd(Z\eta^\dagger)\varpi^{-b}  & \textrm{ if }  \one_{\trd(W)\mid \nrd(\eta)\varpi^{v(\delta)-b}},\\
    \nrd(\eta)\varpi^{v(\delta)-b} & \textrm{ if }  \one_{\trd(W) \nmid \nrd(\eta)\varpi^{v(\delta)-b}},
\end{cases}
\end{align}
so $v(a) = \min(v(\trd(Z\eta^\dagger)),v(\delta)+v(\nrd(\eta)))-b
\ge 0,$ i.e., $a\in \OO_F\cap F^\times.$
Since $v\nmid 2$, we may use \eqref{gauss-sum-non-vanishing-condition} and \eqref{evaluate-off-center-gauss-sum} to evaluate $J_1$ and obtain\footnote{If $v\mid 2$, this approach would require extra (doable) casework.
But we note that a less explicit version of the current method, based on changes of variables of the form $Y\mapsto Y+k+R$ restricted to $\trd(WR)=0$, would likely work at once for all $v$.
However, we find the present iterative approach much more enlightening.}
\begin{equation}
\label{J1-formula}
\begin{split}
J_1(Y,Z)
&= \prod_{1\le i\le n} \mathcal{G}(\coeff_ia,\nrd(\eta)\varpi^{v(\delta)-b}, \trd(2\coeff_iWY_i+\rho_i)) \\
&= \one_{a\OO_F}(\trd(2\coeff_iWY_i+\rho_i))\,
\psi\left(-\sum_{1\le i\le n}\frac{\trd(2\coeff_iWY_i+\rho_i)^2}
{4\coeff_ia\,\nrd(\eta)\varpi^{v(\delta)-b}}\right)
G_1(W),
\end{split}
\end{equation}
where $G_1(W) \defeq \prod_{1\le i\le n} \mathcal{G}(\coeff_ia,\nrd(\eta)\varpi^{v(\delta)-b}, 0)$.
Here $G_1(W)$ depends (implicitly) on $\delta$ and $Z$,
but not on $Y$.
Moreover, by \eqref{eq:gaussbd}, we have
\begin{equation}
\label{G1W-product-gauss-bound}
\abs{G_1(W)}
\le \prod_{1\le i\le n} \frac{\abs{\nrd(\eta)\varpi^{v(\delta)-b}}^{1/2}}
{\abs{\coeff_ia}^{1/2}}
= \frac{q^{v(a)n/2}\abs{\nrd(\eta)}^{n/2}}{(q^{v(\delta)-b})^{n/2}}.
\end{equation}

Substituting the formula \eqref{J1-formula} for $J_1$
into \eqref{scalar-shifted-I0}, we get
\begin{equation}\label{aftergauss}
I_0(Z,\delta,\gamma) = \int_{\OO_{D}^n}
S_1(Y,Z)
\psi\left(\frac{\trd(W P(Y) + \rho\cdot Y)}{\nrd(\eta)\varpi^{v(\delta)-b}}
-\sum_{1\le i\le n}\frac{\trd(2\coeff_iWY_i+\rho_i)^2}{4\coeff_ia\, \nrd(\eta)\varpi^{v(\delta)-b}}\right) \,dY,
\end{equation}
where
\begin{equation}
\label{define-S1}
S_1(Y,Z) :=
\one_{\nrd(\eta)\OO_D}(\coeff_iZ\eta^\dagger Y_i\eta + \gamma'_i\eta)\,
\one_{a\OO_F}(\trd(2\coeff_iWY_i+\rho_i))\,
G_1(W).
\end{equation}
Thus
\begin{align}\label{eq:triv}
\begin{split}
    \abs{I_0(Z,\delta,\gamma)}
&\leq \int_{\OO_{D}^n}
\abs{S_1(Y,Z)} \,dY\\
&=|G_1(W)|\int_{\OO_D^n} \one_{\nrd(\eta)\OO_D}(\coeff_iZ\eta^\dagger Y_i\eta + \gamma'_i\eta)\,
\one_{a\OO_F}(\trd(2\coeff_iWY_i+\rho_i))\, dY.
\end{split}
\end{align}

\begin{lem} \label{lem:cong:bound} 
One has
$$   
\int_{\OO_D} \one_{\nrd(\eta)\OO_D}(\coeff_iZ\eta^\dagger Y_i\eta + \gamma'_i\eta)
\,dY_i \leq q^{h}|\nrd(\eta)|,
$$
where $h := \min(v(ZN),v(\nrd(\eta)))$. 
\end{lem}

\begin{proof}

By Lemma~\ref{LEM:delta-module-arithmetic},
the map
$Y_i \mapsto \coeff_iZ\eta^\dagger Y_i\eta \bmod{\nrd(\eta)}$ has image $\OO_F \coeff_iZN \bmod{\nrd(\eta)}$, which is of order $|\nrd(\eta)|^{-1}q^{-h}$ since $\gcd(\coeff_iZN,\nrd(\eta)) = \varpi^h$.
Thus the kernel of the map has order $|\nrd(\eta)|^{-3}q^h.$
If $\gamma'_i\eta$ is in the image, the integral is thus $|\nrd(\eta)|q^{h},$ otherwise it is zero.
\end{proof}


\begin{lem}
\label{LEM:h=min(b,l)}
We have $h=\min(b,v(\nrd(\eta)))$.
\end{lem}

\begin{proof}
Both sides are zero if $v(\nrd(\eta))=0$.
Now suppose $v(\nrd(\eta))>0$.
Then $$I_2
= \begin{psmatrix} 1-\nrd(\eta) & 0 \\ 0 & \nrd(\eta) \end{psmatrix}
+ \begin{psmatrix} \nrd(\eta) & 0 \\ 0 & 1-\nrd(\eta) \end{psmatrix}
\in \OO_D^\times\eta\OO_D^\times + \OO_D^\times\eta\OO_D^\times,$$
by the Cartan decomposition.
Therefore, $I_2$ lies in the \emph{two-sided ideal} $\OO_D\eta\OO_D$ of $\OO_D$.
Since $N$ generates $\eta^\dagger\frac{\OO_D}{\nrd(\eta)\OO_D}\eta$ by definition, we conclude that
\begin{align*}
h&= \min(v(ZN),v(\nrd(\eta)))= \min_{Y\in Z\eta^\dagger\OO_D\eta}( v(Y),v(\nrd(\eta)))\\
&= \min_{Y\in Z\eta^\dagger\OO_D\eta\OO_D}(v(Y),v(\nrd(\eta)))= \min( v(Z\eta^\dagger),v(\nrd(\eta)))=  \min(b,v(\nrd(\eta))).
\end{align*}
\end{proof}

\begin{lem}
\label{alt:cong:bound}
One has $$\int_{\OO_D} \one_{a\OO_F}(\trd(2\coeff_iWY_i+\rho_i))
\,dY_i = \abs{a}.$$
\end{lem}

\begin{proof}
By the definitions of $b$ in \eqref{eq:b} and $W$ in \eqref{eq:rhoW},
we have $$\min(v(\coeff_iW),v(\nrd(\eta))+v(\delta)-b) = 0.$$
By \eqref{eq:a} we have $v(a)\le v(\nrd(\eta))+v(\delta)-b$,
whence $\min(v(\coeff_iW),v(a)) = 0$.
Thus $\min(v(2\coeff_iW),v(a)) = 0$ since $v\nmid 2$.
Therefore the linear congruence $\trd(2\coeff_iWY_i+\rho_i) \equiv 0 \bmod{a}$ has exactly $\abs{a}^{-3}$ solutions $Y_i\bmod{a}$.
Multiplying by $\abs{a}^4$ gives the desired equality.
\end{proof}

By Lemmas~\ref{lem:cong:bound}, \ref{LEM:h=min(b,l)}, and \ref{alt:cong:bound}, substituting \eqref{G1W-product-gauss-bound} into \eqref{eq:triv} gives
\begin{equation}
\label{trivialbound}
\abs{I_0(Z,\delta,\gamma)}
\leq  (q^{-t})^n
\cdot \frac{q^{v(a)n/2}\abs{\nrd(\eta)}^{n/2}}{q^{(v(\delta)-b)n/2}},
\end{equation}
where 
\begin{align}\label{eq:t}
q^{-t} := \min(q^{h}|\nrd(\eta)|, \abs{a})=\min(q^b|\nrd(\eta)|,\abs{a}).
\end{align}
Equivalently,
\begin{equation}
\label{EQN:alternative-t-form}
\begin{split}
t =  \max(v(\nrd(\eta))-b,v(a)) = v(\nrd(\eta))+v(a) - \min(v(\nrd(\eta)),v(a)+b).
\end{split}
\end{equation}
This completes our study of cancellation in the exponential sum using changes of variables of the form $Y\mapsto Y+k.$

Let us now consider changes of variables of the form $Y\mapsto Y+kR$ in \eqref{aftergauss}, for some $R\in \OO_D$.
Such a change of variables will preserve $S_1(Y,Z)$ if we have $Z\eta^\dagger R\eta\equiv 0\bmod{\nrd(\eta)}$ and $\trd(2WR)\equiv 0\bmod{a}$.  The analogue of the argument above will then introduce Gauss sums  over $k_i\in \OO_F$ with modulus $\nrd(\eta)\varpi^{v(\delta)-b}$
and leading quadratic coefficient
\begin{equation*}
\coeff_i\trd(WR^2) - \trd(2\coeff_iWR)^2/(4\coeff_ia)
=\coeff_i\, (\trd(WR^2) - \trd(WR)^2/a),
\end{equation*}
for each $i$.
Thus for any $R \in \OO_D$ such that $Z\eta^\dagger R\eta\equiv 0\bmod{\nrd(\eta)}$ and $\trd(2WR)\equiv 0\bmod{a}$,
we have 
\begin{align} \label{before:opti}
\abs{I_0(Z,\delta,\gamma)} \leq  \int_{\OO_D^n}\frac{|\nrd(\eta)|^{n/2}|S_1(Y,Z)|}{(q^{v(\delta)-b}|\trd(WR^2) - \trd(WR)^2/a|)^{n/2}}\,dY
\end{align}
by \eqref{eq:gaussbd} applied to the $n$ Gauss sums just described. Applying our previous estimates on $|S_1(Y,Z)|$ using Gauss sum estimates and Lemmas ~\ref{lem:cong:bound} and \ref{alt:cong:bound} as in \eqref{trivialbound} we obtain
\begin{align} \label{before:opti2}
\abs{I_0(Z,\delta,\gamma)}
\leq  \frac{q^{v(a)n/2}|\nrd(\eta)|^{n}}{q^{(t+v(\delta)-b)n}|\trd(WR^2) - \trd(WR)^2/a)|^{n/2}}.
\end{align}
Let us now try to maximize  $|\trd(WR^2) - \trd(WR)^2/a|$.

\begin{lem}
\label{LEM:optimize-final-clean-Gauss-direction-R}
Assume $v(a)<v(\delta)+v(\nrd(\eta))-b$.
Then there exists $R\in \OO_D$,
with $Z\eta^\dagger R\eta\equiv 0\bmod{\nrd(\eta)}$ and $\trd(WR)=0$,
such that $v(\trd(WR^2) - \trd(WR)^2/a) = t$.
\end{lem}

\begin{proof}
Throughout this proof,
we let
\begin{equation*}
m \defeq \min(v(\nrd(\eta)),v(a)+b).
\end{equation*}
By the Cartan decomposition,
there exists an element $\eta_0\in \OO_D\cap D^\times$ such that
$v(\nrd(\eta_0)) = m$
and $\eta \equiv \eta_0 \bmod{\varpi^{m}}$.\footnote{Indeed, if $\eta = \begin{psmatrix}\nrd(\eta) & 0 \\ 0 & 1\end{psmatrix}$, then we could take $\eta_0 = \begin{psmatrix}\varpi^{m} & 0 \\ 0 & 1\end{psmatrix}$.}
This ensures that
\begin{equation}
\label{simplify-W-trace-condition}
Z\eta_0^\dagger + \eta_0 Z^\dagger
= \trd(Z\eta_0^\dagger)
\equiv \trd(Z\eta^\dagger)
= \varpi^b\trd(W)
\equiv 0 \bmod{\varpi^{m}},
\end{equation}
by the definition of $W$ in \eqref{eq:rhoW}
and the inequality $v(a)\le v(\trd(W)).$

Let $R_0 \defeq \eta\eta_0^{-1}
\in (\eta_0+\OO_D\nrd(\eta_0))\eta_0^{-1}
= 1 + \OO_D\eta_0^\dagger$.
Then $\eta = R_0\eta_0$, and $R_0\in \OO_D$, with
\begin{equation}
\label{compute-nrd-R0}
v(\nrd(R_0))
= v(\nrd(\eta)) - v(\nrd(\eta_0))
= v(\nrd(\eta)) - m.
\end{equation}
By our assumption $v(a)<v(\delta)+v(\nrd(\eta))-b$, we have
\begin{equation}
\label{grind-assumption-W-trace-equality}
v(\trd(W)) = v(a)<v(\delta)+v(\nrd(\eta))-b.
\end{equation}
In particular, $W\ne 0$.
So $WR_0\ne 0$, since $R_0\in D^\times$ by \eqref{compute-nrd-R0}.
Let $W'\defeq WR_0/\varpi^{v(WR_0)}$.

Next, we construct $S\in \OO_D$ with $\trd(WR_0S)=0$ and $\nrd(S)\in \OO_F^\times$. 
For any $W'$ we have a hyperplane
$\trd(W'S)=0$ in $\mathbb{P}^3_{\OO_F/\varpi\OO_F}$,
while $\nrd(S)=0$ defines an irreducible quadric surface.
Because no conic passes through every point of $\mathbb{P}^2(\OO_F/\varpi\OO_F)$,\footnote{Indeed, $\#\mathbb{P}^2(\OO_F/\varpi\OO_F) = q^2+q+1$, while a (possibly degenerate) conic has at most $2q+1$ points valued in $\OO_F/\varpi\OO_F$.}
we conclude that exists $S_0\in \OO_D$ such that $\trd(W'S_0)\in \varpi\OO_F$ and $\nrd(S_0)\in \OO_F^\times$.
By Hensel's lemma, there exists $S\in S_0+\varpi\OO_D$ with $\trd(W'S) = 0$.
Then $\nrd(S) \equiv \nrd(S_0) \not\equiv 0 \bmod{\varpi}$, so $S\in \OO_D^\times$.

Finally, let $R \defeq R_0S\in \OO_D$.
Then by \eqref{eq:rhoW},
we have $\trd(Z\eta^\dagger R)=\varpi^b\trd(WR)=0$, so
\begin{equation*}
(Z\eta^\dagger R)\eta
= -(R^\dagger \eta Z^\dagger)\eta
= -(S^\dagger R_0^\dagger)
(R_0\eta_0) Z^\dagger \eta
= -S^\dagger \nrd(R_0)
(\eta_0 Z^\dagger) \eta.
\end{equation*}
By \eqref{simplify-W-trace-condition}, \eqref{compute-nrd-R0}, and the congruence $\eta \equiv \eta_0 \bmod{\varpi^{m}}$, we conclude that
\begin{equation*}
Z\eta^\dagger R\eta
\equiv S^\dagger \nrd(R_0)
(Z\eta_0^\dagger) \eta_0
= S^\dagger \nrd(R_0)
Z \nrd(\eta_0)
\equiv 0
\bmod{\nrd(\eta)},
\end{equation*}
since $\nrd(R_0) \nrd(\eta_0) = \nrd(\eta)$.
Finally, by Cayley--Hamilton,
\begin{equation*}
\trd(WR^2) - \trd(WR)^2/a
= \trd(R)\trd(WR) - \nrd(R)\trd(W) - 0
= -\nrd(R)\trd(W),
\end{equation*}
since $\trd(WR)=0$.
Since $S\in \OO_D^\times$, we have $v(\nrd(R)) = v(\nrd(R_0))$, so
\begin{equation*}
v(\nrd(R)) + v(\trd(W))
= v(\nrd(\eta)) - m + v(a)
= t,
\end{equation*}
by \eqref{compute-nrd-R0}, \eqref{grind-assumption-W-trace-equality}, and \eqref{EQN:alternative-t-form}.
Thus $v(\trd(WR^2) - \trd(WR)^2/a) = t$ and we are done.
\end{proof}

\subsection{Proof of the inequality \eqref{eq:I0Zbound}}

Consider first the case $v(a)< v(\nrd(\eta))+v(\delta)-b$, then by \eqref{before:opti2} and Lemma \ref{LEM:optimize-final-clean-Gauss-direction-R}, we have 
\begin{align} \label{before:rem:w}
\abs{I_0(Z,\delta,\gamma)}
\leq \frac{q^{v(a)n/2}|\nrd(\eta)|^{n}}{q^{(t/2+v(\delta)-b)n}}
=q^{(v(\nrd(\eta))+v(a)-t+2b)n/2}|\nrd(\eta)|^{3n/2}\norm{\delta}^n.
\end{align}
The same bound applies to the case $v(a)\ge v(\nrd(\eta))+v(\delta)-b$. Indeed, in this case we have $t=v(a) = v(\delta)+v(\nrd(\eta))-b$ by \eqref{eq:a} and \eqref{EQN:alternative-t-form}. Therefore, by \eqref{trivialbound}
\begin{align*}
    \abs{I_0(Z,\delta,\gamma)}
\leq  (q^{-t})^n
\cdot \frac{q^{v(a)n/2}\abs{\nrd(\eta)}^{n/2}}{q^{(v(\delta)-b)n/2}}=\frac{q^{v(a)n/2}|\nrd(\eta)|^{n}}{q^{(t/2+v(\delta)-b)n}}.
\end{align*}

We claim 
\begin{align*}
    q^{(v(\nrd(\eta))+v(a)-t+2b)n/2}\le \max(\norm{(\gamma'-\gamma'^{\dagger})\eta},|\nrd(\eta)|)^{-n/2}\max(\norm{\gamma'},\norm{\delta}|\nrd(\eta)|)^{-n}
\end{align*}
and thus \eqref{eq:I0Zbound} whence Theorem \ref{thm:estimate} follows from \eqref{before:rem:w}. By \eqref{eq:b} and \eqref{DIV:b-gamma'-divisibility},
we have 
$$
b \le \min(v(\delta)+v(\nrd(\eta)), v(\gamma')).
$$
On the other hand, by Lemma~\ref{LEM:tricky-traceless-part-kappa-divisibility} below if $I_0(Z,\delta,\gamma)\ne 0$ then
\begin{equation*}
\min(v((\gamma' - \gamma^{\prime\dagger})\, \eta),v(\nrd(\eta)))
\geq \min(v(\nrd(\eta)),v(a)+b)
= v(\nrd(\eta))+v(a)-t,
\end{equation*}
where the final equality holds by \eqref{EQN:alternative-t-form}. This justifies the claim.

\begin{lem}
\label{LEM:tricky-traceless-part-kappa-divisibility}
Assume $I_0(Z,\delta,\gamma) \ne 0$.
Then the following congruence holds:
\begin{equation*}
(\gamma'_i - \gamma_i^{\prime\dagger})\, \eta
\equiv 0 \bmod{\varpi^{\min(v(\nrd(\eta)),b+v(a))}\OO_D}.
\end{equation*}
\end{lem}

\begin{proof}
By \eqref{aftergauss}, we have $S_1(Y,Z) \ne 0$ for some $Y$.
Then by \eqref{define-S1} ,
we have $\coeff_iZ\eta^\dagger Y_i\eta + \gamma'_i\eta \equiv 0 \bmod{\nrd(\eta)\OO_D}$,
$\trd(\coeff_iZ\eta^\dagger) = \coeff_i\varpi^b\trd(W) \equiv 0\bmod{\varpi^ba\OO_F}$ (because $v(a)\le v(\trd(W))$),
and $\trd(2\coeff_iZ\eta^\dagger Y_i+\gamma'_i) = \varpi^b\trd(2\coeff_iWY_i+\rho_i) \equiv 0\bmod{\varpi^ba\OO_F}$.
Thus the conditions of Lemma~\ref{LEM:delta-module-arithmetic:4} below hold with parameters
``$(\delta,A,M,Z,K) \defeq
(\eta,Y_i,\gamma'_i,\coeff_iZ,\varpi^{\min(v(\nrd(\eta)),b+v(a))})$'',
and the desired congruence follows by Lemma~\ref{LEM:delta-module-arithmetic:4}.
\end{proof}

\begin{lem}\label{LEM:delta-module-arithmetic:4}
Let $A,M,Z\in \OO_D$ and $K\in \OO_F$ with $\nrd(\delta)\in K\OO_F$.
Suppose $Z\delta^\dagger A\delta + M\delta \equiv 0 \bmod{K\OO_D}$,
$\trd(Z\delta^\dagger) \equiv 0 \bmod{K\OO_F}$,
and $\trd(2Z\delta^\dagger A + M) \equiv 0 \bmod{K\OO_F}$.
Then $(M-M^\dagger) \delta \equiv 0 \bmod{K\OO_D}$.
\end{lem}

\begin{proof}
Since $Z\delta^\dagger + \delta Z^\dagger = \trd(Z\delta^\dagger)I_2 \equiv 0 \bmod{K\OO_D}$,
we have $Z\delta^\dagger \equiv -\delta Z^\dagger \bmod{K\OO_D}$.
On replacing $Z\delta^\dagger$ with $-\delta Z^\dagger$ in the assumption, we get
$M\delta \equiv \delta Z^\dagger A\delta \bmod{K\OO_D}$
and $\trd(M) \equiv \trd(2\delta Z^\dagger A) \bmod{K\OO_F}$.
Right-multiplying the latter by $\delta$ produces the congruence
\begin{equation*}
(M+M^\dagger)\delta
\equiv (2\delta Z^\dagger A + 2A^\dagger Z\delta^\dagger)\delta
\equiv 2\delta Z^\dagger A\delta
\equiv 2M\delta \bmod{K\OO_D}.
\end{equation*}
Subtracting $2M\delta$ from both sides gives the result.
\end{proof}

\qed

\section{Applying the geometry of numbers}
\label{sec:geo:num}

Let $F$ be a number field. Assume $P(\gamma)=\sum_{1\le i\le n} \coeff_i \gamma_i^2$, where $\coeff_i\in \OO_F^\times$ are fixed.
Recall the $\gamma\ne 0$ contribution $E_1(X)$ to \eqref{renormalize-after-poisson},
defined as follows:
\begin{equation*}
E_1(X) :=
c_{\Phi,X}X^{4(n-1)}
\sum_{\gamma \in D^n - \{0\}} \sum_{\delta \in D^\times}
\left(I_{1}\left(\frac{\delta}{X},\frac{\gamma}{X}\right) -  I_0\left(\frac{\delta}{X},\frac{\gamma}{X} \right)\right).
\end{equation*}
This section is devoted to proving the following result:

\begin{thm}
\label{THM:off-center-total-bound}
Assume $F=\QQ$ and $D_\infty,D_2$ nonsplit. Fix  compactly supported smooth functions $f=\otimes_vf_v,$ $\Phi=\otimes_v\Phi_v$ where $f_{v}=\one_{\OO_{D_v}^n}$ and $\Phi_v=\one_{\OO_{D_v}^2}$ whenever $D_v$ is split.
Then $E_1(X) \ll_\epsilon X^{3n+\epsilon}$ for all $\epsilon>0$.
\end{thm}

The idea is to use good exponential sum vanishing results and estimates where possible,
and to control the likelihood of bad events on average over $\gamma$
using Lemma~\ref{lem:badcount} based on the geometry of numbers. The assumption $F=\QQ$ is only to avoid the difficulty of treating multiple infinite places.

To ease the notation, in this section for $M\in \OO_D^r$ we let $\norm{M} \defeq \norm{M}_\infty$ be its box norm at the archimedean place.
\subsection{Proof of Theorem \ref{THM:off-center-total-bound}}
Throughout this proof we say a place $v$ is split (resp.~nonsplit) if $D_v$ is split (resp.~nonsplit). By the symmetry property \eqref{EQN:general-I0-I1-symmetry} 
relating $I_{1}$ to $I_{0}$ and the fact that at split places $v$ one has $\one_{\OO_{D_v}^n}=\one_{\OO_{D_v}^n}\circ\dagger$ and $\one_{\OO_{D_v}^2}=\one_{\OO_{D_v}^2}^{\mathrm{sw}}\circ\dagger,$ it suffices to show
\begin{align}\label{eq:goal}
    \sum_{\gamma \in D^n - \{0\}} \sum_{\delta \in D^\times}
\left|I_{0}\left(\frac{\delta}{X},\frac{\gamma}{X}\right) \right|\ll_\epsilon X^{4-n+\epsilon}.  
\end{align}

Recall that since $f$ and $\Phi$ are pure tensors,
\begin{align*}
    I_{0}\left(\frac{\delta}{X},\frac{\gamma}{X}\right)=I_{0,\infty}\left(\frac{\delta}{X},\frac{\gamma}{X}\right)I_{0}^\infty\left(\delta,\gamma\right)=I_{0,\infty}\left(\frac{\delta}{X},\frac{\gamma}{X}\right)\prod_{v \neq \infty}I_{0,v}(\delta,\gamma).
\end{align*}

\begin{prop}\label{prop:adelic bound}
There exists $C\in \OO_F\cap F^\times$ with $|C|_v=1$ at all split places $v$ such that $I_{0}^\infty\left(\delta,\gamma\right)= 0$ unless $C\delta\in \OO_D\cap D^\times$. Furthermore, there is a positive integer $J$ such that
\begin{align*}
    C\delta&=J\eta, \,\textnormal{ where $\eta\in \OO_D\cap D^\times$ is primitive,}\\
    C^2\gamma&=J\gamma', \textnormal{ where }\gamma'\in \OO_D^n-\{0\},
\end{align*}
and 
\begin{align}\label{eq:supportcond}
    C\gamma'\eta \in \mu M_0\eta+\nrd(\eta)\OO_D^n
\end{align}
for some $M_0\in \OO_D$ and primitive $\mu\in \OO^n_F-\{0\}$. Furthermore, $\left|I_{0}^\infty\left(\delta,\gamma\right)\right|$ is dominated by
\begin{align*}  
\prod_{\substack{v\neq \infty \\ \textrm{ nonsplit }}}\frac{\abs{\nrd(\delta)}_v^n}{\max(\norm{\gamma}_v,\abs{\nrd(\delta)}_v)^n}\prod_{\substack{v\textrm{ split }}}\frac{\mathcal{W}_v(M_0,\eta)\,|\nrd(\eta)|^{3n/2}_v\norm{\delta}_v^n}{\max(\norm{(\gamma'-\gamma'^{\dagger})\eta}_v,|\nrd(\eta)|_v)^{n/2}\max(\norm{\gamma'}_v,\norm{\delta}_v|\nrd(\eta)|_v)^n},
\end{align*}
where $\mathcal{W}_v(M_0,\eta)$ is defined in \eqref{EQN:define-Z-measure-bound-W}.
\end{prop}

\begin{proof}
If we ignore \eqref{eq:supportcond}, the existence of $C$ and $J$ follows from Proposition \ref{PROP:nonsplit-finite-place-integral-bounds}(1) at nonsplit primes and Theorem \ref{thm:estimate} together with \eqref{EQN:standard-I0-I1} at split primes. 

For \eqref{eq:supportcond}, as we can always take $\mu$ to be primitive by scaling $M_0$ and $\mu$ appropriately, it suffices to construct any pair $(\mu,M_0)\in (\OO^n_F-\{0\},\OO_D)$ such that \eqref{eq:supportcond} holds. Observe that the constraint \eqref{eq:supportcond} over a place $v$ is trivially satisfied when $|\nrd(\eta)|_v=1,$ i.e. no restrictions on $\mu$ and $M_0$ need to be imposed. Therefore, it suffices to verify the constraint is satisfied locally at a place $v$. For split $v,$ it follows from Theorem \ref{thm:estimate}. For nonsplit $v,$ as $\eta$ is primitive, we have $v(\nrd(\eta)) \in \{0,1\}$ by Lemma \ref{LEM:nonsplit-primitive-norms}. Then \eqref{eq:supportcond} is satisfied for an appropriate $\mu$ and $M_0$ as long as $v(C)\ge 1,$ which we can always assume to be true. 

Finally, the bound of $\left|I_{0}^\infty\left(\delta,\gamma\right)\right|$ follows from Proposition \ref{PROP:nonsplit-finite-place-integral-bounds}(2) and Theorem \ref{thm:estimate}.
\end{proof}

Since $\Phi$ is compactly supported, if $I_0\left(\frac{\delta}{X},\frac{\gamma}{X}\right) \neq 0,$ then $\norm{\delta}\ll X$ and
\begin{align*}
    I_{0,\infty}\left(\frac{\delta}{X},\frac{\gamma}{X}\right)
\ll_{A,\epsilon}  (1 + \norm{\gamma}/\norm{\delta})^{-A}
(1 + \norm{X\gamma/\nrd(\delta)})^{\epsilon-n}
\end{align*}
for any $A\ge 0$ and $0<\epsilon<1/2$ by Proposition~\ref{PROP:infinite-place-integral-bounds}. Assume first that $\delta$ and $\gamma$ satisfy $\norm{\gamma} \geq X^\epsilon \norm{\delta}.$ Then
$I_{0,\infty}\left(\frac{\delta}{X},\frac{\gamma}{X}\right) \ll_{\epsilon,A} X^{-A}.$  On the other hand, by Proposition \ref{prop:adelic bound}, $I_{0}^\infty(\delta,\gamma)$ is supported in a compact subset of $D_{\A_F^\infty} \times D_{\A_F^\infty}^n$ and one trivially has $I_0^\infty(\delta,\gamma) \ll 1.$
Hence the contribution to $E_1(X)$ of $\delta$ and $\gamma$ satisfying $\norm{\gamma} \geq X^\epsilon \norm{\delta}$ is $\ll_{\epsilon,A} X^{-A}$ for any $A  \in \RR.$ 

 We henceforth assume $\norm{\gamma} \le X^\epsilon \norm{\delta},$ so 
\begin{align}\label{eq:inftyboundp2}
    I_{0,\infty}\left(\frac{\delta}{X},\frac{\gamma}{X}\right)
&\ll_\epsilon
(1 + \norm{X\gamma/\nrd(\delta)})^{\epsilon-n}.
\end{align}
Let us rewrite the sum over $(\delta,\gamma)$ using the constraint \eqref{eq:supportcond}. Let $m=\nrd(\eta)\in \OO_F\cap F^\times.$ For $H,K\in \OO_F\cap F^\times,$ define the $\OO_F$-lattice
\begin{equation*}
\begin{split}
\Lambda(H,K,m,\eta,M_0) :=
\left\{M\in H\OO_D:
\begin{array}{ll}
(M-M^\dagger) \eta \in K\OO_D,\\
M\eta \in \OO_F M_0\eta + m\OO_D
\end{array}\right\}.
\end{split}
\end{equation*}
Note that for a given $\eta,$ $\Lambda(H,K,m,\eta,M_0)$ depends only on the class of $M_0$ in the quotient set
\begin{equation}
\label{global-M_0-quotient-set}
(\OO_F/m\OO_F)^\times \backslash (\OO_D/\OO_D\eta^\dagger)
\cong \prod_{\textrm{finite }v\mid m} (\OO_{F_v}^\times \backslash (\OO_{D_v}/\OO_{D_v}\eta^\dagger)).
\end{equation}

Take $H := \gcd(\gamma', Jm).$ Then we have $C\gamma' \in \Lambda(H,1,m,\eta,M_0)^n$ by definition.  Let $K=K(\gamma',\eta,M_0)$ be the largest divisor of $m$ such that
\begin{equation*}
C\gamma' \in \Lambda(H,K,m,\eta,M_0)^n.
\end{equation*}
Let $S=S(\eta)\subseteq \OO_D$ be a complete set of representatives for the quotient set \eqref{global-M_0-quotient-set}.
Then by replacing $M_0$ with
the representative for its image in \eqref{global-M_0-quotient-set},
we may assume $M_0\in S$, while keeping $\mu$ primitive.

If $v$ is split, we have 
\begin{equation*}
\begin{split}
\frac{\mathcal{W}_v(M_0,\eta)\,|\nrd(\eta)|^{3n/2}_v\norm{\delta}_v^n}{\max(\norm{(\gamma'-\gamma'^{\dagger})\eta}_v,|\nrd(\eta)|_v)^{n/2}\max(\norm{\gamma'}_v,\norm{\delta}_v|\nrd(\eta)|_v)^n} = \frac{\mathcal{W}_v(M_0,\eta)\,|J^2m^3|_v^{n/2}}{|KH^2|^{n/2}_v}.
\end{split}
\end{equation*}
When $v$ is finite nonsplit,
\begin{equation*}
 \frac{\abs{\nrd(\delta)}_v^n}{\max(\norm{\gamma}_v,\abs{\nrd(\delta)}_v)^n}
= \frac{\abs{J^2m/C^2}_v^n}{\max(\norm{J\gamma'/C^2}_v,\abs{J^2m/C^2}_v)^n}
\asymp_v \frac{|J^2m^3|_v^{n/2}}{|KH^2|^{n/2}_v},
\end{equation*}
since $0\le v(K)\le v(m)=v(\nrd(\eta))\le 1$.
By Proposition \ref{prop:adelic bound} and the product formula,  we have
\begin{equation}
\label{global--product-I_v-bound}
I_{0}^\infty(\delta,\gamma)
\ll \frac{K^{n/2} H^n}
{J^n m^{3n/2}}
\mathcal{W}(M_0,\eta).
\end{equation}
where $\mathcal{W}(M_0,\eta) := \prod_{\textnormal{split }v\mid Jm} \mathcal{W}_v(M_0,\eta)$.  Say 
\begin{align*}
\norm{\gamma'} \in [R,2R),
\end{align*}
where $R$ is \emph{dyadic}, i.e.~$R\in \{2^t: t\in \ZZ\}$. Then by \eqref{eq:inftyboundp2} we have \begin{equation}
\begin{split}
\label{dyadic-I_infty-bound}
I_{0,\infty}\left(\frac{\delta}{X},\frac{\gamma}{X}\right)
&\ll_\epsilon
(1 + \norm{X\gamma/\nrd(\delta)})^{\epsilon-n} \\
&\asymp (1 + (XJR)/(J^2m))^{\epsilon-n}
=(1 + XR/(Jm))^{\epsilon-n}.
\end{split}
\end{equation}

Since $\norm{\delta}\ll X$, we have $\norm{\eta} \ll X/J$ and by definition $m \asymp \norm{\eta}^2$ as $D_\infty$ is nonsplit. Observe that $$1\le H\ll \norm{\gamma'} \asymp R \asymp \norm{\gamma/J}
\le X^\epsilon \norm{\delta/J} \asymp X^\epsilon \norm{\eta}.$$ In summary, $H\ll R\ll X^\epsilon \norm{\eta}$. Then combining \eqref{global--product-I_v-bound} and \eqref{dyadic-I_infty-bound}, we have
\begin{equation}
\begin{split}
\label{INEQ:first-E(X)-bound}
&\sum_{\gamma \in D^n - \{0\}} \sum_{\delta \in D^\times}
\left|I_{0}\left(\frac{\delta}{X},\frac{\gamma}{X}\right)\right|\\
&\ll_\epsilon \sum_{J\ge 1} \sum_{\norm{\eta} \ll X/J}
\sum_{\substack{H\mid Jm \\ K\mid m}}
\sum_{\substack{\textnormal{dyadic }R \\ H\ll R\ll X^\epsilon \norm{\eta}}}
\sum_{M_0\in S}
\sum_{\substack{\gamma' \in C^{-1}\Lambda^n \\ 0<\norm{\gamma'} \ll R}}
\frac{\mathcal{W}(M_0,\eta)  K^{n/2} H^n}
{(1 + XR/(Jm))^{n-\epsilon} J^n m^{3n/2}}, \\
&\ll_\epsilon X^{\epsilon}\sum_{J\ge 1} \sum_{\norm{\eta} \ll X/J}
\sum_{\substack{H\mid Jm \\ K\mid m}}
\sum_{\substack{\textnormal{dyadic }R \\ H\ll R\ll X^\epsilon \norm{\eta}}}
\sup_{M_0\in \OO_D}
\sum_{\substack{\gamma' \in C^{-1}\Lambda^n \\ 0<\norm{\gamma'} \ll R}}
\frac{ K^{n/2} H^n}
{(1 + XR/(Jm))^{n-\epsilon} J^n m^{3n/2}},
\end{split}
\end{equation}
where $\Lambda\defeq \Lambda(H,K,m,\eta,M_0)$,
and where in the final step we sum $\mathcal{W}(M_0,\eta)$ over $M_0\in S$
via Lemma~\ref{LEM:Z-equivalence-classes} at split places $v\mid m$
and via the bound $\#(\OO_{F_v}^\times \backslash (\OO_{D_v}/\OO_{D_v}\eta^\dagger)) \le \#(\OO_{D_v}/m\OO_{D_v}) = \abs{m}_v^{-4} \le q_v^4$ at nonsplit places $v\mid m$.

We postpone the proof of the following key lemma:
\begin{lem}
\label{lem:badcount}
Assume $F=\QQ$ and $D$ nonsplit.
Let $H,K,m\ge 1$ be integers.
Suppose $K\mid m \mid \nrd(\eta)$, where $\eta\in \OO_D\cap D^\times$ is primitive.
Let $M_0\in \OO_D$.
Then for all $R>0$, we have
\begin{equation*}
\begin{split}
&\#\{M \in \Lambda(H,K,m,\eta,M_0): \norm{M} \le R\} \\
&\ll 1 + \frac{R}{H}
+ \frac{(R/H)^2}{(K')^{1/2}}
+ \frac{(R/H)^3}{(K'm')^{1/2}}
+ \frac{(R/H)^4}{K'm'},
\end{split}
\end{equation*}
where $K' := K/\gcd(K,H)$ and $m' := m/\gcd(m,H)$.
\end{lem}

Assuming Lemma~\ref{lem:badcount}, we have
\begin{equation*}
\sup_{M_0\in \OO_D}
\sum_{\substack{\gamma' \in C^{-1}\Lambda^n \\ 0<\norm{\gamma'} \ll R}} 1
\ll \left(1 + \frac{R}{H}
+ \frac{(R/H)^2}{(K')^{1/2}}
+ \frac{(R/H)^3}{(K'm')^{1/2}}
+ \frac{(R/H)^4}{K'm'}\right)^n,
\end{equation*}
Since $R \ll X^\epsilon \norm{\eta}\asymp X^\epsilon m^{1/2},$
\begin{equation*}
\frac{(R/H)^4}{K'm'}
\ll X^\epsilon \frac{(R/H)^3}{(K'm')^{1/2}}
\ll X^{2\epsilon} \frac{(R/H)^2}{(K')^{1/2}},
\end{equation*}
by comparing consecutive terms using the bound $R/X^\epsilon \ll m^{1/2} \le H(m')^{1/2}$.
Since $1\ll \frac{R}{H}$ as well, the penultimate display now implies the simpler bound
\begin{equation}
\label{INEQ:simplified-geom-nums}
\sup_{M_0\in \OO_D}
\sum_{\substack{\gamma' \in C^{-1}\Lambda^n \\ 0<\norm{\gamma'} \ll R}} 1
\ll X^{2n\epsilon} \left(\frac{R}{H}
+ \frac{R^2}{H^2 (K')^{1/2}}\right)^n.
\end{equation}

A direct application of \eqref{INEQ:simplified-geom-nums},
for each $\Lambda$, gives (for any $R$ with $1\ll R\ll X^\epsilon \norm{\eta}$)
\begin{equation*}
\begin{split}
\sum_{K\mid m} \sum_{\substack{1\le H\ll R \\ H\mid Jm}} \sup_{M_0\in \OO_D}
\sum_{\substack{\gamma' \in  C^{-1}\Lambda^n \\ 0<\norm{\gamma'} \ll R}}
K^{n/2} H^n
&\ll X^{2n\epsilon} \sum_{K\mid m} \sum_{\substack{1\le H\ll R \\ H\mid Jm}}
K^{n/2} H^n
\left(\frac{R^n}{H^n}
+ \frac{R^{2n}}{H^{2n} (K')^{n/2}}\right) \\
&= X^{2n\epsilon} \sum_{K\mid m} \sum_{\substack{1\le H\ll R \\ H\mid Jm}}
\left(K^{n/2}R^n
+ \frac{\gcd(K,H)^{n/2}R^{2n}}{H^n}\right) \\
&\ll_\epsilon (Xm)^{2n\epsilon} \sum_{\substack{1\le H\ll R \\ H\mid Jm}}
\left(m^{n/2}R^n
+ \frac{\gcd(m,H)^{n/2}R^{2n}}{H^n}\right),
\end{split}
\end{equation*}
by the divisor bound on $m$.
Now, by the divisor bound on $Jm$, this is dominated by
\begin{equation*}
\begin{split}
& (XJm)^{4n\epsilon} \max_{\substack{1\le H\ll R \\ H\mid Jm}} \left(m^{n/2}R^n
+ \frac{\gcd(m,H)^{n/2}R^{2n}}{H^n}\right) \\
&\le (XJm)^{4n\epsilon} (m^{n/2}R^n + R^{2n}) \\
&\ll_\epsilon (XJm)^{5n\epsilon} m^{n/2}R^n,
\end{split}
\end{equation*}
since $R \ll_\epsilon X^\epsilon m^{1/2}$.
Substituting this into our previous bound \eqref{INEQ:first-E(X)-bound}, we get
\begin{align*}
&X^{O(\epsilon)}\sum_{J\ge 1} \sum_{\norm{\eta} \ll X/J}
\sum_{\substack{\textnormal{dyadic }R \\ 1\ll R\ll X^\epsilon \norm{\eta}}}
\frac{ m^{n/2}R^n}{(1 + XR/(Jm))^{n-\epsilon} J^n m^{3n/2}}\\
&\ll X^{-n+O(\epsilon)}\sum_{J\ge 1} \sum_{\norm{\eta} \ll X/J}
\sum_{\substack{\textnormal{dyadic }R \\ 1\ll R\ll X^\epsilon \norm{\eta}}} 1\\
&\ll X^{-n+O(\epsilon)}\sum_{J\ge 1} \sum_{\norm{\eta} \ll X/J}
1.
\end{align*}
However,
\begin{equation*}
\sum_{J\ge 1} \sum_{\norm{\eta} \ll X/J} 1
\ll \sum_{J\ge 1} (X/J)^4
\ll X^4,
\end{equation*}
so \eqref{eq:goal} holds as desired.
\qed

We are left with proving Lemma \ref{lem:badcount} which concerns the $\OO_F$-lattices 
\begin{equation}
\begin{split}
\label{EQN:define-key-lattice}
\Lambda(H,K,m,\eta,M_0) :=
\left\{M\in H\OO_D:
\begin{array}{ll}
(M-M^\dagger) \eta \in K\OO_D,\\
M\eta \in \OO_F M_0\eta + m\OO_D
\end{array}\right\}.
\end{split}
\end{equation}
for  $H,K,m\in \OO_F\cap F^\times$, primitive $\eta\in \OO_D\cap D^\times$ such that $K\mid m \mid \nrd(\eta)$, and $M_0\in \OO_D.$ We need to show for $R\ge 0$
\begin{equation}
\begin{split}
\label{INEQ:main-geom-numbers}
&\#\{M \in \Lambda(H,K,m,\eta,M_0): \norm{M} \le R\} \\
&\ll 1 + \frac{R}{H}
+ \frac{(R/H)^2}{(K')^{1/2}}
+ \frac{(R/H)^3}{(K'm')^{1/2}}
+ \frac{(R/H)^4}{K'm'},
\end{split}
\end{equation}
where $K' := K/\gcd(K,H)$ and $m' := m/\gcd(m,H)$. 

We remark that in the proof of Theorem \ref{THM:off-center-total-bound} we apply the lemma with $m=\nrd(\eta).$
However, allowing general $m\mid \nrd(\eta)$ turns out to be useful for the statement and proof of the lemma.
The proof of Lemma~\ref{lem:badcount} proceeds by techniques from Geometry of Numbers. While the messier case of $H\ne 1$  plays an important role in the proof of Theorem \ref{THM:off-center-total-bound}, we show that the general case follows from the case $H=1$ in Lemma \ref{LEM:remove-awkward-H-condition} below, so that one can actually neglect the awkward condition $M\in H\OO_D$ in their first reading.

\subsection{Proof of Lemma \ref{lem:badcount}}



Since $\eta$ is primitive, we have $0 \leq v(\nrd(\eta)) \leq 1$ if $D_{v}$ is nonsplit,
by Lemma~\ref{LEM:nonsplit-primitive-norms}.
Since $K\mid m\mid \nrd(\eta)$,
to prove Lemma \ref{lem:badcount} it suffices to prove the lemma in the special case where $v(K)=v(m)=0$ if $v$ is a finite place where $D_v$ is nonsplit.
In this case, we will need the following preliminary observations.

\begin{prop}
\label{PROP:prepgeom}
    Assume $F=\QQ$. Let $\eta\in\OO_{D}\cap D^\times$ be primitive and let $K\ge 1$ be such that $K\mid\nrd(\eta)$,
where $D$ is split at all primes dividing $K$. Then the following statements hold.
\begin{enumerate}
\item Suppose that $A\in\OO_{D}$ and $A\eta\equiv0\bmod{K\OO_D}$. Then $K\mid\nrd(A)$.
\item 
One has
$$
\#\{ \theta \in \OO_D/K\OO_D:\eta\theta\equiv0\bmod{K\OO_D}\}=K^2.$$
\item There exists $\theta\in\OO_D$ such that $\eta\theta\equiv0\bmod{K\OO_D}$ and $0<\norm{\theta}\ll K^{1/2}$.
\item There exists $\eta'\in\eta\OO_D+K\OO_D$ such that $0<\norm{\eta'}\ll K^{1/2}$.
\end{enumerate}
\end{prop}
\begin{proof}
For (1) and (2), we observe that, by factoring $K$ into prime powers, and using the assumption that $D$ splits at all prime divisors of $K$, an application of the Chinese Remainder Theorem
renders it sufficient to prove the following local statements, in both of which we assume that $\eta\in M_2(\OO_{F_v})\cap \GL_2(F_v)$ is primitive, and $K\in\OO_{F_v}$ is such that $K\mid\nrd(\eta)$:
\begin{enumerate}[label=(\arabic*')]
\item Suppose that $A\in M_2(\OO_{F_v})$ is such that $A\eta\equiv0\bmod{K}$. Then $K\mid\nrd(A)$.
\item The congruence
$$\eta\theta\equiv0\bmod{K}\text{,}$$
to be solved in $\theta\in M_2(\OO_{F_v})$, has exactly $K^2$ solutions modulo $K$.
\end{enumerate}
To prove (1'), by a Cartan decomposition argument we may assume that $\eta = \begin{psmatrix}\varpi^k & 0 \\ 0 & 1\end{psmatrix}$ for some integer $k\geq v(K)$.
If $A=\begin{psmatrix} s&u\\t&v\end{psmatrix}$, then the condition $A\eta\equiv0\bmod{K}$ becomes $K\mid u,v$, whence $K\mid sv-tu=\nrd(A)$.
To prove (2'), again by a Cartan decomposition argument we may assume that $\eta = \begin{psmatrix}\varpi^k & 0 \\ 0 & 1\end{psmatrix}$, whence if $\theta=\begin{psmatrix}
a&b\\c&d
\end{psmatrix}$, the congruence to be solved becomes $c\equiv d\equiv0\bmod{K}$, after which the result is clear. This finishes the proofs of (1') and (2'), and hence of (1) and (2).

To prove (3), consider the lattice
$$L=\{\theta\in\OO_D:\eta\theta\equiv0\bmod{K}\}\subseteq D\text{.}$$
 Note that $K\OO_D\subseteq L\subseteq \OO_D$, and by (2) we have $[L:K\OO_D]=K^2$. Therefore,
$$\covol(L) \asymp [\OO_D:L]
=\frac{[\OO_D:K\OO_D]}{[L:K\OO_D]}=\frac{K^4}{K^2}=K^2\text{.}$$
Hence, by Minkowski's Theorem, $L$ contains a nonzero vector of magnitude at most $O(\covol(L)^{1/4})=O(K^{1/2})$, which translates into the desired result.

Finally, to prove (4), consider the lattice
$$R=\eta\OO_D+K\OO_D\text{.}$$
Note that $K\OO_D\subseteq R\subseteq \OO_D$, and
\begin{equation}
\label{covolR}
\covol(R) \asymp [\OO_D:R]
=\frac{[\OO_D:K\OO_D]}{[R:K\OO_D]}
=\frac{K^4}{[R:K\OO_D]}\text{.}
\end{equation}
Consider the $\ZZ$-module homomorphism
$$\begin{array}{ccc}
\varphi: & \OO_D/K\OO_D & \longrightarrow \OO_D/K\OO_D \\
& \theta+K\OO_D & \longmapsto \eta\theta+K\OO_D
\end{array}$$
and observe that $[R:K\OO_D]=|\im(\varphi)|$.
Note that $|\ker(\varphi)|=K^2$ by (2).
By the First Isomorphism Theorem it follows that
$$[R:K\OO_D]=\frac{|\OO_D/K\OO_D|}{|\ker(\varphi)|}=\frac{K^4}{K^2}=K^2\text{.}$$
Inserting this into \eqref{covolR}, it follows that
$$\covol(R) \asymp \frac{K^4}{K^2}=K^2\text{.}$$
Now Minkowski's Theorem implies that $R$ contains a nonzero vector of magnitude at most $O(\covol(R)^{1/4})=O(K^{1/2})$, as desired.
\end{proof}

The following lemma will help us deal with the awkward condition $M\in H\OO_D$ in \eqref{EQN:define-key-lattice}.

\begin{lem}
\label{LEM:remove-awkward-H-condition}
Let $F,D,H,K,m,\eta,M_0,K',m'$ be as in Lemma~\ref{lem:badcount}.
Assume that $v(K)=v(m)=0$ if $D_v$ is nonsplit.
Then there exist
integers $K^\star,m^\star\ge 1$,
a vector $Y\in \OO_D$,
and a primitive vector $\sigma\in \OO_D\cap D^\times$,
such that the following hold:
\begin{enumerate}
\item $K^\star\mid m^\star\mid \nrd(\sigma)$;
\item $K^\star \asymp K'$ and $m^\star \asymp m'$;
\item $\norm{\sigma} \asymp (m^\star)^{1/2}$ and $\abs{\nrd(\sigma)} \asymp m^\star$;
\item there is an injection
\begin{align*}\phi: \Lambda(H,K,m,\eta,M_0) &\lto \Lambda(1,K^\star,m^\star,\sigma,Y)\\
M &\longmapsto M/H.
\end{align*}
\end{enumerate}
\end{lem}

\begin{proof}
For convenience, we define the map
\begin{equation*}
\Psi: H\OO_D \to \OO_D/m\OO_D,
\quad M \mapsto M\eta+m\OO_D.
\end{equation*}
For each $A\in \OO_D$, we also define the map
\begin{equation*}
T_A: \OO_F/m\OO_F \to \OO_D/m\OO_D,
\quad \lambda \mapsto \lambda A\eta+m\OO_D.
\end{equation*}
Observe that for any $M\in H\OO_D$, we have
\begin{equation*}
M\eta\in \OO_F A\eta+m\OO_D
\Leftrightarrow
\Psi(M)\in \im(T_A).
\end{equation*}

Both $\Psi$ and $T_A$ are $\OO_F$-module homomorphisms.
We will use these maps to gradually alter $M_0$ into a friendlier form.
Consider the ideal
\begin{equation*}
I := T_{M_0}^{-1}(\im(\Psi))
= \{x\in \OO_F/m\OO_F: T_{M_0}(x)\in \im(\Psi)\}
\le \OO_F/m\OO_F.
\end{equation*}
All ideals of $\OO_F/m\OO_F$ are principal,
so $I$ is generated by some $\lambda_0\in \OO_F$.
Let $M_1 := \lambda_0 M_0 \in \OO_D$.
Then
\begin{equation*}
\im(T_{M_1})
= T_{M_0}(I)
= \im(\Psi) \cap \im(T_{M_0}).
\end{equation*}
So for any $M\in H\OO_D$, we have
\begin{equation*}
\Psi(M)\in \im(T_{M_0})
\Leftrightarrow
\Psi(M)\in \im(\Psi) \cap \im(T_{M_0})
\Leftrightarrow
\Psi(M)\in \im(T_{M_1}).
\end{equation*}
Thus $\Lambda(H,K,m,\eta,M_0) = \Lambda(H,K,m,\eta,M_1)$.
Yet $T_{M_1}^{-1}(\im(\Psi)) = \OO_F/m\OO_F$, since $\im(T_{M_1}) \le \im(\Psi)$.
Replacing $M_0$ with $M_1$ if necessary, we may thus assume that $I=\OO_F/m\OO_F$.

Since $1\in I$, there exists $M_2\in H\OO_D$ such that
$T_{M_0}(1) = \Psi(M_2)$, i.e.~$$M_0\eta\equiv M_2\eta\bmod{m\OO_D}.$$
But then $\im(T_{M_0}) = \im(T_{M_2})$,
so $\Lambda(H,K,m,\eta,M_0) = \Lambda(H,K,m,\eta,M_2)$.
Replacing $M_0$ with $M_2$ if necessary, we may therefore assume that $M_0\in H\OO_D$.
Let $M_3 := M_0/H\in \OO_D$.

Since $\lcm(H,K) = HK'$ and $\lcm(H,m) = Hm'$, we have
$$\Lambda(H,K,m,\eta,M_0) = \Lambda(H,HK',Hm',\eta,M_0)$$
by \eqref{EQN:define-key-lattice},
because $H\mid M_0$.
Now, dividing the conditions in \eqref{EQN:define-key-lattice} through by $H$,
we see that the formula $M\mapsto M/H$ maps $\Lambda(H,K,m,\eta,M_0)$ into $\Lambda(1,K',m',\eta,M_3)$.

Clearly $K'\mid m'\mid \nrd(\eta)$, since $K\mid m\mid \nrd(\eta)$.
By Proposition \ref{PROP:prepgeom}(4) with ``$(\eta,K)\defeq (\eta,m')$'', there exists a nonzero vector $\eta' \in \eta\OO_D + m'\OO_D$ such that $\norm{\eta'} \ll (m')^{1/2}$.
Since $D$ is nonsplit, we have $\eta'\in D^\times$, i.e.~$\nrd(\eta')\ne 0$.
Say $\eta' \in \eta\alpha + m'\OO_D$, where $\alpha\in \OO_D$.
Then multiplying the conditions in \eqref{EQN:define-key-lattice} by $\alpha$ on the right, we find that
\begin{equation*}
\Lambda(1,K',m',\eta,M_3) \subseteq \Lambda(1,K',m',\eta',M_3).
\end{equation*}
Moreover, $\nrd(\eta') \equiv \nrd(\eta\alpha) \equiv 0 \bmod{m'}$, since $m'\mid \nrd(\eta)$.
In particular, $\abs{\nrd(\eta')} \asymp m'$ and $\norm{\eta'} \asymp (m')^{1/2}$, since $\norm{\eta'} \ll (m')^{1/2}$.

We are almost done, but $\eta'\in \OO_D\cap D^\times$ need not be primitive.
Write $\eta' = z\sigma$, where $z\ge 1$ and $\sigma\in \OO_D\cap D^\times$ is primitive.
Then $\gcd(m',z)\mid \eta\alpha$, since $\eta' \in \eta\alpha + m'\OO_D$.
So by Proposition~\ref{PROP:prepgeom}(1) with ``$(\eta,K,A)\defeq (\eta^\dagger,\gcd(m',z),\alpha^\dagger)$'', we have $\gcd(m',z)\mid \nrd(\alpha)$.
Thus
\begin{equation*}
\begin{split}
\nrd(\eta') = (\eta') (\eta')^\dagger
&\in (\eta\alpha + m'\OO_D) ((\eta\alpha)^\dagger + m'\OO_D) \\
&\subseteq \nrd(\eta\alpha) + \gcd(m',z)m'\OO_D
= \gcd(m',z)m'\OO_D,
\end{split}
\end{equation*}
since $m'\mid \nrd(\eta)$.
But $\nrd(\eta') = z^2\nrd(\sigma)$, so it follows that
\begin{equation*}
\gcd(m',z)m'
\mid \gcd(z^2,\gcd(m',z)m') \nrd(\sigma)
= \gcd(m',z)^2 \nrd(\sigma),
\end{equation*}
whence $m' \mid \gcd(m',z) \nrd(\sigma) \mid z\, \nrd(\sigma)$.
Yet $\norm{\sigma} = \norm{\eta'}/z \asymp (m')^{1/2}/z$, so we conclude that
\begin{equation*}
m' \le z\, \abs{\nrd(\sigma)}
\ll z \norm{\sigma}^2
\asymp z(m'/z^2) = m'/z.
\end{equation*}
Thus $z\ll 1$.
In other words, $\eta'$ is either primitive, or nearly primitive.

Let $K^\star\defeq K'/\gcd(z^2,K')$ and $m^\star\defeq m'/\gcd(z^2,m')$.
Since $z\ll 1$, we have
\begin{equation*}
K^\star\asymp K',
\quad
\abs{\nrd(\sigma)} \asymp \abs{\nrd(\eta')} \asymp m' \asymp m^\star,
\quad
\norm{\sigma} \asymp \norm{\eta'} \asymp (m')^{1/2} \asymp (m^\star)^{1/2}.
\end{equation*}
Also, since $K'\mid m'\mid \nrd(\eta') = z^2\nrd(\sigma)$,
we have $K'\mid m'\mid \gcd(z^2,m')\nrd(\sigma)$,
whence $K^\star\mid m^\star\mid \nrd(\sigma)$.
Moreover, since $\eta' = z\sigma$, it is clear that
\begin{equation*}
\Lambda(1,K',m',\eta',M_3) \subseteq \Lambda(1,K^\star,m^\star,\sigma,M_3),
\end{equation*}
by \eqref{EQN:define-key-lattice}.
Let $Y := M_3 \in \OO_D$, and take $\phi$ to be the composition
\begin{equation*}
\Lambda(H,K,m,\eta,M_0)
\to \Lambda(1,K',m',\eta,M_3)
\to \Lambda(1,K',m',\eta',M_3)
\to \Lambda(1,K^\star,m^\star,\sigma,M_3),
\end{equation*}
where the first arrow is $M\mapsto M/H$
and the second and third arrows are inclusion.
This finishes the job, with conditions (1)--(4) all satisfied.
\end{proof}

\begin{proof}
[Proof of Lemma~\ref{lem:badcount}]
In the notation of Lemma~\ref{LEM:remove-awkward-H-condition}, we have
\begin{equation*}
\#\{M \in \Lambda(H,K,m,\eta,M_0): \norm{M} \le R\}
\le \#\{M \in \Lambda(1,K^\star,m^\star,\sigma,Y): \norm{M} \le R/H\}.
\end{equation*}
If \eqref{INEQ:main-geom-numbers} holds for $\Lambda(1,K^\star,m^\star,\sigma,Y)$, then the right-hand side above is
\begin{equation*}
\ll 1 + \frac{R}{H}
+ \frac{(R/H)^2}{(K^\star)^{1/2}}
+ \frac{(R/H)^3}{(K^\star m^\star)^{1/2}}
+ \frac{(R/H)^4}{K^\star m^\star}.
\end{equation*}
This would imply \eqref{INEQ:main-geom-numbers} for $\Lambda(H,K,m,\eta,M_0)$, because $K^\star \asymp K'$ and $m^\star \asymp m'$.

Therefore, we are reduced to proving \eqref{INEQ:main-geom-numbers} for $\Lambda(1,K^\star,m^\star,\sigma,Y)$.
In other words, by Lemma~\ref{LEM:remove-awkward-H-condition},
it suffices to prove Lemma~\ref{lem:badcount} under the conditions
\begin{equation}
H=1,
\quad
\norm{\eta}\asymp m^{1/2},
\quad
\abs{\nrd(\eta)} \asymp m,
\end{equation}
which we henceforth assume.
Then $(K',m') = (K,m)$, since $H=1$.

For the rest of the proof, we abbreviate $\Lambda(H,K,m,\eta,M_0)$ to $\Lambda$.
Trivially,
\begin{equation}
\label{trivial-lattice-bounds}
m\OO_D \subseteq \Lambda \subseteq \OO_D,
\end{equation}
by \eqref{EQN:define-key-lattice} and the condition $K\mid m$.
In particular, $\Lambda$ is a full lattice in $D_\infty\cong \RR^4$.
We claim that its successive minima $\lambda_i>0$ (for $1\le i\le 4$) satisfy
\begin{enumerate}
\item $\lambda_1 \gg 1$;
\item $\lambda_2 \gg K^{1/2}$;
\item $\lambda_4 \ll (Km)^{1/2}$;
\item $\lambda_1\lambda_2\lambda_3\lambda_4
\asymp \covol(\Lambda)
\gg Km$.
\end{enumerate}
Assuming the claim, the inequality 
~\eqref{INEQ:main-geom-numbers}
follows from
the classical fact,
recorded in \cite[Lemma~3.6 in the published version; Lemma~3.5 in {\tt arXiv:2006.02356v1}]{modernlattice}
and based on \cite[Lemma~2]{schmidt},
that if $\Lambda\subseteq\mathbb{R}^n$ is a full lattice and $\lambda_1\leq\cdots\leq\lambda_n$ denote its successive minima, and $S$ denotes the ball of radius $R$ centered at the origin, then
\begin{equation}\label{latticebound}
|S\cap\Lambda|
\ll_n
\sum_{j=0}^n\frac{R^j}{\lambda_1\cdots\lambda_j}.
\end{equation}
Indeed, $\lambda_1\lambda_2\lambda_3\gg (Km)^{1/2}$ by (3) and (4).
Thus we are reduced to proving the claim.

The proof of (1) is immediate from the containment $\Lambda\subseteq \OO_D$.
For (2) we begin by observing that if $U\neq0$ and $U\eta\equiv0\bmod{K}$, then $\norm{U}\gg K^{1/2}$: this follows by Proposition \ref{PROP:prepgeom}(1), which implies that under these assumptions $K\mid\nrd(U)\ll\norm{U}^2$, where $\nrd(U)\ne 0$ because $D$ is nonsplit.
Now if $M\in\Lambda$ is not a scalar multiple of the identity, then $M-M^\dagger\neq0$ and hence, by the previous observation, $\norm{M-M^\dagger}\gg K^{1/2}$.
Since $\norm{M-M^\dagger}\ll\norm{M}$, it follows that for $M\in\Lambda$ not a scalar multiple of the identity we have $\norm{M}\gg K^{1/2}$, whence $\lambda_2\gg K^{1/2}$.

We moreover have $\lambda_1\lambda_2\lambda_3\lambda_4\asymp\covol(\Lambda)$,
by Minkowski's second theorem, as is conveniently recorded in \cite[(3.5)]{modernlattice}, say.
By \eqref{trivial-lattice-bounds}, we may write
\begin{equation}
\label{vol}
\covol(\Lambda) \asymp [\OO_D:\Lambda]
= \frac{[\OO_D:m\OO_D]}{[\Lambda:m\OO_D]}
= \frac{m^4}{[\Lambda:m\OO_D]}
= \frac{m^4}{[\tau\Lambda:\nrd(\eta)\OO_D]},
\end{equation}
where $\tau \defeq \nrd(\eta)/m\in \ZZ$.
By the final condition in \eqref{EQN:define-key-lattice},
elements of
\begin{equation*}
\tau\Lambda/\nrd(\eta)\OO_D
= \Lambda(\tau, \tau K, \tau m, \eta, \tau M_0)/\nrd(\eta)\OO_D
\end{equation*}
can be written in the form $M=aM_0+A\eta^\dagger$, where $a\in\tau\OO_F/\nrd(\eta)\OO_F$ and $A\in\OO_D/\OO_D\eta$.
The condition $(M-M^\dagger)\eta\equiv0\bmod{\tau K}$ implies
$(M-M^\dagger)\eta\equiv0\bmod{K}$, which translates into
\begin{equation}
\label{aA-cong}
a(M_0-M_0^\dagger)\eta - \eta A^\dagger\eta\equiv0\bmod{K}.
\end{equation}
Note that $\eta A^\dagger\eta\bmod{K}$ is well-defined for $A\in\OO_D/\OO_D\eta$, since $K\mid\nrd(\eta)=\eta\eta^\dagger$.
Since $\OO_{D_v}^\times \eta \OO_{D_v}^\times=\OO_{D_v}^\times \eta^\dagger \OO_{D_v}^\times$ for all $v\mid K$,
we find that by Lemma~\ref{LEM:delta-module-arithmetic},
the kernel of $A\mapsto\eta A^\dagger\eta\bmod{K}$ has index $K$ in $\OO_D$ (or in $\OO_D/\OO_D\eta$).
This implies that
once $a$ is specified,
the congruence \eqref{aA-cong} has either $0$ or $\frac{|\OO_D/\OO_D\eta|}{K}$ solutions $A\in \OO_D/\OO_D\eta$.
On summing over $a\in \tau\OO_F/\nrd(\eta)\OO_F$, we get
$$[\tau\Lambda:\nrd(\eta)\OO_D]
\le \sum_{a\in \tau\OO_F/\nrd(\eta)\OO_F} \frac{|\OO_D/\OO_D\eta|}{K}
\ll \frac{m^3}{K},$$
since $|\tau\OO_F/\nrd(\eta)\OO_F| = m$ and $|\OO_D/\OO_D\eta| \asymp \nrd(\eta)^2=O(m^2)$.
Hence by \eqref{vol} we obtain $\covol(\Lambda)\gg mK$, which yields (4).

Finally, we establish (3), for which it suffices to show that $\Lambda$ contains a full sublattice spanned by elements of magnitude $O((Km)^{1/2})$.
For this we observe that there exists $\theta\in\OO_D$ such that $0<\norm{\theta}\ll K^{1/2}$ and $\eta\theta\equiv0\bmod{K}$, by Proposition~\ref{PROP:prepgeom}(3).
We necessarily have $\theta\in D^\times$, because $D$ is nonsplit.
Now consider the lattice
$$\Lambda' \defeq \{B\theta^\dagger\eta^\dagger:B\in\OO_D\}
\qquad \text{(a full lattice in $D_\infty$, because $\theta,\eta\in D^\times$).}$$
Clearly $\Lambda'\subseteq\Lambda$,
by \eqref{EQN:define-key-lattice},
since $K\mid \nrd(\eta)$ and $m\mid \nrd(\eta)$.
Moreover, if $e_1,e_2,e_3,e_4$ denotes an integral basis of $\OO_D$, with each $e_i$ having magnitude $O(1)$, then $e_1\theta^\dagger\eta^\dagger,e_2\theta^\dagger\eta^\dagger,e_3\theta^\dagger\eta^\dagger,e_4\theta^\dagger\eta^\dagger$ is an integral basis of $\Lambda'$ formed by vectors of magnitude $O(\norm{\theta}\norm{\eta})=O((Km)^{1/2})$.
\end{proof}

\begin{rem}
In key ranges of parameters, we expect Lemma~\ref{lem:badcount} to be optimal for some $M_0$, such as the scalar $1$, but not necessarily for all $M_0$.
For example, one may usually expect to have $1\notin \Lambda$, in which case our lower bound for $\lambda_1$ could be improved.
It is an interesting question to try to improve our lattice bounds, at least for typical choices of $M_0$.
Indeed, this is one plausible route to decreasing the $n$ required in Theorem~\ref{thm:simplest-case}.
\end{rem}

\bibliography{refs}{}
\bibliographystyle{alpha}

\end{document}


\title{Nonabelian circle method supplements}

\author{Nuno Arala}
\address{Mathematics Institute, Zeeman Building, University of Warwick, Coventry CV4 7AL, United Kingdom}
\email{Nuno.Arala-Santos@warwick.ac.uk}

\author{Jiaqi Hou}
\address{Department of Mathematics\\
University of Wisconsin-Madison\\
Madison, WI 53706}
\email{jhou39@wisc.edu}

\author{Chun-Hsien Hsu}
\address{Department of Mathematics\\
Duke University\\
Durham, NC 27708}
\email{simonhsu@math.duke.edu}

\author{Huajie Li}
\address{Department of Mathematics\\
Johns Hopkins University\\
Baltimore, MD 21218}
\email{hli213@jhu.edu}

\author{Victor Y. Wang}
\address{Department of Mathematics\\
IST Austria}
\email{victor.wang@ist.ac.at}

\begin{abstract}
We discuss supplementary topics of classical interest concerning quaternion zeros of the equation $\gamma_1^2 \pm \dots \pm \gamma_n^2 = 0$.
None of the material in this paper is necessary for our previous paper on the topic, joint with Getz.
\end{abstract}

\maketitle


Notation is as in our previous paper \cite{Main},
joint with Getz,
whom we thank for helpful discussions.
Appendix~\ref{SEC:birch-singular-locus} computes some singular loci associated to the system
\begin{equation}
\label{EQN:explicit-system}
\sum_{1\le i\le n} (w_i^2 - x_i^2-y_i^2-z_i^2)= \sum_{1\le i\le n} w_ix_i = \sum_{1\le i\le n} w_iy_i = \sum_{1\le i\le n} w_iz_i = 0.
\end{equation}
Appendix~\ref{SEC:prime-case-formulas} gives integral formulas that could be useful for Kloosterman-type averaging;
it also provides an extra sanity check for the local integral estimates proven in \cite[\S~8]{Main}.
Appendix~\ref{SEC:abelian} sketches an abelian approach to the $\gamma=0$ calculation in some simple cases,
and explains why this simple approach fails in general.

\appendix 

\section{Birch singular locus}
\label{SEC:birch-singular-locus}

Recall that \eqref{EQN:explicit-system} in $\ZZ[1/2]$
is equivalent to
$\gamma_1^2+\dots+\gamma_n^2=0$ in $\mathbb{H} \otimes \ZZ[1/2]$.

\begin{prop}
\label{PROP:estimate-singular-locus}
Assume $n\ge 6$.
Then the system \eqref{EQN:explicit-system}
defines an affine complete intersection of dimension $4n-4$ over $\CC$,
with a singular locus of dimension $2n+1$.
Moreover, the Birch singular locus has dimension $2n+2$.
\end{prop}

\begin{proof}
Write $\mathbb{H} \otimes \CC \cong M_2(\CC)$.
For $A\in M_2(\CC)^n$, let $P(A)\defeq A_1^2+\dots+A_n^2$.
For any $W\in M_2(\CC)$, view $\trd(WP)$ as a polynomial in the $4n$ scalar entries of $A$.
By definition, the Birch singular locus
consists of $A\in M_2(\CC)^n$ for which there exists $W\in M_2(\CC) - \{0\}$ with $\grad{\trd(WP)}=0$.
Moreover, by the Jacobian criterion, the (non-Birch) singular locus
is the intersection of the Birch singular locus with the zero locus $P(A)=0$.
But $P$ is diagonal, so
\begin{equation*}
\grad{\trd(WP)} = (\grad{\trd(WA_1^2)},\dots,\grad{\trd(WA_n^2)}).
\end{equation*}
Therefore, the Birch singular locus equals $\bigcup_{W\in M_2(\CC) - \{0\}} L(W)^n$, where
\begin{equation*}
L(W) \defeq \{A\in M_2(\CC): \grad{\trd(WA^2)}=0\}
= \{A\in M_2(\CC): WA+AW=0\}.
\end{equation*}
The relation $A\in L(W)$ is symmetric, i.e.~$A\in L(W)\Leftrightarrow W\in L(A)$.

We claim that $\dim(L(W)) = \max(2\cdot \one_{\trd(W)=0}, \one_{\det(W)=0})$ for all $W\ne 0$.
We check this by casework on Jordan normal form, up to scaling, as follows.
If $W = \begin{psmatrix}\lambda & 0 \\ 0 & 1\end{psmatrix}$, then $L(W) = \{A: \lambda a_{11} = (\lambda+1)a_{12} = (1+\lambda)a_{21} = 2a_{22} = 0\}$, which has dimension $2\cdot \one_{\lambda=-1} + \one_{\lambda=0}$.
If $W = \begin{psmatrix}1 & k \\ 0 & 1\end{psmatrix}$, then $L(W) = \{A: 2a_{11}+ka_{21} = 2a_{12}+k(a_{11}+a_{22}) = 2a_{21} = 2a_{22}+ka_{21} = 0\} = 0$.
If $W = \begin{psmatrix}0 & 1 \\ 0 & 0\end{psmatrix}$, then $L(W) = \{A: a_{21} = a_{22}+a_{11} = 0\}$, which has dimension $2$.

For the eventual (non-Birch) singular locus estimate, we note that the previous paragraph also shows that if $\trd(W)=0$, then $L(W)\subseteq \trd^{-1}(0)$,
even if $\det(W)=0$.
Moreover, if $\trd(W)=0$, then $\det\vert_{L(W)} \ne 0$, i.e.~there exists $A\in L(W)$ with $\det(A)\ne 0$.

Next, we claim that if $W_1,W_2\ne 0$ and $L(W_1)\ne L(W_2)$, then $\dim(L(W_1)\cap L(W_2)) \le 1$.
Indeed, $\dim(L(W_i))\le 2$, so if $\dim(L(W_1)\cap L(W_2)) \ge 2$, then necessarily $L(W_1)=L(W_2)$.

Finally, we claim that if $W_1,W_2\ne 0$ and $L(W_1)=L(W_2)\cong \CC^2$, then $W_1\in \CC^\times W_2$.
Indeed, suppose not.
If $L(W_1)$ is spanned by $A_1,A_2\ne 0$, then $W_1,W_2\in L(A_1)\cap L(A_2)$, whence $L(A_1)=L(A_2) = \CC W_1 + \CC W_2$.
But $\dim(L(W_1))=\dim(L(A_1))=2$ implies $\trd(W_1)=\trd(A_1)=0$, so $L(W_1),L(A_1)\subseteq \trd^{-1}(0)\cong \CC^3$, whence $L(W_1)\cap L(A_1)\ni B\ne 0$, say.
Then $W_1,W_2,A_1,A_2\in L(B)$.
Since $\dim(L(B))\le 2$, we conclude that $\CC A_1+\CC A_2 = L(B) = \CC W_1+\CC W_2$.
By definition of $A_1,A_2$, this means $W_1,W_2\in L(W_1)=L(W_2)$.
%
%
So for all $A\in L(W_1)$, we have $A^2=0$, whence $\det(A)=0$.
This contradicts a previous observation.

It follows that if $W_1\ne W_2\in \mathbb{P}^3$, then $\dim(L(W_1)\cap L(W_2)) \le 1$.
So $\bigcup_{W_1\ne W_2\in \mathbb{P}^3} (L(W_1)^n\cap L(W_2)^n)$ has dimension $\le 3\cdot 2 + n \le 2n$, since $n\ge 6$.
Therefore, the Birch singular locus has dimension exactly
$\dim(V_{\mathbb{P}^3}(\trd)) + 2n = 2n+2$,
by inclusion-exclusion.
Moreover, the intersection of the (non-Birch) singular locus with $L(W)^n$ is precisely the set $\{A\in L(W)^n: \sum_i A_i^2=0\} = \{A\in L(W)^n: \sum_i \det(A_i)=0\}$, whenever $\trd(W)=0$.
So the singular locus has dimension exactly
$\dim(V_{\mathbb{P}^3}(\trd)) + (2n-1) = 2n+1$,
by inclusion-exclusion.
The full system \eqref{EQN:explicit-system} is the union of its singular part, of dimension $2n+1$, and a smooth part, of dimension $4n-4$.
Thus each irreducible component must have dimension exactly $4n-4$.
\end{proof}

\quash{
\section{Van der Corput lemma over unramified extensions}

Let $K$ be an unramified extension of $\QQ_p$ with residue field $\FF_q$. Let $\psi:K\to \CC^\times$ be an additive character such that $\psi$ is trivial on $\OO_K$ but nontrivial on $p^{-1}\OO_K$, and the Haar measure on $K$ is chosen so that $\OO_K$ has volume $1$. Let $\calT_n=\calT_n(K)$ be the Tate algebra (the ring of convergent power series) in $n$ indeterminates over $K$. 

\begin{defn}
We say $f\in \calT_1$ is \textbf{special} if writing $f(x)=\sum_{i=0}^\infty a_ix^i$, one has $a_1\neq 0$ and $|a_j|<|a_1|$ for all $j>1$. More generally, if there is an integer $r\in \ZZ_{\ge 0}$ such that for all $b\in \OO_K$ and $0\neq c\in p^r\OO_K,$ the series
\begin{align*}
    f_{b,c}(x):=\frac{1}{c}f(b+cx) 
\end{align*}
is special, then we say $f$ is \textbf{$r$-special}.
\end{defn}

\begin{rem}
One checks easily $f$ is special iff $f$ is $0$-special.
\end{rem}

\begin{lem}\label{lem:basic}
Let $f\in \calT_1$. 
\begin{enumerate}
    \item If $|f'(x)|\ge 1$ for all $x\in \OO_K$, then $f$ is $r$-special for some $r$ satisfying
    \begin{align*}
        q^{r-1}\le \norm{f-f(0)}.
    \end{align*}
    Moreover, for any $r$ such that $f$ is $r$-special we have $\norm{f_{b,c}-f_{b,c}(0)}\ge 1$ for all $b\in \OO_K$ and $0\neq c\in p^r\OO_K$.
    \item If $f'$ is special and nonvanishing on $\OO_K,$ then $f$ is $1$-special. If moreover $\mathrm{deg}(f)<p$, then $f$ is special.
\end{enumerate}
\end{lem}

\begin{proof}
The first statement is \cite[Lemma 2.4]{Cluckers:vdC} and the second statement follows from  \cite[Corollary 2.9]{Cluckers:vdC} as $K$ is an unramified extension over $\QQ_p$.
\end{proof}

\begin{prop}\label{prop:unrvdC}
    Let $f\in \calT_1$. Suppose for some $k\in \ZZ_{>0}$, $|f^{(k)}(x)|\ge 1$ for all $x\in \OO_K$. Then for all $y\in K^\times$, 
\begin{align*}
    \Bigg|\int_{\OO_K }\psi(yf(x)) dx \Bigg| \ll_{k} q^{k-1}q^{r}|y|^{-1/k}.
\end{align*}
where $r$ is any number such that $f^{(k-1)}$ is $r$-special. In particular, 
\begin{align*}
    \Bigg|\int_{\OO_K }\psi(yf(x)) dx \Bigg| \ll_{k} q^{k}\norm{f-f(0)}|y|^{-1/k}.
\end{align*}
If $q$ is large in terms of $k$, implicit constants can be chosen to be $2^{2(k-1)}$.
\end{prop}

\begin{proof}
The case $k=1$ follows from the proof of \cite[Proposition 3.3]{Cluckers:vdC}. Thus we assume $k\ge 2$. Since $|f^{(k)}(x)|\ge 1 $ for all $x\in \OO_K,$ Lemma \ref{lem:basic} (1) implies $f^{(k-1)}$ is $r$-special for some $r$. Observe that the assertion is trivial unless $|y|\ge q^{k(k-1+r)},$ which we assume henceforth. As argued in \textit{loc.~ cit.}, it suffices to prove the claimed bound is an upper bound of
\begin{align*}
    \Bigg|\int_{\OO_K} \psi(p^{kr} y g(x))dx\Bigg|
\end{align*}
for any $g\in \calT_1$ such that $g^{(k-1)}$ is special $|g^{(k-1)}(x)|\ge |x|$ for all $x\in \OO_K$. 

Choose $m\in \ZZ_{\ge 1}$. Partition $\OO_K$ into the ball $B_0=p^m\OO_K$ and $n$ balls $B_j=d_j+n_j\OO_K$ where $d_j\not\in B_0$ and $n_j=pd_j$. One easily computes $n=(q-1)m$ and
\begin{align}\label{eq:djsum}
    \sum_{j=1}^{(q-1)m} |d_j|^{-\frac{1}{k-1}}\le q^{\frac{m}{k-1}}\frac{q-1}{q^{1/(k-1)}-1}\ll_k q^{\frac{m+k-2}{k-1}}
\end{align}
where the implicit constant can be chosen to be $2$ if $q$ is large in terms of $k$. By the argument in \textit{loc.~ cit.}, we have
\begin{align*}
    \Bigg|\int_{\OO_K} \psi(p^{kr} y g(x))dx\Bigg|\le q^{-m}+\sum_{j=1}^{(q-1)m} |n_j| \Bigg|\int_{\OO_K} \psi(p^{kr} y g_{d_j,n_j}(x))dx\Bigg|.
\end{align*}
Moreover $|g_{d_j,n_j}^{(k-1)}(x)|\ge |n_j^{k-1}d_j|>0$ for all $x\in \OO_K$, and $g_{d_j,n_j}^{(k-2)}$ is $1$-special.

We now proceed by induction on $k$. Assume the assertion holds for $k-1$. Then by Lemma \ref{lem:basic} (2)
\begin{align*}
    |n_j| \Bigg|\int_{\OO_K} \psi(p^{kr} y g_{d_j,n_j}(x))dx\Bigg|\ll_k q^{k-1}|d_jp^{kr}y|^{-\frac{1}{k-1}}
\end{align*}
and hence
\begin{align*}
     \Bigg|\int_{\OO_K} \psi(p^{kr} y g(x))dx\Bigg|\ll_k q^{-m}+q^{k-1}|y|^{-\frac{1}{k-1}}q^{\frac{kr+m+k-2}{k-1}}.
\end{align*}
Since $|y|\ge q^{k(k-1+r)},$ we can choose $m\ge 1$ such that 
\begin{align*}
    |y|^{1/k}\le q^{\frac{(k-1)^2+k-2}{k}}q^{r+m}\le q^{\frac{k-1}{k}}|y|^{1/k},
\end{align*}
and we have
\begin{align*}
    &q^{-m}+q^{k-1}|y|^{-\frac{1}{k-1}}q^{\frac{kr+m+k-2}{k-1}}=q^{-m}+q^r|y|^{-\frac{1}{k-1}}q^{\frac{(k-1)^2+r+m+k-2}{k-1}}\\
    &\le 2q^r|y|^{-\frac{1}{k}}q^{\frac{(k-1)^2+k-2}{k}+\frac{1}{k}}=2q^{k-1}q^r|y|^{-\frac{1}{k}}.
\end{align*}
\end{proof}

\begin{prop}\label{prop:unrvdCm}
Let $f\in \calT_n$. Suppose there exists a nonzero $\alpha\in \ZZ_{\ge 0}^n$ such that $|\partial^\alpha f(x)|\ge 1$ for all $x\in \OO_K^n$. Then for all $y\in K^\times$, 
\begin{align*}
    \Bigg|\int_{\OO_K^n }\psi(yf(x)) dx \Bigg| \ll_{k} q^{2k+1}\norm{f-f(0)}^{k+2}|y|^{-1/k} 
\end{align*}
where $k=|\alpha|$.
\end{prop}

\begin{proof}
The assertion follows from the proof of \cite[Proposition 3.6]{Cluckers:vdC} together with explicit bound obtained in Proposition \ref{prop:unrvdC}.
\end{proof}
}











\section{Prime-case formulas}
\label{SEC:prime-case-formulas}


Let $F$ be a non-Archimedean local field of characteristic zero and $\FF_q \defeq \OO_F/\varpi\OO_F$ be its residue field.
Assume $D\cong M_2(F)$ is split.
Let $f=\one_{\OO_D^n}$ and $\Phi=\one_{\OO_D^2}$.
Let $\delta\in \OO_D\cap D^\times$ and $\gamma\in \OO_D^n$.
For $\alpha,\beta\in \OO_{D}^\times$ we have the symmetry
\begin{equation}
\label{EQN:partial-Cartan-symmetry}
I_0(\delta,\gamma) = I_0(\alpha\delta\beta,\alpha\gamma\alpha^{-1}).
\end{equation}
For the rest of Appendix~\ref{SEC:prime-case-formulas},
suppose $v(\nrd(\delta))=1$, i.e.~that $\nrd(\delta)$ is \emph{prime}.
Let
\begin{equation*}
\begin{split}
S_2(\delta)&\defeq \#\{Y\in (\OO_D/\varpi\OO_D)^n: \delta^\dagger P(Y) \equiv 0\bmod{\varpi}\}, \\
S_3(\delta,\gamma) &\defeq \#\{Y\in (\OO_D/\varpi\OO_D)^n: \delta^\dagger P(Y) \equiv \trd(\gamma\cdot Y)\equiv 0\bmod{\varpi}\}.
\end{split}
\end{equation*}
For notational simplicity, we assume $\coeff_1=\dots=\coeff_n=1$, so $P(\gamma)=\sum_{1\le i\le n} \gamma_i^2$.
For $s\in \FF_q^n$, let
\begin{equation*}
X_2(s) \defeq \#\{a\in \FF_q^n: \sum_{1\le i\le n} a_i^2 = \sum_{1\le i\le n} s_ia_i = 0\}.
\end{equation*}
Then $X_2(0) = q^{n-1}+O(q^{\lfloor{n/2}\rfloor})$,
and $X_2(s) = q^{n-2} + O(q^{1+\lfloor{(n-2)/2}\rfloor})
= q^{n-2} + O(q^{\lfloor{n/2}\rfloor})$ if $s\ne 0$.

\begin{prop}
Let $\tilde{\gamma} \defeq \gamma\bmod{\nrd(\delta)} \in M_2(\FF_q)^n$.
In the setting above, where $\nrd(\delta)$ is a uniformizer, the following hold:
\begin{enumerate}

\item $q^{4n}(1-q^{-1}) I_0(\delta,\gamma) = S_3(\delta,\gamma) - q^{-1} S_2(\delta)$.


\item Suppose $\delta = \begin{psmatrix}\varpi & 0 \\ 0 & 1\end{psmatrix}$.
Write $\tilde{\gamma} = \begin{psmatrix}s & u \\ t & v\end{psmatrix}$, where $s,t,u,v\in \FF_q^n$.
Then
\begin{equation*}
\begin{split}
S_3(\delta,\gamma) &= \bm{1}_{(u,v)\ne 0}[q^{2n-1} X_2(0) + q^{3n-3} (q^n - q)]
+ \bm{1}_{u\ne 0} \bm{1}_{v\notin \FF_q u} q^{3n-3} (q-1) \\
&+ \bm{1}_{u\ne 0} \bm{1}_{v\in \FF_q u} q^{2n-2} (\#\{a,\lambda: \sum_{1\le i\le n} a_i^2 = \sum_{1\le i\le n} (s_i-v_i)a_i\lambda+t_iu_i\lambda^2\} - X_2(0)) \\
&+\bm{1}_{u=0} \bm{1}_{v\ne 0} q^{2n-2} (\#\{a,\lambda: \sum_{1\le i\le n} (s_i-v_i)a_i+t_iv_i\lambda = 0\} - \#\{a: \sum_{1\le i\le n} (s_i-v_i)a_i = 0\}) \\
&+ \bm{1}_{u=v=0}[q^{2n} X_2(s)
+ q^{2n-2} (\#\{a,b: \sum_{1\le i\le n} s_ia_i+t_ib_i = 0\} - \#\{a: \sum_{1\le i\le n} s_ia_i = 0\})],
\end{split}
\end{equation*}
where the variables in the point counts range over $a,b\in \FF_q^n$ and $\lambda\in \FF_q$.

\item In (2), we have the following useful formulas.
First,
\begin{equation*}
\begin{split}
\bm{1}_{(u,v)\ne 0} = \bm{1}_{\tilde{\gamma}\delta\ne 0},
\quad
\bm{1}_{u=v=0}
= \one_{\tilde{\gamma}\delta=0},
\quad
\bm{1}_{u=0} \bm{1}_{v\ne 0}
= \one_{\delta^\dagger\tilde{\gamma}\delta=0} \one_{\tilde{\gamma}\delta\ne 0}, \\
\bm{1}_{u\ne 0} \bm{1}_{v\in \FF_q u}
= \one_{\delta^\dagger\tilde{\gamma}\delta\ne 0} \one_{\tilde{\gamma}\delta\in \FF_q^nM_2(\FF_q)},
\quad
\bm{1}_{u\ne 0} \bm{1}_{v\notin \FF_q u}
= \one_{\delta^\dagger\tilde{\gamma}\delta\ne 0} \one_{\tilde{\gamma}\delta\notin \FF_q^nM_2(\FF_q)}.
\end{split}
\end{equation*}
If $u=v=0$, then
\begin{equation*}
s_i = \trd(\tilde{\gamma}_i),
\quad
\#\{a,b: \sum_{1\le i\le n} s_ia_i+t_ib_i = 0\}
= q^{2n-1} + (q^{2n}-q^{2n-1})\one_{\tilde{\gamma}=0}.
\end{equation*}
If $u=0$, then
\begin{equation*}
\begin{split}
\#\{a: \sum_{1\le i\le n} (s_i-v_i)a_i = 0\} = q^{n-1} + (q^n-q^{n-1}) \one_{s-v=0}, \\
\#\{a,\lambda: \sum_{1\le i\le n} (s_i-v_i)a_i+t_iv_i\lambda = 0\}
= q^n + (q^{n+1}-q^n) \one_{s-v=\sum_{1\le i\le n} t_iv_i=0}, \\
(s_i-v_i)^2 = \trd(\tilde{\gamma}_i)^2 - 4\nrd(\tilde{\gamma}_i),
\quad \one_{s-v=0} = \one_{(\tilde{\gamma}-\tilde{\gamma}^\dagger)\delta=0}.
\end{split}
\end{equation*}
If $u=s-v=0$,
then $$s_i=v_i=\frac{\trd(\tilde{\gamma}_i)}{2}\one_{2\nmid q} + \nrd(\tilde{\gamma}_i)^{1/2}\one_{2\mid q},
\quad \one_{\sum_{1\le i\le n} t_iv_i=0} = \one_{\sum_{1\le i\le n} (\tilde{\gamma}_i-v_i)v_i = 0}.$$
If $u\ne 0$ and $v\in \FF_q u$,
then the map $a_i\mapsto a_i-s_i\lambda$ defines an isomorphism
\begin{equation}
\label{future-decisive-quadric-iso}
\{a,\lambda: \sum_{1\le i\le n} a_i^2 = \sum_{1\le i\le n} (s_i-v_i)a_i\lambda+t_iu_i\lambda^2\}
\cong \{a,\lambda: \sum_{1\le i\le n} a_i^2+\trd(\tilde{\gamma}_i)a_i\lambda+\nrd(\tilde{\gamma}_i)\lambda^2 = 0\}.
\end{equation}

\item Generally, by \eqref{EQN:partial-Cartan-symmetry}, we may compute $S_3(\delta,\gamma)$ intrinsically by combining (2) and (3).

\item $q^{4n} I_0(\delta,0)
= S_2(\delta)
= q^{3n-2}(q^n-1) + q^{2n} X_2(0)
= q^{4n-2} + q^{3n-1} - q^{3n-2} + O(q^{5n/2})$.

\item If $\tilde{\gamma}\delta\notin \FF_q^nM_2(\FF_q)$, then $I_0(\delta,\gamma) = 0$.
If $\tilde{\gamma}\delta\in \FF_q^nM_2(\FF_q)$ and $\delta^\dagger\tilde{\gamma}\delta\ne 0$,
then
\begin{equation*}
q^{4n} I_0(\delta,\gamma) \ll q^{2n-2+\lfloor{(n+1)/2}\rfloor}
\end{equation*}
unless $2\nmid q$ and $\sum_{1\le i\le n} \nrd(\tilde{\gamma}_i) = \frac14\sum_{1\le i\le n} \trd(\tilde{\gamma}_i)^2$,
in which case
\begin{equation*}
q^{4n} I_0(\delta,\gamma) \ll q^{2n-1+\lfloor{n/2}\rfloor}.
\end{equation*}
If $\delta^\dagger\tilde{\gamma}\delta=0$ and $\tilde{\gamma}\delta\ne 0$, then $q^{4n} I_0(\delta,\gamma)\ll q^{3n-1} \one_{(\tilde{\gamma}-\tilde{\gamma}^\dagger)\delta=0}$.
\end{enumerate}

\end{prop}

\begin{proof}
(1):
Note that $I_0(\delta,\gamma)$ depends only on $\delta$ and $\tilde{\gamma}$.
By abuse of notation, let $I_0(\delta,\tilde{\gamma}) \defeq I_0(\delta,\gamma)$.
Furthermore, $I_0(\delta,\tilde{\gamma}) = I_0(\delta,a\tilde{\gamma})$ for all $a\in \FF_q^\times$ by the homogeneity of $P$.
Therefore, on summing over $a\in \FF_q^\times$, we find that 
\begin{align*}
    (q-1) I_0(\delta,\gamma)
&= \left(\sum_{a\in \FF_q} I_0(\delta,a\tilde{\gamma})\right) - I_0(\delta,0)\\
&= q\int_{\OO_F}\int_{\OO_D^n}\one_{\varpi\OO_D}(\delta^\dagger P(Y)) \psi(a\trd(\gamma\cdot Y)/\nrd(
\delta
))dY da-q^{-4n}S_2(\delta)\\
&=q\int_{\OO_D^n}\one_{\varpi\OO_D}(\delta^\dagger P(Y)) \one_{\varpi\OO_F}(\trd(\gamma\cdot Y))dY -q^{-4n}S_2(\delta)\\
&=q^{-4n}(qS_3(\delta,\gamma)-S_2(\delta)).
\end{align*}

(2):
In this case, $\delta^\dagger = \begin{psmatrix}1 & 0 \\ 0 & \varpi\end{psmatrix}$.
Write $Y = \begin{psmatrix}a & b \\ c & d\end{psmatrix}$.
Then $S_3$ is the number of solutions to
\begin{equation*}
\sum_{1\le i\le n} a_i^2+b_ic_i = \sum_{1\le i\le n} a_ib_i+b_id_i
= \sum_{1\le i\le n} s_ia_i+t_ib_i+u_ic_i+v_id_i = 0
\end{equation*}
in $\FF_q^{4n}$.
We decompose the point count into several disjoint cases, as follows.

Contribution from $b=0$:
If $(u,v)\ne 0$, then there are $X_2(0)$ choices for $a$ satisfying the first equation,
and $q^{2n-1}$ choices for $(c,d)$.
If $u=v=0$, then instead there are $X_2(s)$ consistent choices for $a$,
and $q^{2n}$ choices for $(c,d)$.
In total, $q^{2n} X_2(s) \bm{1}_{u=v=0} + q^{2n-1} X_2(0) \bm{1}_{(u,v)\ne 0}$.

Contribution from $b,u$ linearly independent:
$\bm{1}_{u\ne 0} q^n (q^n - q) q^{n-1} q^{n-2}$ (choose $a, b, d$ satisfying the second equation, then choose $c$ satisfying the remaining two equations).


Contribution from $b,u$ linearly dependent
with $b = \lambda u \ne 0$:
For each consistent $a,\lambda,d$, we get $q^{n-1}$ choices for $c$,
for a total of $\bm{1}_{u\ne 0} q^{n-1} \#\{a,\lambda,d: \lambda\ne 0,\; \sum_{1\le i\le n} a_i^2 = \lambda\sum_{1\le i\le n} s_ia_i+t_i\lambda u_i+v_id_i,\; \sum_{1\le i\le n} a_iu_i+u_id_i = 0\}$.
If $v\notin \FF_q u$ the latter count is $q^n (q-1) q^{n-2}$ (by choosing $a,\lambda$ first and then choosing $d$),
and if $v = \kappa u$ with $\kappa\in \FF_q$, then it is instead $q^{n-1} \#\{a,\lambda: \lambda\ne 0,\; \sum_{1\le i\le n} a_i^2 = \lambda\sum_{1\le i\le n} s_ia_i+t_i\lambda u_i-\kappa a_iu_i\}$ (since for each consistent $a,\lambda$ there are $q^{n-1}$ choices for $d$).
The last point count can be written as $\#\{a,\lambda: \sum_{1\le i\le n} a_i^2 = \lambda\sum_{1\le i\le n} s_ia_i+t_i\lambda u_i-v_ia_i\} - \#\{a: \sum_{1\le i\le n} a_i^2 = 0\}$.


Contribution from $u=0$, $b\ne 0$:
For each $a,b,d$, we get $q^{n-1}$ choices for $c$,
for a total of $\bm{1}_{u=0} q^{n-1} \#\{a,b,d: b\ne 0,\; \sum_{1\le i\le n} a_ib_i+b_id_i = 0,\; \sum_{1\le i\le n} s_ia_i+t_ib_i+v_id_i = 0\}$.
We break the count over $a,b,d$ into pieces, as follows.
\begin{itemize}
\item Contribution from $b,v$ linearly independent:
$\bm{1}_{v\ne 0} q^n (q^n-q) q^{n-2}$ (freely choose $a,b$, then choose $d$ satisfying the two independent linear conditions).

\item Contribution from $b,v$ linearly dependent
with $b = \lambda v \ne 0$:
For each consistent $a,\lambda$, we get $q^{n-1}$ choices for $d$,
for a total of $\bm{1}_{v\ne 0} q^{n-1} \#\{a,\lambda: \lambda\ne 0,\; \sum_{1\le i\le n} s_ia_i+t_i\lambda v_i = \sum_{1\le i\le n} a_iv_i\}$.

\item Contribution from $v=0$, $b\ne 0$:
There are $q^{n-1}$ choices for $d$ given $a,b$ satisfying the second equation;
this leads to a total of $\bm{1}_{v=0} q^{n-1} \#\{a,b: b\ne 0,\; \sum_{1\le i\le n} s_ia_i+t_ib_i = 0\}$.
\end{itemize}




(3):
This is routine.
For \eqref{future-decisive-quadric-iso}, we use the identities $a_i^2-(s_i-v_i)a_i\lambda = (a_i-s_i\lambda)(a_i+v_i\lambda) + s_iv_i\lambda^2$, $s_i+v_i=\trd(\tilde{\gamma}_i)$, and $\nrd(\tilde{\gamma}_i)=s_iv_i-t_iu_i$.

(4):
This is clear from \eqref{EQN:partial-Cartan-symmetry}, because (3) lets us rewrite the formula in (2) in terms of invariants of the action
$(\delta,\gamma,\alpha,\beta) \longmapsto (\alpha \delta \beta,\alpha \gamma \alpha^{-1})$.

(5):
This follows from (1) and (4).

(6):
Use (1), (4), the Weil conjectures for quadrics, and the formula for $S_2(\delta)$ in (5).
\end{proof}

\section{Abelian main terms}
\label{SEC:abelian}

In this appendix, we briefly discuss some special cases admitting an elementary abelian approach to estimating the $\gamma=0$ contribution in \cite{Main}, i.e.~
\begin{equation*}
\Sigma_0(X) :=
c_{\Phi,X}X^{4(n-1)}
\sum_{\delta \in D^\times}
\left(I_{1}\left(\frac{\delta}{X},0\right) -  I_0\left(\frac{\delta}{X},0\right)\right).
\end{equation*}
As in \cite[Theorem 1.2]{Main}, we make the following simplifying assumptions:
Suppose $F=\QQ$ and $D$ is ramified at $\infty,2$.
Suppose functions $f=\otimes_v f_v, \Phi=\otimes_v \Phi_v$ satisfy  
\begin{enumerate}
\item $f_p=\one_{\OO_{D_p}^n}$ and $\Phi_p=\one_{\OO_{D_p}^2}$ at all finite places $p$.
\item $\Phi_\infty(x,y) = \Phi_1(\nrd(x)) \Phi_2(\nrd(y))$, for some $\Phi_1,\Phi_2\in \mathcal{S}(\RR)$.
\end{enumerate}

We write $0\neq\delta\in \OO_D$ in the form $\delta = J\eta$ where $\eta\in \OO_D$ is primitive and $J\ge 1$. Note that $\one_{\OO_{\ZZ_p}}\circ \nrd=\one_{\OO_{D_p}}$ at ramified $p$. Then by the Cartan decomposition and \eqref{EQN:partial-Cartan-symmetry} at unramified $p$, we find that
$\prod_{p} I_{0,p}(\delta,0)
= A_0(J,\nrd(\eta))$ for some function $A_0:\ZZ_{\ge 1}^2 \to \CC$.
Also, we have
\begin{equation*}
I_{0,\infty}\left(\frac{\delta}{X},0\right)
= \int_{D_v^n} f(t) \Phi(\delta/X, \delta^{-1}XP(t)) \, dt
= A_{\infty,0}(\nrd(\delta/X)),
\end{equation*}
for some $A_{\infty,0}\in \mathcal{S}(\RR)$. Since $I_{0,p}(\delta,0) = I_{1,p}(\delta^\dagger,0)$ at all finite places $p,$ 
we get
\begin{equation*}
\Sigma_0(X)
= c_{\Phi,X} X^{4n-4} \sum_{J\ge 1}
\sum_{\substack{0\neq \eta\in \OO_D\\ \textnormal{$\eta$ primitive}}}
A_0(J,\nrd(\eta))
\cdot (A_{\infty,1}-A_{\infty,0})(\nrd(\tfrac{J\eta}{X})).
\end{equation*}

Suppose $\OO_D=\mathbb{H}$.
In this case, $\OO_D$ is a principal ideal domain,
so by \cite[Proposition~26.3.9 and Lemma~26.4.1]{voight2021quaternion},
we have a particularly simple formula
\begin{equation*}
r_D(m) \defeq
\frac{\#\{\textnormal{primitive }\eta\in \OO_D\cap D^\times: \nrd(\eta) = m\}}
{\#\OO_D^\times}
= \one_{v_2(m) \le 1} \prod_{2\nmid p\mid m} (p^{v_p(m)} + p^{v_p(m)-1})
\end{equation*}
with no contribution from cusp forms.
Therefore,
\begin{equation*}
\frac{\Sigma_0(X)}{\#\OO_D^\times}
= c_{\Phi,X} X^{4n-4} \sum_{J\ge 1} \sum_{m\ge 1} r_D(m) A_0(J,m)
\cdot (A_{\infty,1}-A_{\infty,0})(J^2m/X^2).
\end{equation*}
Since $r_D(m) A_0(J,m)$ is multiplicative,
Mellin inversion should now reduce us to understanding the double Dirichlet series
$\sum_{J\ge 1} \sum_{m\ge 1} \frac{r_D(m) A_0(J,m)}{J^r m^s}$
and the Mellin transforms of $A_{\infty,j}$.

For general $D$, this approach fails, because $r_D(m)$ likely need not be Eulerian,
either classically or adelically.
Similarly, the approach fails for general $f$, even if $\OO_D=\mathbb{H}$.

\bibliography{refs}{}
\bibliographystyle{alpha}